\documentclass[10pt,a4paper]{article}

\usepackage{amsmath}
\usepackage{graphicx,latexsym,euscript,makeidx,color}
\usepackage{amsmath,amsfonts,amssymb,amsthm,mathrsfs,dsfont}
\usepackage{enumerate}
\usepackage{multirow}

\usepackage{amssymb}
\usepackage{amsmath}
\usepackage{amsfonts}
\usepackage{graphicx}
\usepackage{graphicx}          
\usepackage{color}

\usepackage[colorlinks,linkcolor=blue,anchorcolor=green,citecolor=red]{hyperref}

\usepackage{geometry}
\def\ds{\displaystyle}

\def\a{\alpha}           
            
\def\z{\zeta}

           \def\F{\Phi}     
             
          \def\f{\varphi}

\font\tenbb=msbm10 \font\sevenbb=msbm7 \font\fivebb=msbm5
\newfam\bbfam
\scriptscriptfont\bbfam=\fivebb \textfont\bbfam=\tenbb
\scriptfont\bbfam=\sevenbb

\newtheorem{theorem}{\indent Theorem}[section]
\newtheorem{definition}[theorem]{\indent Definition}
\newtheorem{proposition}[theorem]{\indent Proposition}

\newtheorem{remark}[theorem]{\indent Remark}
\newtheorem{example}[theorem]{\indent Example}

\makeatletter
   
   \@addtoreset{equation}{section}
\makeatother

\allowdisplaybreaks[4]

\begin{document}

\title{\bf A Nash-Type Fictitious Game Framework to Time-Inconsistent Stochastic Control Problems\thanks{This work is supported in part by the National Natural Science Foundation of China (61773222, 11871369, 61973172, 62173191).}
}
\author{Yuan-Hua Ni\thanks{College of Artificial Intelligence, Nankai University, Tianjin 300350, P.R. China. Email: {\tt yhni@nankai.edu.cn}. }~~~~Binbin Si\thanks{College of Artificial Intelligence, Nankai University, Tianjin 300350, P.R. China. Email: {\tt binbinsi96@163.com}.} ~~~~Xinzhen Zhang\thanks{School of Mathematics, Tianjin University, Tianjin 300352, P.R. China. Email: {\tt xzzhang@tju.edu.cn}.}}
\maketitle

{\bf Abstract:} In this paper, a Nash-type fictitious game framework is introduced to handle time-inconsistent linear-quadratic optimal control problems.
The Nash-type game in this framework is called fictitious as it is between the decision maker (called real player) and an auxiliary control variable (called fictitious player) with the real player and fictitious player  looking for time-consistent policy and precomitted optimal policy, respectively.
Namely, the fictitious-game framework is actually an auxiliary-variable-based mechanism where the fictitious player is our particular design.
Noting that the real player's cost functional is revised in accordance with that of fictitious player, the equilibrium policy of real player is called an open-loop self-coordination control of original linear-quadratic problem.
As a generalization, a time-inconsistent nonzero-sum stochastic linear-quadratic dynamic game is investigated, where one player is to look for precommitted optimal policy and the other player is to search time-consistent policy.
Necessary and sufficient conditions are presented to ensure the existence of open-loop equilibrium of the nonzero-sum game, which resort to a set of Riccati-like equations and linear equations.
By applying the developed theory of nonzero-sum game, open-loop self-coordination control of the linear-quadratic optimal control is fully characterized, and multi-period mean-variance portfolio selection is also investigated.
Finally, numerical simulations are presented, which show the efficiency of the proposed fictitious-game framework.


\section{Introduction}

Dynamic programming is a fundamental and powerful approach to solving optimal control problems; the basic idea is to consider a family of problems with different initial times and states, and to establish relationships among these problems. Bellman's principle of optimality is the core of this approach which states using Bellman's words \cite{Bellman} as ``An optimal policy has the property that whatever the initial state and initial decision are, the remaining decisions must constitute an optimal policy with regard to the state resulting from the first decision." This property is termed as the time consistency of optimal control.
However, in many practical situations, the time-consistency fails quite often and the corresponding optimal control problems are time-inconsistent.

Actually, there are several factors that ruin the time consistency of optimal control. The first one is the nonlinear terms of conditional expectations of state/control that appear in the objective functional. As there is no nonlinear version of the tower property of conditional expectation, the controller at different time instants is facing with different objectives, which are not consistent with the global objective. Therefore, the time inconsistency comes from the conflicts between global optimal control on the lifetime horizon and local optimal control on the tail time horizon.
Such a kind of time-inconsistent problems are classified as mean-field optimal controls, which have gained considerable attention during the last few years \cite{Ni-Zhang-Krstic2017,Yong-2013}.
Another factor is the non-exponential discounting in objective functional, which does not
possess the property of group or separability any longer. Though exponential discounting is of great importance to model people's time preference \cite{Samuelson1937}, empirical researches over the last half century have documented the inadequacy of constant discount rate. Among others, hyperbolic discounting is a known anomaly and is often used to describe the case with a declining discount rate \cite{Time-discounting}.
Further, it is hypothesized \cite{Strotz} that people are born with tendency to overvalue current consumption and that more discounting occurs between the present and the near future than between periods in the more distant future, namely, the discounting function in the  objective functional is non-exponential.
The aforementioned two factors that ruin the time consistency reflect people's risk preferences and time preference in some nontraditional ways, both of which are of the phenomenon of ``changing tastes" in intertemporal choices \cite{Auer}, namely, today's preference conflicts with tomorrow's preference. 
%

\subsection{Existing methodologies on handling the time inconsistency}

To handle the time inconsistency, there are several different approaches in existing literature and a rule of selecting the preferred solution is called  as a choice mechanism \cite{Auer}.
The first one is the precommitment choice for which the initial policy is implemented on the lifetime horizon. This approach neglects the time inconsistency, and the optimal policy is optimal only when viewed at the initial time.
The second mechanism is naive choice or myopic choice: at each time instant a naive
agent embarks on the option that currently seems best, namely, this agent sticks to the local objective and completely ignores the global interest. However, the naive policy makes no sense of optimality, and simple example \cite{Auer} shows that it might be the worst one of all the policies viewed from the initial time instant.

\subsubsection{Strotz's time-consistent solution}

Another mechanism is sophisticated/time-consistent choice proposed by Strotz \cite{Strotz}. In the viewpoint of Strotz,
the decision maker at different time instants is regarded as different selves, and the time inconsistency suggests a conflict between different these selves.
At any time instant the current self takes account of future selves' decisions, and the equilibrium of this intertemporal game is called a sophisticated policy, or a time-consistent policy.
Inspired by Strotz's idea, hundreds of works have sought to tackle practical problems in economics and finance;
see, for example,
\cite{{Cui-2017},{Cui-2012},Ekeland,Ekland-2,Goldman,Krusell,Laibson,{Ni-Li-Zhang-Krstic-2020},Palacios} and
the references therein.
Moreover, accompanying the appearance of time-inconsistent mean-field optimal control, recent years have witnessed the rapid  progresses of extending Strotz's idea in the theoretical control
community \cite{Bjork,Bjork2017,Hu-jin-Zhou,Hu-jin-zhou-2,{Ni-Li-Zhang-Krstic-2019},{Ni-Zhang-Krstic2017},Pun,Wang-haiyang,Wang-Tianxiao1, Wang-Tianxiao2,Wang-Tianxiao3,Wei-Yu-Yong2017,Yong-1,Yong-0,Yong-2013}; in particular, the open-loop time-consistent control, feedback time-consistent strategy and mixed time-consistent solution are elaborately studied.
It should be noted that Strotz's solution is essentially a closed-loop time-consistent strategy.
So far, the sophisticated policy and precommitted policy are two extreme solutions; namely, the sophisticated policy recovers the time consistency and ignores the global optimality, while the precomitted solution does not care about the time consistency and just pays attention to the global optimality on the lifetime horizon. 

\subsubsection{Self-control of Thaler and Shefrin}


Different from Strotz's formulation and to balance the global optimality and time consistency, the economists Thaler and Shefrin \cite{Thaler} introduce a two-dimensional self-control model: the individual at any instant in time is assumed to be both a farsighted planner and a myopic doer. This division into conflicting subselves is how psychologists think about self-control, and the notion of self-control is paradoxical without it.
The doer at each moment in time exists only for one period and is completely selfish, or myopic \cite{Thaler}, namely, the objective of each doer is independent of past and future variables that are concerned. On the contrary, the planner is concerned with the lifetime objective, which is derived from the objectives of all the doers.
Interestingly, pointed out by \cite{Thaler} and due to the myopic nature of the doers, the conflict between the planner and doers is fundamentally similar to the agency relationship between the employer/principal and employees/agents of a firm.
In fact, this two-self model to understand the savings behavior of individuals and households is one of four Thaler's contributions in behavioral economics to win 2017 Nobel Prize \cite{Thaler2017}.
For the recent progresses on self-control, we are referred to \cite{Benabou,Cui-2017-2,Cui-Li-Shi-2017,Fudenberg2006,Fudenberg2012,Gul,ODonoghue} and references therein.

Actually, the idea of two-self model can be traced back to the work of Adam Smith \cite{A-Smith} in 1759; and \cite{Thaler} is the first systematic and formal treatment of a two-self economic man, which integrates economics with psychology. A key feature of this planner-doer modelling is that the planner is also allowed to bear some influence on doers' behavior. For this, the doers are  given the discretion to either modify their preferences or alter the incentives (rewards, punishments, \emph{etc}). Specifically, through incorporating the costly control of a ``preference modification parameter" (selected by the planner) into the doers' utilities, the behaviors of the planner and doers can be mutually influenced, and the planner's utility is simultaneously modified.
By finding the equilibrium of this intrapersonal game, a balance between the lifetime objective and myopic objectives is achieved. Noting that the planner does not actually consume, the policy of this game selected by the doers is the one that is executed by the individual.

\subsubsection{Two-tier game framework of Cui, Li and Shi}

Note that the doers' local objectives of standard self-control schemes are fully myopic, namely, there are no conflicts among different local objectives. To meet the time-inconsistent stochastic decision problems with conflicting non-myopic local objectives and non-expectation operators, Cui, Li and Shi propose a two-tier planner-doer game framework \cite{Cui-2017-2} to reconcile the global and local interests, where a sequential game among the doers is involved that significantly extends existing self-control schemes \cite{Fudenberg2006,Fudenberg2012,Gul,ODonoghue,Thaler}. 
Through commitment by punishment, the proposed mechanism of \cite{Cui-2017-2}
revises the original preference of each individual doer by adding a penalty term, while the expected
total penalty in turn modifies the planner's preference. Then, both a sequential game among the doers in the lower tier and a leader-follower game between the planner and doers in the upper tier are constructed.
Given any planned policy, the best-response policies of the doers form a Nash equilibrium of the low tier game
and are time-consistent on the lifetime horizon. Therefore, the two-tier game is indeed a game between precommitted policy and time-consistent policy, and the doers' equilibrium time-consistent policy of this game is called a self-coordination policy of the original time-inconsistent decision problem.

Furthermore, the proposed self-coordination mechanism is applied to the dynamic mean-variance portfolio selection, and an explicit self-coordination policy is obtained together with a detailed sensitivity analysis \cite{Cui-2017-2}; this enable investors to understand the trade-off between global and local interests and coordinate among various selves.
To the best of the authors and within the realm of dynamic games, the work \cite{Cui-2017-2} is the first to study the game between precommitted policy and time-consistent policy, which clearly merits further investigations. The followings are some of questions that might be studied.

\begin{itemize}

\item[i).] There exist several notions on time-consistent equilibrium solution in existing literature, such as the open-loop equilibrium control, closed-loop equilibrium strategy and mixed equivalent equations \cite{Hu-jin-Zhou,Ni-Li-Zhang-Krstic-2019,Yong-2013}. Note that the best-response time-consistent policy of the doers is closed-loop. So, it is desirable to study other types of best-response time-consistent policy.

\item[ii).] If we go beyond, the game between the planner and doers can be settled within more general framework of games between precommitted policy and time-consistent policy, whose study must enrich the game theory. In particular, it is desirable to study the Nash-type games, namely, the precommitted policy and time-consistent policy have equal status.

\item[iii).] Linear-quadratic (LQ, for short) optimal control is pioneered by Kalman in 1960s and is now a classical yet fundamental problem in control theory.  Though general time-inconsistent decision problems has been considered in \cite{Cui-2017-2, Cui-Li-Shi-2017}, the pretty structure of LQ problem has not been fully explored within the games between the precommitted policy and time-consistent policy.

\end{itemize}


\subsection{Nash-type fictitious game framework}

\subsubsection{The framework}

In this paper, a time-inconsistent LQ optimal control is studied via a method of Nash-type fictitious game. 
%
%
%
%
%
Specifically, consider the system
\begin{eqnarray}\label{system-1}
\left\{\begin{array}{l}
X_{k+1}=\big{(}{A}^0_{k}X_k+B^0_{k}v_k\big{)}+\sum_{i=1}^{p}\big{(}{C}^{0i}_{k}X_k+{D}^{0i}_{k}v_k\big{)}w^{i}_k, \\[1mm]
X_t=x,~~k\in  \mathbb{T}_t,~~t\in \mathbb{T},
\end{array}
\right.
\end{eqnarray}
where $\mathbb{T}=\{0, \dots, N-1\}$, $\mathbb{T}_t=\{t, \dots, N-1\}$, and $A^0_{k}, {C}^{0i}_{k}\in \mathbb{R}^{n\times
n}$, $B^0_{k}, {D}^{0i}_{k}\in
\mathbb{R}^{n\times m}$ are deterministic.
Letting $w_k=(w_k^1, \dots,w_k^{p})^T$,
the noise process $\{w_k, k\in
\mathbb{T}\}$  is assumed to be a vector-valued martingale difference sequence
defined on a probability space $(\Omega, {\mathbb{F}}, \mathbb{P})$ with
\begin{eqnarray}\label{w-moment-0-0}
\mathbb{E}_{k}[w_{k}]=0,~~\mathbb{E}_{k}[w_{k}w_k^T]=\Delta_k,~~k\in \mathbb{T},
\end{eqnarray}
where $\Delta_{k}=(\delta_{k}^{ij})_{p\times p}, k\in \mathbb{T}$, are assumed to be deterministic.
$\mathbb{E}_k[\,\cdot\,]$ in (\ref{w-moment-0-0}) is the conditional mathematical expectation
$\mathbb{E}[\,\cdot\,|\,\mathbb{F}_k]$ with respect to $\mathbb{F}_k=\sigma\{w_l, l=0, 1,\dots,k-1\}$, and $\mathbb{F}_{0}$ is understood as
$\{\emptyset, \Omega\}$.
Let $l^2_{\mathbb{F}}(t; \mathbb{R}^n)$ and $l^2_{\mathbb{F}}(\mathbb{T}_{t}; \mathbb{R}^m)$ be defined as
\begin{eqnarray}\label{L^2-t-LQ}
&&l^2_\mathbb{F}(t; \mathbb{R}^{n})=\Big{\{}\zeta \in \mathbb{R}^{n}\,\big{|}\, \zeta\mbox{ is }\mathbb{F}_t\mbox{-measurable}, \mathbb{E}|\zeta|^2<\infty \Big{\}},
\\
\label{L^2-T-LQ}
&&l^2_\mathbb{F}(\mathbb{T}_{t}; \mathbb{R}^{m})=\left\{\nu=\{\nu_k, k\in \mathbb{T}_{t}\}\,\Big{|}\,
\nu_k\in \mathbb{R}^{m}\mbox{ is }\mathbb{F}_k\mbox{-measurable},
\mathbb{E}|\nu_k|^2<\infty, k\in\mathbb{T}_{t}\right\};
\end{eqnarray}
and $x$ of (\ref{system-1}) belongs to $l^2_{\mathbb{F}}(t; \mathbb{R}^n)$. 
%
%
%
%
%
The objective functional is
\begin{eqnarray}\label{cost-1}
&&\hspace{-3.5em}J(t,x;v)
 =\sum_{k=t}^{N-1}\mathbb{E}_t\Big{\{}X_k^T{Q}^0_{t,k}X_k +(\mathbb{E}_tX_k)^T
 \bar{{Q}}^0_{t,k}\mathbb{E}_tX_k+ v_k^T{R}^0_{t,k}v_k+(\mathbb{E}_tv_k)^T \bar{{R}}^0_{t,k}
 \mathbb{E}_tv_k \Big{\}}
\nonumber\\
&&\hspace{-3.5em}\hphantom{J(t,x;u)=}+\mathbb{E}_t\big{[}(X_N)^T{G}^0_tX_N\big{]}
+(\mathbb{E}_tX_N)^T \bar{{G}}^0_t\mathbb{E}_tX_{N}, 
\end{eqnarray}
where ${Q}^0_{t,k}, \bar{{Q}}^0_{t,k},{R}^0_{t,k},\bar{{R}}^0_{t,k}, k\in
\mathbb{T}_t$, ${G}^0_t, \bar{{G}}^0_t$ are deterministic symmetric matrices
of appropriate dimensions.
Then, the LQ problem is stated as follows.

\textbf{Problem (LQ).}  \emph{Letting $t\in \mathbb{T}$ and $x\in l^2_{\mathbb{F}}(t; \mathbb{R}^n)$, find a
${v}^*\in l^2_{\mathbb{F}}(\mathbb{T}_t; \mathbb{R}^m)$ such that
\begin{eqnarray}\label{Problem-LQ-precommitted}
J(t,x;{v}^*) = \inf_{u\in l^2_{\mathbb{F}}(\mathbb{T}_t; \mathbb{R}^m)}J(t,x;v).
\end{eqnarray}
}

Problem (LQ) is time-inconsistent as the objective functional (\ref{cost-1}) contains nonlinear terms of conditional expectation and the weighting matrices of (\ref{cost-1}) depend on the initial time.
$u^{*}$ of (\ref{Problem-LQ-precommitted}) is called an open-loop precommitted optimal control, or simply precommitted solution/policy, for the initial pair $(t,x)$, which totally adheres to the global interest on the lifetime horizon $\mathbb{T}_t$. Noting that $u^*$ neglects the time inconsistency, the following notion yet pays attention to the open-loop time-consistent solution of Problem (LQ). 

\begin{definition}\label{Definition-time-consistency}
A control $v\in l^2_{\mathbb{F}}(\mathbb{T}_t;\mathbb{R}^m)$ is called an open-loop time-consistent equilibrium control of Problem (LQ) for the initial pair $(t,x)$, if for any $k\in \mathbb{T}_t$ and any $\bar{v}_k\in
l^2_{\mathbb{F}}(k; \mathbb{R}^m)$,
\begin{eqnarray}
J\big{(}k, X_k; v|_{\mathbb{T}_k}\big{)}\leq J\big{(}k, X_k; (\bar{v}_k,v|_{\mathbb{T}_{k+1}})\big{)}.
\end{eqnarray}
Here, $v|_{\mathbb{T}_k}$ and
$v|_{\mathbb{T}_{k+1}}$  are the restrictions of $v$
on $\mathbb{T}_k$ and $\mathbb{T}_{k+1}$, respectively; and $X_k$ is given by
\begin{eqnarray}\label{system-1-real-player}
\left\{\begin{array}{l}
X_{k+1}=\big{(}{A}^0_{k}X_k+B^0_{k}v_k\big{)}+\sum_{i=1}^{p}\big{(}{C}^{0i}_{k}X_k+{D}^{0i}_{k}v_k\big{)}w^{i}_k, \\[1mm]
X_t=x,~~k\in  \mathbb{T}_t,~~t\in \mathbb{T}.
\end{array}
\right.
\end{eqnarray}

\end{definition}

Since Strotz's work, time-consistent solutions have gained much attention in the areas of economics, finance \emph{etc}.
It should be noted that time-consistent policy and precommitted optimal policy are two extreme and irreconcilable solutions for time-inconsistent optimal control. More precisely, time-consistent policy recovers the time consistency and ignores the global optimality, while the precomitted solution does not care about the time consistency and just pays attention to the global optimality on the lifetime horizon.
In this paper, we will extend the notion of time-consistent solution in order to strike a balance between the time consistency and global optimality. The main idea is to introduce an auxiliary fictitious player that plays games with the decision maker $v$ in Definition \ref{Definition-time-consistency} (this $v$ is called a real player throughout this paper).

The Nash-type fictitious game framework of this paper is divided into the following three steps.

\emph{Step 1}. Introduce an auxiliary control variable $u\in l^2_{\mathbb{F}}(\mathbb{T}_t;\mathbb{R}^m)$, which is called a fictitious player throughout this paper.
The cost functional of $u$ is $J(t, x; u)$ with the internal state
\begin{eqnarray}\label{widehat-X}
\left\{
\begin{array}{l}
\widehat{X}_{k+1}=\big{(}A^0_{k}\widehat{X}_k+B^0_{k}u_k\big{)} +\sum_{i=1}^p\big{(}{C}^{0i}_{k}\widehat{X}_k+{D}^{0i}_{k}u_k\big{)}w^i_k,\\[1mm]
\widehat{X}_t=x,~~~ k\in \mathbb{T}_t,
\end{array}
\right.
\end{eqnarray}
namely,
\begin{eqnarray*}\label{cost-1-u}
&&\hspace{-3.5em}J(t, x; u)
 =\sum_{k=t}^{N-1}\mathbb{E}_t\Big{\{}\widehat{X}_k^T{Q}^0_{t,k}\widehat{X}_k +(\mathbb{E}_t\widehat{X}_k)^T
 \bar{{Q}}^0_{t,k}\mathbb{E}_t\widehat{X}_k+ u_k^T{R}^0_{t,k}u_k+(\mathbb{E}_tu_k)^T \bar{{R}}^0_{t,k}
 \mathbb{E}_tu_k \Big{\}}
\nonumber\\
&&\hspace{-3.5em}\hphantom{J(t,x;u)=}+\mathbb{E}_t\big{[}(\widehat{X}_N)^T{G}^0_t\widehat{X}_N\big{]}
+(\mathbb{E}_t\widehat{X}_N)^T \bar{{G}}^0_t\mathbb{E}_t\widehat{X}_{N}. 
\end{eqnarray*}
Different from the real player ($v$ of Definition \ref{Definition-time-consistency}), the fictitious player $u$ would like to find a precommitted optimal solution.

\emph{Step 2}. Modify the cost functionals of $v$ and $u$ to
\begin{eqnarray}
&&\hspace{-2em}\label{bar-J}\bar{J}(k, \bar{x}; u|_{\mathbb{T}_k}, v|_{\mathbb{T}_k}) =J(k,\bar{x};v|_{\mathbb{T}_k})+\mu_k(u_k-v_k)^T\Psi_k(u_k-v_k),~~~k\in \mathbb{T}_t,\\
&&\hspace{-2em}\label{widehat-J}\widehat{J}(t, x; u, v)=J(t,x; u)+\sum_{k=t}^{N-1}\mathbb{E}_t\big{[}\mu_k(u_k-v_k)^T \Psi_k(u_k-v_k)\big{]}
\end{eqnarray}
with $\Psi_k\in \mathbb{R}^{m\times m}$ being symmetric and $\mu_k\in \mathbb{R}, k\in \mathbb{T}_t$. Here, the inner states in $\bar{J}(k,\bar{x};u|_{\mathbb{T}_k}, v|_{\mathbb{T}_{k}})$ and $\widehat{J}(t,x;u, v)$ are given, respectively, by (\ref{system-1-real-player}) and (\ref{widehat-X}).
%
%
%
%
%
%

%
%

\emph{Step 3}. Solve the fictitious game:
\begin{itemize}
\item[]\textbf{Problem (LQ)$_g$}. 
Find $({u}^*,v^*)\in l^2_{{\mathbb{F}}}(\mathbb{T}_t; \mathbb{R}^{m})\times l^2_{{\mathbb{F}}}(\mathbb{T}_t; \mathbb{R}^{m})$ such that
\begin{eqnarray}\label{Problem-LQ-0}
&&\hspace{-4em}\widehat{J}(t,x; u^*, v^*)\leq \widehat{J}(t,x;u, v^*),~~~\forall u\in l^2_{{\mathbb{F}}}(\mathbb{T}_t, \mathbb{R}^{m}), \\
&&\hspace{-4em}\label{Problem-LQ-0-2}\bar{J}(k,{X}^{*}_k; u^*|_{\mathbb{T}_k}, v^*|_{\mathbb{T}_k})\leq \bar{J}(k,{X}^{*}_k;u^*|_{\mathbb{T}_k}, (v_k, v^*|_{\mathbb{T}_{k+1}})),~~~\forall k\in \mathbb{T}_t,~ \forall v_k\in l^2_{{\mathbb{F}}}(k, \mathbb{R}^{m})
\end{eqnarray}
hold with ${X}^*_k$ computed via
\begin{eqnarray}\label{System-equilibrium-0}
\left\{
\begin{array}{l}
{X}^{*}_{k+1}=\big{(}A^0_{k}{X}^{*}_k+B^0_{k}v^*_k\big{)} +\sum_{i=1}^{p}\big{(}{C}^{0i}_{k}{X}^{*}_k+{D}^{0i}_{k}v^*_k\big{)}{w}^{i}_k,\\[1mm]
{X}^{*}_t=x,~~~ k\in \mathbb{T}_t.
\end{array}
\right.
\end{eqnarray}
%
%

\end{itemize}
Here, $v^*$ is called an open-loop self-coordination control of Problem (LQ) for the initial pair $(t,x)$ and $\{\mu_k\Psi_k, k\in \mathbb{T}_t\}$.

\begin{remark}
$\{\Psi_k, k\in \mathbb{T}_t\}$ and $\{\mu_k, k\in \mathbb{T}_t\}$ are called the punishment direction and punishment intensity, respectively.
Note that modifying the objective functionals (\ref{bar-J}) (\ref{widehat-J}) of Step 2 is similar to that of \cite{Cui-2017-2}. Actually, it is \cite{Cui-2017-2} that motivates the study of this paper, and the term ``self-coordination policy" is introduced firstly by \cite{Cui-2017-2} that is of closed loop indeed. As the formulation of this paper is somewhat similar to that of \cite{Cui-2017-2}, the ``self-coordination" is borrowed here and the equilibrium policy of real player is called an open-loop self-coordination control of Problem (LQ).
Furthermore, Problem (LQ)$_g$ is called a fictitious game as it is a game between a real player and a fictitious player.

%
%


\end{remark}

\subsubsection{Motivation}

The Nash-type  fictitious game framework can be viewed as an auxiliary-variable method.
The $v$ of (\ref{system-1}) is the real controller and its policy is the one that is actually performed,
and the fictitious player $u$ is an auxiliary control variable with system equation (\ref{widehat-X}). On the fictitious game of Problem (LQ)$_g$, the real player's policy is obtained that is called an open-loop self-coordination control of Problem (LQ).
Namely, all the derivations from (\ref{widehat-X}) to (\ref{System-equilibrium-0}) are our particular design in order to obtain the open-loop self-coordination control.
Furthermore, the method of fictitious game differs from the two-tier planner-doer game \cite{Cui-2017-2} in the following two points. Firstly, the two-tier game framework is a self-control scheme, where the decision maker at any instant is assumed to have conflicting subselves---planner and doer. On the contrary, there is only one real controller in the fictitious game framework, and the fictitious player is an auxiliary variable. Secondly, the game of \cite{Cui-2017-2} is of lead-follower type, where the planner is a leader and the doers are the followers. In contrast, the game between real player and fictitious player of this paper is Nash-type, namely, the real player and fictitious player have equal status.

Auxiliary-variable method is somewhat usual in controller design \cite{auxiliary-1,auxiliary-2}; in this case, some internal auxiliary variables are introduced by which the concerned controllers can be constructed. In the fictitious game framework, the fictitious player is an internal auxiliary control variable, and through the fictitious game the real player's policy is obtained. Furthermore, the best-response policies of real player and fictitious player are time-consistent and globally optimal, respectively. In other words, the fictitious game is indeed a game between a real time-consistent policy and a fictitious precommitted policy, which shares some similarity with that of \cite{Cui-2017-2}. To the best of the authors, the work \cite{Cui-2017-2} is the first to study the game between time-consistent policy and precommitted policy, where the conflicting subselves have hierarchical status. Under the idea of introducing auxiliary control variable, it is possible to conduct a fictitious Nash-type game between time-consistent policy and precommitted policy.
Similarly to standard Nash equilibrium, the Nash-type equilibrium of fictitious game is also  non-cooperative, namely, if any one player stays at its equilibrium policy, the other player's equilibrium policy is optimal in the sense of (\ref{Problem-LQ-0}) or (\ref{Problem-LQ-0-2}).

%
%
%

Furthermore, on the status of precommitted policy and time-consistent policy in the fictitious game, let us go back to the pioneer work \cite{Strotz} too.
%
%
Though consistent planning is proposed firstly by \cite{Strotz},
%
%
%
it seems that no prejudice on precommitted policy is found in \cite{Strotz}. If the conflicts between the global interest and local interests are recognized, it ``may be solved either by (a) a strategy of precommitment, or (b) a strategy of consistent planning"; this is noted in the Summary of \cite{Strotz}.
More interestingly, in the very recent year, Caliendo and Findley \cite{Caliendo2019} present some positive results about that the precommitted policy acts better than time-consistent policy: ``In some prominent, well-studied examples with intertemporal tradeoffs (like the choice between investing in a project now or later, doing an unpleasant task now or procrastinating it until later, and eating a cake), we find that the commitment allocation can multiself Pareto dominate the non-cooperative equilibrium allocation if the number of time-dated selves exceeds a low threshold." Here, the commitment allocation and non-cooperative equilibrium allocation are the precommitted policy and time-consistent policy, respectively.
Besides, 
commitment mechanism is widely accepted in the fields such as decision science, economics and finance. To mention a few, see, for example, \cite{Ariely2002,Casari2009,Kivetz2002,Kurth2012} and references therein.
Namely, though nowadays consistent planning has gained much attention in the control community, 
the above mentioned sample of papers provides many positive evidences of studying precommittmed policy, which ought not to be ignored for its practical values in some situations. Therefore, for the general study, it is reasonable to treat the precommitted policy and time-consistent policy equally; this is the case of \cite{Strotz}. Hence,
to balance the global optimality and time consistency, an alternative way to that of \cite{Cui-2017-2} is conducting a Nash-type game between time-consistent policy and precommitment policy. This is another motivation of the fictitious game of this paper.

Additionally, concerned with standard nonzero sum games, Nash solution and Stackelberg solution are two standard noncooperative equilibrium solutions in the sense that no player can achieve an improvement if she attempts to deviate from her strategies. Nash solution ensures simultaneously that at the same time each player will not benefit from changing their strategy, and Stackelberg solution is in a sequential manner. Though it is established in \cite{Simaan-2} that the leader in the Stackelberg solution is at least as good and possibly lower cost than in the Nash solution, yet  nothing can be said about the follower who may or may not do better than the Nash solution \cite{Simaan}. Hence, for the self-coordination policy of \cite{Cui-2017-2} and open-loop self-coordination control of this paper, generally it is hard to say that one acts better than the other one. Interestingly, according to the examples of Section \ref{example}, the scheme of \cite{Cui-2017-2} looks for self-coordination policy between open-loop precommitted optimal control and open-loop time-inconsistent equilibrium control, while to some extent this paper goes beyond open-loop precommitted optimal control and open-loop time-consistent equilibrium control. Therefore, our formulation adds a new dimension to handle time-inconsistent optimal control problems.

\subsubsection{Generalization}

Letting $X^a_k=[\widehat{X}_k^T~~{X}_k^T]^T$, we have
\begin{eqnarray}\label{X-a}
\left\{
\begin{array}{l}
{X}^a_{k+1}=\left(
\begin{array}{cc}
A^0_k&0\\0&A^0_k
\end{array}
\right){X}^a_k+\left(
\begin{array}{c}
B^0_k\\0
\end{array}
\right)u_k+\left(
\begin{array}{c}
0\\ B^0_k
\end{array}
\right)v_k\\
\hphantom{\widetilde{X}_{k+1}=}\ds +\sum_{i=1}^p\left\{
\left(
\begin{array}{cc}
{C}^{0i}_k&0\\0&{C}^{0i}_k
\end{array}
\right){X}^a_k+\left(
\begin{array}{c}
{D}^{0i}_k\\0
\end{array}
\right)u_k+\left(
\begin{array}{c}
0\\ {D}^{0i}_k
\end{array}
\right)v_k
\right\}w_k^i,\\
{X}^a_{t}=\left(
\begin{array}{c}
x\\ x
\end{array}
\right),~~~k\in \mathbb{T}_t.
\end{array}
\right.
\end{eqnarray}
Noting $\widehat{X}_k=[I~~0]{X}^a_k, {X}_k=[0~~I]{X}^a_k$, the inner states of $\widehat{J}(t,x;u,v)$ and $\bar{J}(k,\bar{x};u|_{\mathbb{T}_k},v|_{\mathbb{T}_k})$ can be replaced by $\{{X}^a_k, k\in \mathbb{T}_t\}$. Hence, (\ref{widehat-J}) (\ref{bar-J}) are expressed as
\begin{eqnarray}\label{widehat-J-2}
&&\hspace{-4em}\widehat{J}(t,x;u,v)=\sum_{k=t}^{N-1}{\mathbb{E}}_t\Bigg{\{}(X_k^a)^T\left(
\begin{array}{cc}
{Q}^0_{t,k}&0\\0&0
\end{array}
\right)X^a_k
+
({\mathbb{E}}_tX^a_k)^T\left(
\begin{array}{cc}
\bar{{Q}}^0_{t,k}&0\\0&0
\end{array}
\right){\mathbb{E}}_tX^a_k\nonumber\\
&&
\hspace{-4em}\hphantom{J(t,x;u,v)=}
+
{\mathbf{u}}_k^T\left(
\begin{array}{cc}
{R}^0_{t,k}+\mu_k\Psi_k&-\mu_k\Psi_k\\-\mu_k\Psi_k&\mu_k\Psi_k
\end{array}
\right){\mathbf{u}}_k
+
({\mathbb{E}}_t{\mathbf{u}}_k)^T\left(
\begin{array}{cc}
\bar{{R}}^0_{t,k}&0\\0&0
\end{array}
\right)
{\mathbb{E}}_t{\mathbf{u}}_k
\Bigg{\}}\nonumber\\
&&\hspace{-4em}\hphantom{J(t,x;u,v)=}
+
{\mathbb{E}}_t\Bigg{\{}
(X_N^a)^T\left(
\begin{array}{cc}
{G}^0_{t}&0\\0&0
\end{array}
\right)X^a_N
+
({\mathbb{E}}_tX_N^a)^T\left(
\begin{array}{cc}
\bar{{G}}^0_{t}&0\\0&0
\end{array}
\right){\mathbb{E}}_tX^a_N
\Bigg{\}},
\end{eqnarray}
and
\begin{eqnarray}\label{bar-J-2}
&&\hspace{-4em}\bar{J}(k,{X}^a_k;u|_{\mathbb{T}_k},v|_{\mathbb{T}_k})
=\sum_{\ell=k}^{N-1}{\mathbb{E}}_k\Bigg{\{}(X^a_{\ell})^T\left(
\begin{array}{cc}
0&0\\0&{Q}^0_{k,\ell}
\end{array}
\right)X^a_{\ell}
+
({\mathbb{E}}_kX^a_{\ell})^T\left(
\begin{array}{cc}
0&0\\0&\bar{{Q}}^0_{k,\ell}
\end{array}
\right){\mathbb{E}}_kX^a_{\ell}\nonumber\\
&&\hspace{-4em}\hphantom{\bar{J}(k,{X}^a_k;u|_{\mathbb{T}_k},v|_{\mathbb{T}_k})=}
+
{\mathbf{u}}_{\ell}^T\left(
\begin{array}{cc}
0&0\\0&{R}^0_{k,\ell}
\end{array}
\right){\mathbf{u}}_{\ell}+
({\mathbb{E}}_k{\mathbf{u}}_{\ell})^T\left(
\begin{array}{cc}
0&0\\0&\bar{{R}}^0_{k,\ell}
\end{array}
\right)
{\mathbb{E}}_k{\mathbf{u}}_{\ell}
\Bigg{\}}\nonumber
\\
&&\hspace{-4em}\hphantom{\bar{J}(k,{X}^a_k;u|_{\mathbb{T}_k},v|_{\mathbb{T}_k})=}
+
\mu_k{\mathbf{u}}_{k}^T\left(
\begin{array}{cc}
\Psi_k&-\Psi_k\\ -\Psi_k&\Psi_k
\end{array}
\right){\mathbf{u}}_{k}+
{\mathbb{E}}_k\Bigg{\{}
(X_N^a)^T\left(
\begin{array}{cc}
0&0\\0& {G}^0_{k}
\end{array}
\right)X^a_N \nonumber\\
&&\hspace{-4em}\hphantom{\bar{J}(k,{X}^a_k;u|_{\mathbb{T}_k},v|_{\mathbb{T}_k})=}
+
({\mathbb{E}}_kX^a_N)^T\left(
\begin{array}{cc}
0&0\\0&\bar{{G}}^0_{k}
\end{array}
\right){\mathbb{E}}_kX^a_N
\Bigg{\}}
\end{eqnarray}
with ${\mathbf{u}}_k=[u^T_k~\,v^T_k]^T, k\in \mathbb{T}_t$. By the above notations, the game finding the open-loop self-coordination control of Problem (LQ) is sublimed to solve a generalized time-inconsistent nonzero-sum LQ dynamic game (Problem (GLQ) below). Namely, finding the open-loop self-coordination control is a motivation to study the generalized time-inconsistent nonzero-sum LQ dynamic game.




Specifically, consider the system
\begin{eqnarray}\label{system-2}
\left\{\begin{array}{l}
X_{k+1}=\big{(}A_{k}X_k+B^1_{k}u_k+{B}^2_{k}v_k\big{)}+\sum_{i=1}^p\big{(}C^i_{k}X_k+D^{1i}_{k}u_k+{D}^{2i}_{k}v_k\big{)}w^i_k, \\[1mm]
X_t=y\in \mathbb{R}^{\widetilde{n}},~~k\in  \mathbb{T}_t,~~t\in \mathbb{T}.
\end{array}
\right.
\end{eqnarray}
Here, $\{X_{k}, k\in
\widetilde{{\mathbb{T}}}_t\}\triangleq X$, $\{u_k, k\in
\mathbb{T}_t\}\triangleq u$ and $\{v_k, k\in
\mathbb{T}_t\}\triangleq v$ with $\widetilde{\mathbb{T}}_t=\{t,...,N\}$ are the state process and control
processes, respectively;
the system matrices $A_{k},C^i_{k}\in \mathbb{R}^{\widetilde{n}\times
\widetilde{n}}$, $B^1_{k},D^{1i}_{k}\in
\mathbb{R}^{\widetilde{n}\times m_1}$, $B^2_{k},D^{2i}_{k}\in
\mathbb{R}^{\widetilde{n}\times m_2}$ of (\ref{system-2}) are deterministic.
In (\ref{system-2}),
$y$ belongs to $l^2_\mathbb{F}(t; \mathbb{R}^{\widetilde{n}})$, which is defined as
\begin{eqnarray}\label{L^2-t}
l^2_\mathbb{F}(t; \mathbb{R}^{\widetilde{n}})=\Big{\{}\zeta \in \mathbb{R}^{\widetilde{n}}\,\big{|}\, \zeta\mbox{ is }\mathbb{F}_t\mbox{-measurable}, \mathbb{E}|\zeta|^2<\infty \Big{\}}.
\end{eqnarray}
The cost functionals associated with (\ref{system-2}) are
\begin{eqnarray}\label{J-1}
&&\hspace{-3em}J_1(t,y;u,v)=\sum_{k=t}^{N-1}\mathbb{E}_t\left[\left(
\begin{array}{c}
X_k\\ \textbf{u}_k
\end{array}\right)^T\left(
\begin{array}{cc}
Q_{t,k}^1&(S^1_{t,k})^T\\S^1_{t,k}&R^1_{t,k}
\end{array}\right)\left(
\begin{array}{c}
X_k\\{\textbf{u}}_k
\end{array}\right)\right.\nonumber\\
&&\hspace{-3em} \hphantom{J^1(t,x;u,v)=}+\left.
\left(
\begin{array}{c}
\mathbb{E}_tX_k\\\mathbb{E}_t\textbf{u}_k
\end{array}\right)^T\left(
\begin{array}{cc}
\bar{Q}_{t,k}^1&(\bar{S}^1_{t,k})^T\\\bar{S}^1_{t,k}&\bar{R}^1_{t,k}
\end{array}\right)\left(
\begin{array}{c}
\mathbb{E}_tX_k\\ \mathbb{E}_t\textbf{u}_k
\end{array}\right)+2(q_{t,k}^1)^TX_k+2(\rho_{t,k}^1)^T\textbf{u}_k
\right] \nonumber\\
&&\hspace{-3em}\hphantom{J^1(t,x;u,v)=}+\mathbb{E}_t\big{[}(X_N)^TG^1_tX_N\big{]}
+(\mathbb{E}_tX_N)^T \bar{G}^1_t\mathbb{E}_tX_{N}+2(g^1_t)^T\mathbb{E}_tX_{N}, 
\end{eqnarray}
and
\begin{eqnarray}\label{J-2}
&&\hspace{-3em}J_2(t,y;u,v)=\sum_{k=t}^{N-1}\mathbb{E}_t\left[\left(
\begin{array}{c}
X_k\\{\textbf{u}}_k
\end{array}\right)^T\left(
\begin{array}{cc}
Q_{t,k}^2&(S^2_{t,k})^T\\S^2_{t,k}&R^2_{t,k}
\end{array}\right)\left(
\begin{array}{c}
X_k\\{\textbf{u}}_k
\end{array}\right)\right.\nonumber\\
&&\hspace{-3em} \hphantom{J^2(t,x;u,v)=}+\left.
\left(
\begin{array}{c}
\mathbb{E}_tX_k\\\mathbb{E}_t{\textbf{u}}_k
\end{array}\right)^T\left(
\begin{array}{cc}
\bar{Q}_{t,k}^2&(\bar{S}^2_{t,k})^T\\\bar{S}^2_{t,k}&\bar{R}^2_{t,k}
\end{array}\right)\left(
\begin{array}{c}
\mathbb{E}_tX_k\\ \mathbb{E}_t{\textbf{u}}_k
\end{array}\right)+2(q_{t,k}^2)^TX_k+2(\rho_{t,k}^2)^T{\textbf{u}}_k
\right] \nonumber\\
&&\hspace{-3em}\hphantom{J^2(t,x;u,v)=}+\mathbb{E}_t\big{[}(X_N)^TG^2_tX_N\big{]}
+(\mathbb{E}_tX_N)^T \bar{G}^2_t\mathbb{E}_tX_{N}+2(g^2_t)^T\mathbb{E}_tX_{N},
\end{eqnarray}
where
\begin{eqnarray*}
&&\hspace{-1em}{\textbf{u}}_k=\left(\begin{array}{c}
u_{k}\\v_{k}
\end{array}\right),~
S^{1}_{t,k}=\left(\begin{array}{c}
S^{1(1)}_{t,k}\\[1mm] S^{1(2)}_{t,k}
\end{array}\right),~\bar{S}^{1}_{t,k}=\left(\begin{array}{c}
\bar{S}^{1(1)}_{t,k}\\[1mm] \bar{S}^{1(2)}_{t,k}
\end{array}\right), ~S^{2}_{t,k}=\left(\begin{array}{c}
S^{2(1)}_{t,k}\\[1mm] S^{2(2)}_{t,k}
\end{array}\right),~\bar{S}^{2}_{t,k}=\left(\begin{array}{c}
\bar{S}^{2(1)}_{t,k}\\[1mm] \bar{S}^{2(2)}_{t,k}
\end{array}\right), \\[1mm]
&&\hspace{-1em}R^1_{t,k}=\left(\begin{array}{cc}
R^{1(11)}_{t,k}&R^{1(12)}_{t,k}\\[1mm]
R^{1(21)}_{t,k}&R^{1(22)}_{t,k}
\end{array}\right),~\bar{R}^1_{t,k}=\left(\begin{array}{cc}
\bar{R}^{1(11)}_{t,k}&\bar{R}^{1(12)}_{t,k}\\[1mm]
\bar{R}^{1(21)}_{t,k}&\bar{R}^{1(22)}_{t,k}
\end{array}\right),~R^2_{t,k}=\left(\begin{array}{cc}
R^{2(11)}_{t,k}&R^{2(12)}_{t,k}\\[1mm]
R^{2(21)}_{t,k}&R^{2(22)}_{t,k}
\end{array}\right),~\\[1mm]
&&\hspace{-1em}\bar{R}^2_{t,k}=\left(\begin{array}{cc}
\bar{R}^{2(11)}_{t,k}&\bar{R}^{2(12)}_{t,k}\\[1mm]
\bar{R}^{2(21)}_{t,k}&\bar{R}^{2(22)}_{t,k}
\end{array}\right),~\rho^{1}_{t,k}=\left(\begin{array}{c}
\rho^{1(1)}_{t,k}\\[1mm] \rho^{1(2)}_{t,k}
\end{array}\right), \rho^{2}_{t,k}=\left(\begin{array}{c}
\rho^{2(1)}_{t,k}\\[1mm] \rho^{2(2)}_{t,k}
\end{array}\right),~~~~t\in \mathbb{T},~~k\in \mathbb{T}_t.
\end{eqnarray*}
The weighting matrices in (\ref{J-1}) (\ref{J-2}) are deterministic matrices; and $Q^j_{t,k}, \bar{Q}^j_{t,k}, R^j_{t,k}, \bar{R}^j_{t,k}, G_t^j, \bar{G}_t^j, j=1,2$ are symmetric.
Let
\begin{eqnarray}\label{L^2-T}
l^2_\mathbb{F}(\mathbb{T}_{t}; \mathbb{R}^{m_i})=\left\{\nu=\{\nu_k, k\in \mathbb{T}_{t}\}\,\Big{|}\,
\nu_k\in \mathbb{R}^{m_i}\mbox{ is }\mathbb{F}_k\mbox{-measurable},
\mathbb{E}|\nu_k|^2<\infty, k\in\mathbb{T}_{t}\right\},~~~i=1,2.
\end{eqnarray}
As the weighting matrices depend on the initial time and the nonlinear terms of conditional expectation appear in the cost functionals, the considered dynamic optimization problem associated with (\ref{system-2}) (\ref{J-1}) (\ref{J-2}) will be time-inconsistent. 

\textbf{Problem (GLQ)}. \emph{For the initial pair $(t,y)$, find a pair $({u}^*,v^*)\in l^2_\mathbb{F}(\mathbb{T}_t; \mathbb{R}^{m_1})\times l^2_\mathbb{F}(\mathbb{T}_t; \mathbb{R}^{m_2})$ such that
\begin{eqnarray}
&&\label{Problem-LQG-1}\hspace{-3em}J_1(t,y; u^*, v^*)\leq J_1(t,y;u, v^*),~~~\forall u\in l^2_{\mathbb{F}}(\mathbb{T}_t, \mathbb{R}^{m_1}),\\
&&\label{Problem-LQG-2}\hspace{-3em}J_2(k,{X}^{*}_k; u^*|_{\mathbb{T}_k}, v^*|_{\mathbb{T}_k})\leq J_2(k,{X}^{*}_k;u^*|_{\mathbb{T}_k}, (v_k, v^*|_{\mathbb{T}_{k+1}})),~~~\forall k\in \mathbb{T}_t,~ \forall v_k\in l^2_{\mathbb{F}}(k, \mathbb{R}^{m_2}),
\end{eqnarray}
where
\begin{eqnarray}\label{System-equilibrium}
\left\{
\begin{array}{l}
{X}^{*}_{k+1}=\big{(}A_{k}{X}^{*}_k+B^1_{k}u^*_k+{B}^2_{k}v^*_k\big{)} +\sum_{i=1}^p\big{(}C^i_{k}{X}^{*}_k+D^{1i}_{k}u^*_k+{D}^{2i}_{k}v^*_k\big{)}w^i_k,\\[1mm]
{X}^{*}_t=y,~~~ k\in \mathbb{T}_t.
\end{array}
\right.
\end{eqnarray}
}

$(u^*, v^*)$ above is called an {open-loop equilibrium} of Problem (GLQ). By the inequalities (\ref{Problem-LQG-1}) (\ref{Problem-LQG-2}), the best-response policies $u^*$ and $v^*$ are precommitted and time-consistent, respectively. Specifically, in (\ref{Problem-LQG-1}) $u^*$ is compared with all other elements in $l^2_{\mathbb{F}}(\mathbb{T}_t, \mathbb{R}^{m_1})$ which is global optimal on $\mathbb{T}_t$, and in contrast $v^*$ of (\ref{Problem-LQG-2}) is called an equilibrium control and is local optimal in the sense that at any time instant it is optimal pointwisely provided that future equilibrium policies are given. In Section \ref{Section-results}, Problem (GLQ) is firstly studied and the obtained results are applied directly to Problem (LQ) to obtain the open-loop self-coordination control.

\subsection{Contents and findings}

By discrete-time convex variation, the stationary conditions and convex conditions are obtained, which together ensure the existence of open-loop equilibrium of Problem (GLQ). Then, several sets of Riccati-like equations and linear equations are introduced, by which the stationary conditions and convex conditions are equivalently characterized.
Interestingly, the Riccati-like equations (\ref{U}) and linear equations (\ref{V}) characterizing the convex conditions have nothing to do with the Riccati-like equations (\ref{P}) (\ref{T}) that are for the stationary conditions. To the best of the authors, the result, b)-c) of Theorem \ref{Theorem-Problem-GLQ} that characterizes the convexity (Proposition \ref{Thm-convexity-condition}), is the first one for mean-field LQ problems, which is proved by using a technique of control shifting.
Furthermore, sufficient conditions (\ref{W-H-condition}) in terms of Riccati-like equations and linear equations are presented to characterize the open-loop equilibrium of Problem (GLQ), which can be easily checked. Moreover, the uniqueness of open-loop equilibrium is also studied. By applying the developed theory, open-loop self-coordination control of Problem (LQ) is obtained.

As an example, dynamic multi-period mean-variance portfolio selection is investigated, which itself is of much interest. %
By introducing a martingale difference sequence, the wealth equation becomes a special case of general linear stochastic system with multiplicative noises. Therefore, the developed theory of Problem (LQ) can be applied directly to this portfolio selection problem. Firstly, the result that is parallel to the general theory is presented on the open-loop self-coordination control, which resort to a pair of Riccati-like equations with special structure. When the punishment matrices become zero, the open-loop self-coordination control should be an open-loop time-consistent equilibrium control, and further the obtained Riccati-like equations should reduce to and coincide with the ones for open-loop time-consistent equilibrium control. This is validated in Theorem \ref{Theorem-Problem-mv2} indeed. Moreover, a term ``sequently generic" is introduced, and one shows that given any nonnegative definite punishment direction and any initial pair, Problem (MV) admits unique open-loop self-coordination control for sequently generic punishment intensity.

Section \ref{example} presents two examples to validate the theory developed. Numerical simulations indicate a large body of diversity and the following points are manifested. When punishment intensity is small enough, open-loop self-coordination control perform between open-loop precommitted optimal control and open-loop time-consistent equilibrium control; at some late instants, open-loop self-coordination control outperforms open-loop precommitted optimal control and open-loop time-consistent equilibrium control.
In particular, Example \ref{Exam-2} shows: the scheme of \cite{Cui-2017-2} looks for self-coordination policy between open-loop precommitted optimal control and open-loop time-inconsistent equilibrium control, while to some extent this paper goes beyond open-loop precommitted optimal control and open-loop time-consistent equilibrium control. Therefore, our formulation might be viewed as the supplement to that of \cite{Cui-2017-2}.
Furthermore, on the above technical contents, we have additional comments. Firstly, open-loop self-coordination control relives the open-loop time-inconsistent equilibrium control if one lets the punishment matrices be zero.
On the other hand, by adjusting the punishment matrices, open-loop self-coordination control provides many alternatives to handle the time inconsistency. Namely, the necessity to study open-loop self-coordination controls is indicated.

\textbf{Notation}.  For a matrix $M$,  $M^T$, $M^\dagger$ and $\mbox{Ran}(M)$ denote the transpose, the Moore-Penrose inverse and the range, respectively, of $M$. Mentioned above, $\mathbb{T}=\{0,...,N-1\}, \widetilde{\mathbb{T}}=\{0,...,N\}$, and $\mathbb{T}_t=\{t,...,N-1\}$ for $t\in \mathbb{T}$. The spaces $l^2_\mathbb{F}(t; \mathbb{R}^{\widetilde{n}})$ and $l^2_\mathbb{F}(\mathbb{T}_tt; \mathbb{R}^{m_i})$ are given in (\ref{L^2-t}) (\ref{L^2-T}); and $l^2_\mathbb{F}(t; \mathbb{R}^{n })$ and $l^2_\mathbb{F}(\mathbb{T}_tt; \mathbb{R}^{m})$ are similarly defined. If a matrix $M$ is nonnegative definite or positive definite, it will be denoted as $M\succeq0$ and $M\succ 0$.


\section{Main results}\label{Section-results}

This section presents the main results of this paper, whose proofs are given in Section \ref{Section-proof}. As looking for open-loop self-coordination control is a special case of solving Problem (GLQ), the results of this section are stated mostly for Problem (GLQ).

\begin{theorem}\label{Nece-suffic1}
For the initial pair $(t,y)$, the following statements are equivalent.
\begin{itemize}
\item[i)] Problem (GLQ) admits an open-loop equilibrium.

\item[ii)] There exists a $({u}^*,v^*)\in l^2_\mathbb{F}(\mathbb{T}_t; \mathbb{R}^{m_1})\times l^2_\mathbb{F}(\mathbb{T}_t; \mathbb{R}^{m_2})$ such that the stationary conditions
\begin{eqnarray}\label{stationary-condition}
\left\{
\begin{array}{l}
0=S^{1(1)}_{t,k}X^{*}_k+\bar{S}^{1(1)}_{t,k}\mathbb{E}_tX^{*}_k+R^{1(11)}_{t,k}u^*_k+\bar{R}^{1(11)}_{t,k} \mathbb{E}_tu^*_k+R^{1(12)}_{t,k}v^*_k+\bar{R}^{1(12)}_{t,k}\mathbb{E}_tv^*_k\\[1mm]
\hphantom{0=}+(B^1_{k})^T\mathbb{E}_kY^{*}_{k+1}+\sum_{i=1}^p(D_{k}^{1i})^T\mathbb{E}_k(Y^{*}_{k+1}w^i_k) +\rho_{t,k}^{1(1)},~~~~a.s., ~~~~k\in  \mathbb{T}_t,\\[2mm]
0=\mathcal{S}^{2(2)}_{k,k}X^{*}_k+\mathcal{R}^{2(21)}_{k,k}u^*_k+\mathcal{R}^{2(22)}_{k,k}v^*_k +(B_{k}^2)^T\mathbb{E}_kZ^{k,*}_{k+1}\\[1mm]
\hphantom{0=}+\sum_{i=1}^p(D_{k}^{2i})^T\mathbb{E}_k(Z^{k,*}_{k+1}w^i_k)+\rho^{2(2)}_{k,k},~~~~a.s., ~~~~k\in  \mathbb{T}_t,
\end{array}
\right.
\end{eqnarray}
and the convex conditions
\begin{eqnarray}\label{convex}
\left\{
\begin{array}{l}\ds
\inf_{u\in l^2_\mathbb{F}(\mathbb{T}_t; \mathbb{R}^m)} \widetilde{J}_1(t,0;u)\geq 0,~~~a.s.,\\
\ds
\inf_{v_k\in l^2_\mathbb{F}(k; \mathbb{R}^m)} \widetilde{J}_2(k,0;v_k)\geq 0,~~~a.s.,~~~\forall k\in \mathbb{T}_t
\end{array}
\right.
\end{eqnarray}
are satisfied. Here, $Y^{*}_{k+1}, Z^{k,*}_{k+1}$ are computed via the backward stochastic difference equations (BS$\Delta$Es, for short)
\begin{eqnarray}\label{system-y-k}
\left\{
\begin{array}{l}
Y^{*}_{k}=Q_{t,k}^1X^{*}_k+\bar{Q}_{t,k}^1\mathbb{E}_tX^{*}_k+\big{(}S^{1(1)}_{t,k}\big{)}^Tu^*_k +\big{(}\bar{S}_{t,k}^{1(1)}\big{)}^T\mathbb{E}_tu^*_k +\big{(}S^{1(2)}_{t,k}\big{)}^Tv^*_k+\big{(}\bar{S}^{1(2)}_{t,k}\big{)}^T\mathbb{E}_tv^*_k\\[2mm]
\hphantom{Y^{*}_{k}=}+A_{k}^T\mathbb{E}_kY^{*}_{k+1}+\sum_{i=1}^p(C_{k}^i)^T\mathbb{E}_k(Y^{*}_{k+1}w^i_k)+q_{t,k}^1,\\[2mm]
{Y}^{*}_N=G_t^1X_N^{*}+\bar{G}_t^{1}\mathbb{E}_tX_N^{*}+g_t^1,~~~~k\in  \mathbb{T}_t,
\end{array}
\right.
\end{eqnarray}
\begin{eqnarray}\label{system-z-k}
\left\{
\begin{array}{l}
\left\{
\begin{array}{l}
Z^{k,*}_{\ell}=Q_{k,\ell}^2X^{*}_\ell+\bar{Q}_{k,\ell}^2\mathbb{E}_kX^{*}_\ell+\big{(}S^{2(1)}_{k,\ell}\big{)}^Tu^*_\ell +\big{(}\bar{S}_{k,\ell}^{2(1)}\big{)}^T\mathbb{E}_ku^*_\ell +\big{(}S^{2(2)}_{k,\ell}\big{)}^Tv^*_\ell\\[1mm]
\hphantom{Z^{k,*}_{\ell}=}+\big{(}\bar{S}^{2(2)}_{k,\ell}\big{)}^T\mathbb{E}_kv^*_\ell+A_{\ell}^T\mathbb{E}_\ell Z^{k,*}_{\ell+1}+\sum_{i=1}^p(C_{\ell}^i)^T\mathbb{E}_\ell(Z^{k,*}_{\ell+1}w^i_\ell)+q_{k,\ell}^2,\\[1mm]
{Z}^{k,*}_N=G_k^2X_N^{*}+\bar{G}_k^2\mathbb{E}_kX_N^{*}+g_k^2,~~~~\ell\in  \mathbb{T}_k,
\end{array}
\right.\\
k\in \mathbb{T}_t
\end{array}
\right.
\end{eqnarray}
with
\begin{eqnarray*}
\left\{
\begin{array}{l}
{X}^{*}_{k+1}=\big{(}A_{k}{X}^{*}_k+B^1_{k}u^*_k+{B}^2_{k}v^*_k\big{)} +\sum_{i=1}^p\big{(}C^i_{k}{X}^{*}_k+D^{1i}_{k}u^*_k+{D}^{2i}_{k}v^*_k\big{)}w^i_k,\\[1mm]
{X}^{*}_t=y,~~~ k\in \mathbb{T}_t.
\end{array}
\right.
\end{eqnarray*}
$\widetilde{J}_1(t,0;u), \widetilde{J}_2(k,0;v_k)$ of (\ref{convex}) are
\begin{eqnarray}
&&\hspace{-3em}\widetilde{J}_1(t,0;u)=\sum_{k=t}^{N-1}\mathbb{E}_t\Big{[}\alpha_k^T{Q}^1_{t,k}\alpha_k+2u_k^TS^{1(1)}_{t,k}\alpha_k+u_k^TR^{1(11)}_{t,k}u_k +(\mathbb{E}_t\alpha_k)^T\bar Q_{t,k}^1\mathbb{E}_t\alpha_k\nonumber\\
&&\hspace{-3em}\hphantom{\widetilde{J}_1(t,0;u)=}+2(\mathbb{E}_tu_k)^T\bar S^{1(1)}_{t,k}\mathbb{E}_t\alpha_k+(\mathbb{E}_tu_k)^T\bar R^{1(11)}_{t,k}\mathbb{E}_tu_k\Big{]}+\mathbb{E}_t[\alpha_N^TG^1_t\alpha_N]\nonumber\\
&&\hspace{-3em}\hphantom{\widetilde{J}_1(t,0;u)=}+(\mathbb{E}_t\alpha_N)^T \bar{G}^1_t\mathbb{E}_t\alpha_N,
\end{eqnarray}
and
\begin{eqnarray}
&&\hspace{-4em}\widetilde{J}_2(k,0;v_k)=v_k^T\mathcal{R}^{2(22)}_{k,k}v_k+\sum_{\ell=k}^{N-1}\mathbb{E}_k\Big{[}\beta_\ell^TQ_{k,\ell}^2\beta_\ell+(\mathbb{E}_k\beta_\ell)^T\bar Q_{k,\ell}^2\mathbb{E}_k\beta_\ell\Big{]}\nonumber\\
&&\hspace{-4em}\hphantom{\widetilde{J}_2(t,0;v_k)=}+\mathbb{E}_k[\beta_N^TG^2_k\beta_N]+(\mathbb{E}_k\beta_N)^T\bar{G}^2_k\mathbb{E}_k\beta_N
\end{eqnarray}
with $\{\alpha_k, k\in \mathbb{T}_t\}, \{\beta_\ell, \ell\in \mathbb{T}_k\}$ given by the following stochastic difference equations (S$\Delta$Es, for short)
\begin{eqnarray}
\label{system-3}
\left\{\begin{array}{l}
\alpha_{k+1}=\big{(}A_{k}\alpha_k+B^1_{k}u_k\big{)}+\sum_{i=1}^p\big{(}C^i_{k}\alpha_k+D^{1i}_{k}u_k\big{)}w^i_k, \\[1mm]
\alpha_t=0,~~k\in  \mathbb{T}_t,
\end{array}
\right.
\end{eqnarray}
and
\begin{eqnarray}
\label{system-4}
\left\{\begin{array}{l}
\beta_{\ell+1}=A_{\ell}\beta_\ell+\sum_{i=1}^pC^i_{\ell}\beta_\ell w^i_\ell, \\[1mm]
\beta_{k+1}=B^2_{k}v_k+\sum_{i=1}^pD^{2i}_{k}v_kw^i_k, \\[1mm]
\beta_k=0,~~~\ell\in  \mathbb{T}_{k+1}.
\end{array}
\right.
\end{eqnarray}

\end{itemize}
Under i) or ii), $(u^*,v^*)$ of ii) is an open-loop equilibrium of Problem (GLQ).
\end{theorem}

To characterize the stationary conditions (\ref{stationary-condition}), introduce the Riccati-like equations:

\begin{eqnarray}\label{P}
\left\{
\begin{array}{l}
P_{t,k}=Q_{t,k}^1+A_k^TP_{t,k+1}A_k+\sum_{i,j=1}^p\delta_{k}^{ij} (C^{i}_{k})^TP_{t,k+1}C_{k}^{j}\\[2mm]
\hphantom{P_{t,k}=}-\Big{[}(H^{1(1)}_{t,k})^T~\,(H^{1(2)}_{t,k})^T\Big{]}\mathbf{W}_{t,k}^\dagger\left[
\begin{array}{c}
H^{1(1)}_{t,k}\\ \widehat{\mathcal H}^{2(2)}_{k,k}
\end{array}
\right],\\[2mm]
\mathcal{P}_{t,k}=\mathcal Q_{t,k}^1+A_{k}^T\mathcal{P}_{t,k+1}A_{k}+\sum_{i,j=1}^p\delta_{k}^{ij} (C^{i}_{k})^TP_{t,k+1}C_{k}^{j}\\[2mm]
\hphantom{\mathcal{P}_{t,k}=}-\Big{[}(\mathcal H^{1(1)}_{t,k})^T~\,(\mathcal H^{1(2)}_{t,k})^T\Big{]}
\widetilde{\mathbf{W}}_{t,k}^\dagger
\left[
\begin{array}{c}
\mathcal H^{1(1)}_{t,k}\\\mathcal H^{2(2)}_{k,k}
\end{array}
\right], \\[2mm]
\sigma_{t,k}=-\Big{[}(\mathcal H^{1(1)}_{t,k})^T~\, (\mathcal H_{t,k}^{1(2)})^T\Big{]}\widetilde{\mathbf{W}}_{t,k}^\dagger
\left[
\begin{array}{c}
h_{t,k}^1\\h_{k,k}^2
\end{array}
\right]
+A_{k}^T\sigma_{t,k+1}+q_{t,k}^1,\\[2mm]
P_{t,N}=G_t^1,  \mathcal{P}_{t,N}=\mathcal{G}_t^1,  \sigma_{t,N}=g_t^1,   ~~~~k\in  \mathbb{T}_t,
\end{array}
\right.
\end{eqnarray}
\begin{eqnarray}\label{T}
\left\{
\begin{array}{l}
\left\{
\begin{array}{l}
T_{k,\ell}={Q}_{k,\ell}^2+(A_{\ell})^T{T}_{k,\ell+1}A_{\ell}+\sum_{i,j=1}^p\delta_{\ell}^{ij} (C^{i}_{\ell})^T{T}_{k,\ell+1}C_{\ell}^{j}\\[2mm]
\hphantom{T_{k,\ell}=}-\Big{[}(H_{k,\ell}^{2(1)})^T~\,(H_{k,\ell}^{2(2)})^T\Big{]}\mathbf{W}_{t,\ell}^\dagger
\left[
\begin{array}{c}
H^{1(1)}_{t,\ell}\\ \widehat{\mathcal H}^{2(2)}_{\ell,\ell}
\end{array}
\right],\\[2mm]
\mathcal{T}_{k,\ell}=\mathcal{Q}_{k,\ell}^2+(A_{\ell})^T\mathcal{T}_{k,\ell+1}A_{\ell}+\sum_{i,j=1}^p\delta_{\ell}^{ij} (C^{i}_{\ell})^T{T}_{k,\ell+1}C_{\ell}^{j}\\[2mm]
\hphantom{\mathcal{T}_{k,\ell}=}-\Big{[}(\widehat{\mathcal H}_{k,\ell}^{2(1)})^T~\,(\widehat{\mathcal H}_{k,\ell}^{2(2)})^T\Big{]}{\mathbf{W}}_{t,\ell}^\dagger
\left[
\begin{array}{c}
H^{1(1)}_{t,\ell}\\\widehat{\mathcal H}^{2(2)}_{\ell,\ell}
\end{array}
\right],\\[2mm]
\widetilde{T}_{k,\ell}=(A_\ell)^T\widetilde{T}_{k,\ell+1}A_\ell -\Big{[}({\mathcal H}_{k,\ell}^{2(1)})^T~\,({\mathcal H}_{k,\ell}^{2(2)})^T\Big{]}\widetilde{\mathbf{W}}_{t,\ell}^\dagger\left[
\begin{array}{c}
\mathcal H^{1(1)}_{t,\ell}\\ \mathcal H^{2(2)}_{\ell,\ell}
\end{array}
\right]\\[2mm]
\hphantom{\mathcal{O}_{k,\ell}=}+\Big{[}(\widehat{\mathcal H}_{k,\ell}^{2(1)})^T~\,(\widehat{\mathcal H}_{k,\ell}^{2(2)})^T\Big{]}{\mathbf{W}}_{t,\ell}^\dagger \left[
\begin{array}{c}
H^{1(1)}_{t,\ell}\\ \widehat{\mathcal H}^{2(2)}_{\ell,\ell}
\end{array}
\right],\\[2mm]
%
%
\xi_{k,\ell}=-
\Big{[}({\mathcal H}_{k,\ell}^{2(1)})^T~\,({\mathcal H}_{k,\ell}^{2(2)})^T\Big{]}\widetilde{\mathbf{W}}_{t,\ell}^\dagger
\left[
\begin{array}{c}
h_{t,\ell}^1\\h_{\ell,\ell}^2
\end{array}
\right]
+A_{\ell}^T\xi_{k,\ell+1}+q_{k,\ell}^2,\\[2mm]
T_{k,N}=G_k^2,~~  \mathcal{T}_{k,N}=\mathcal{G}_k^2,~~   \widetilde{T}_{k,N}=0,~~  \xi_{k,N}=g_k^2, ~~~~\ell\in  \mathbb{T}_k,
\end{array}
\right.\\
k\in \mathbb{T}_t,
\end{array}
\right.
\end{eqnarray}
where
\begin{eqnarray}\label{W-H}
\left\{
\begin{array}{l}
\mathcal W^{1(1s)}_{t,k}={\mathcal{R}}_{t,k}^{1(1s)}+(B_{k}^{1})^T\mathcal{P}_{t,k+1}B_{k}^s+\sum_{i,j=1}^p\delta_{k}^{ij} (D^{1 i}_{k})^TP_{t,k+1}D_{k}^{sj},\\[2mm]
\mathcal W^{2(2s)}_{k,k}={\mathcal{R}}_{k,k}^{2(2s)}+(B_{k}^2)^T(\mathcal{T}_{k,k+1}+\widetilde{T}_{k,k+1})B_{k}^s +\sum_{i,j=1}^p\delta_{k}^{ij} (D^{2 i}_{k})^T{T}_{k,k+1}D_{k}^{sj},\\[2mm]
W^{1(1s)}_{t,k}={R}_{t,k}^{1(1s)}+(B_{k}^1)^T{P}_{t,k+1}B_{k}^s+\sum_{i,j=1}^p\delta_{k}^{ij} (D^{1 i}_{k})^T{P}_{t,k+1}D_{k}^{sj},\\[2mm]
\widehat{\mathcal{W}}^{2(2s)}_{k,k}={\mathcal{R}}_{k,k}^{2(2s)}+(B_{k}^2)^T\mathcal{T}_{k,k+1}B_{k}^s +\sum_{i,j=1}^p\delta_{k}^{ij} (D^{2 i}_{k})^T{T}_{k,k+1}D_{k}^{sj},\\[2mm]
\mathcal{H}_{t,k}^{1(s)}=\mathcal{S}_{t,k}^{1(s)}+(B_{k}^s)^T\mathcal{P}_{t,k+1}A_{k} +\sum_{i,j=1}^p\delta_{k}^{ij}(D_{k}^{sj})^T P_{t,k+1}C^{i}_{k},\\[2mm]
\mathcal{H}_{k,\ell}^{2(s)}=\mathcal{S}_{k,\ell}^{2(s)}+(B_{\ell}^s)^T (\mathcal{T}_{k,\ell+1}+\widetilde{T}_{k,\ell+1})A_{\ell}+\sum_{i,j=1}^p\delta_{\ell}^{ij}(D_{\ell}^{sj})^T {T}_{k,\ell+1}C^{i}_{\ell},\\[2mm]
\widehat{\mathcal{H}}_{k,\ell}^{2(s)}=\mathcal{S}_{k,\ell}^{2(s)}+(B_{\ell}^s)^T \mathcal{T}_{k,\ell+1}A_{\ell}+\sum_{i,j=1}^p\delta_{\ell}^{ij}(D_{\ell}^{sj})^T {T}_{k,\ell+1}C^{i}_{\ell},\\[2mm]
H_{t,k}^{1(s)}={S}_{t,k}^{1(s)}+(B_{k}^s)^T {P}_{t,k+1}A_{k}+\sum_{i,j=1}^p\delta_{k}^{ij}(D_{k}^{sj})^T {P}_{t,k+1}C^{i}_{k},\\[2mm]
H_{k,\ell}^{2(s)}={S}_{k,\ell}^{2(s)}+(B_{\ell}^s)^T {T}_{k,\ell+1}A_{\ell}+\sum_{i,j=1}^p\delta_{\ell}^{ij}(D_{\ell}^{sj})^T {T}_{k,\ell+1}C^{i}_{\ell},\\[2mm]
t\in \mathbb{T},~~~k\in  \mathbb{T}_t,~~~\ell\in  \mathbb{T}_k,~~~s=1,2,
\end{array}
\right.
\end{eqnarray}
and
\begin{eqnarray}\label{W-bf}
&&\hspace{-4em}\mathbf{W}_{t,k}=\left(
\begin{array}{ll}
W^{1(11)}_{t,k}&W^{1(12)}_{t,k}\\\widehat{\mathcal W}^{2(21)}_{k,k}&\widehat{\mathcal W}^{2(22)}_{k,k}
\end{array}\right),~~~ \widetilde{\mathbf{W}}_{t,k}=\left(
\begin{array}{ll}
\mathcal W^{1(11)}_{t,k}&\mathcal W^{1(12)}_{t,k}\\ \mathcal W^{2(21)}_{k,k}&\mathcal W^{2(22)}_{k,k}
\end{array}\right),~~~k\in  \mathbb{T}_t,\\
&&\hspace{-4em}\label{h}
h_{t,k}^1=(B_{k}^1)^T\sigma_{t,k+1}+\rho_{t,k}^{1(1)}, ~~h_{k,\ell}^2=(B_{\ell}^2)^T\xi_{k,\ell+1}+\rho_{\ell,\ell}^{2(2)},~~~~~k\in \mathbb{T}_t,~\ell\in \mathbb{T}_k.
\end{eqnarray}
Furthermore, the following Riccati equations
\begin{eqnarray}\label{U}
\left\{
\begin{array}{l}
U_{t,k}=Q_{t,k}^1+(A_k)^TU_{t,k+1}A_k+\sum_{i,j=1}^p\delta_k^{ij}(C_k^i)^TU_{t,k+1}C_k^j-M_{t,k}^TO_{t,k}^\dagger M_{t,k},\\[2mm]
\mathcal{U}_{t,k}=\mathcal{Q}_{t,k}^1+(A_k)^T\mathcal{U}_{t,k+1}A_k+\sum_{i,j=1}^p\delta_k^{ij}(C_k^i)^TU_{t,k+1}C_k^j -\mathcal{M}_{t,k}^T\mathcal{O}_{t,k}^\dagger \mathcal{M}_{t,k},\\[2mm]
U_{t,N}=G_t^1,~~\mathcal{U}_{t,N}=\mathcal{G}_t^1,~~~~k\in  \mathbb{T}_t,
\end{array}
\right.
\end{eqnarray}
and linear equations
\begin{eqnarray}\label{V}
\left\{
\begin{array}{l}
\left\{
\begin{array}{l}
V_{k,\ell}=Q_{k,\ell}^2+(A_\ell)^TV_{k,\ell+1}A_\ell+\sum_{i,j=1}^p\delta_\ell^{ij}(C_\ell^i)^TV_{k,\ell+1}C_\ell^j,\\[1mm]
\mathcal{V}_{k,\ell}=\mathcal{Q}_{k,\ell}^2+(A_\ell)^T\mathcal{V}_{k,\ell+1}A_\ell +\sum_{i,j=1}^p\delta_\ell^{ij}(C_\ell^i)^TV_{k,\ell+1}C_\ell^j,\\[1mm]
V_{k,N}=G_k^2,~~\mathcal{V}_{k,N}=\mathcal{G}_k^2,~~~~\ell\in  \mathbb{T}_k,
\end{array}
\right.\\
k\in \mathbb{T}_t
\end{array}
\right.
\end{eqnarray}
are introduced to characterize the convex conditions (\ref{convex}) with
\begin{eqnarray}\label{M-O}
\left\{
\begin{array}{l}
M_{t,k}=S_{t,k}^{1(1)}+(B_k^1)^TU_{t,k+1}A_k+\sum_{i,j=1}^p\delta_k^{ij}(D_k^{1i})^T{U}_{t,k+1}C_k^j,\\[2mm]
\mathcal{M}_{t,k}=\mathcal{S}_{t,k}^{1(1)}+(B_k^1)^T\mathcal{U}_{t,k+1}A_k+\sum_{i,j=1}^p\delta_k^{ij}(D_k^{1i})^T{U}_{t,k+1}C_k^j,\\[2mm]
O_{t,k}=R_{t,k}^{1(11)}+(B_k^1)^TU_{t,k+1}B_k^1+\sum_{i,j=1}^p\delta_k^{ij}(D_k^{1i})^T{U}_{t,k+1}D_k^{1j},\\[2mm]
\mathcal{O}_{t,k}=\mathcal{R}_{t,k}^{1(11)}+(B_k^1)^T\mathcal{U}_{t,k+1}B_k^1+\sum_{i,j=1}^p\delta_k^{ij}(D_k^{1i})^T{U}_{t,k+1}D_k^{1j},\\[2mm]
\mathbb{O}_{k,k}=\mathcal{R}^{2(22)}_{k,k}+(B_k^2)^T\mathcal{V}_{k,k+1}B_k^2 +\sum_{i,j=1}^p\delta_k^{ij}(D_k^{2i})^TV_{k,k+1}D_k^{2j}.
\end{array}
\right.
\end{eqnarray}
Throughout the paper and for a matrix $\Phi$, $\mbox{Ran}(\Phi)$ and $\mbox{Ker}(\Phi)$ denote the range and kernel of  $\Phi$, respectively.

\begin{theorem}\label{Theorem-Problem-GLQ}
For the initial pair $(t,y)$, the following statements are equivalent.
\begin{itemize}
\item[i)]Problem (GLQ) admits an open-loop equilibrium.
\item[ii)]The following assertions hold.

\begin{itemize}

\item[a)] The conditions
\begin{eqnarray}
&&\label{stationary-condition-(1)}
\widetilde{\mathbf{H}}_{t,k}\left(
\begin{array}{c}
\mathbb{E}_tX_{k}^{*}\\\mathbb{E}_tX_{k}^{*}
\end{array}\right)+\mathbf{h}_{t,k}\in \mbox{Ran}\big{(}\widetilde{\mathbf{W}}_{t,k}\big{)},\\
&&\label{stationary-condition-(2)}
\mathbf{H}_{t,k}\left(
\begin{array}{c}
X _{k}^{*}-\mathbb{E}_tX_{k}^{*}\\X _{k}^{*}-\mathbb{E}_tX_{k}^{*}
\end{array}\right)\in \mbox{Ran}\big{(}\mathbf{W}_{t,k}
\big{)},~~~k \in \mathbb{T}_t
\end{eqnarray}
are satisfied, where
\begin{eqnarray}\label{state-equilibrium}
\ds
\left\{
\begin{array}{l}
{X}^{*}_{k+1}=A_{k}{X}^{*}_k+\big{[}B^1_{k}~{B}^2_{k}\big{]}\left[
-\widetilde{\mathbf{W}}_{t,k}^\dagger\left(\widetilde{\mathbf{H}}_{t,k}\left(
\begin{array}{c}
\mathbb{E}_tX_{k}^{*}\\\mathbb{E}_tX_{k}^{*}
\end{array}\right)+\mathbf{h}_{t,k}\right)\right.\\[4mm]
\hphantom{{X}^{*}_{k+1}=}\left.-\mathbf{W}_{t,k}^\dagger \mathbf{H}_{t,k}\left(
\begin{array}{c}
X _{k}^{*}-\mathbb{E}_tX_{k}^{*}\\X _{k}^{*}-\mathbb{E}_tX_{k}^{*}
\end{array}\right)
\right] \\[4mm]
\ds\hphantom{{X}^{*}_{k+1}=}+\sum_{i=1}^p\left\{C^i_{k}{X}^{*}_k+\big{[}D^{1i}_{k}~{D}^{2i}_{k}\big{]}\left[
-\widetilde{\mathbf{W}}_{t,k}^\dagger\left(\widetilde{\mathbf{H}}_{t,k}\left(
\begin{array}{c}
\mathbb{E}_tX_{k}^{*}\\\mathbb{E}_tX_{k}^{*}
\end{array}\right)+\mathbf{h}_{t,k}\right)\right.\right.\\[4mm]
\hphantom{{X}^{*}_{k+1}=}\left.\left.-\mathbf{W}_{t,k}^\dagger \mathbf{H}_{t,k}\left(
\begin{array}{c}
X _{k}^{*}-\mathbb{E}_tX_{k}^{*}\\X _{k}^{*}-\mathbb{E}_tX_{k}^{*}
\end{array}\right)
\right]\right\}w^i_k,\\[2mm]
{X}^{*}_t=y,~~~ k\in \mathbb{T}_t,
\end{array}
\right.
\end{eqnarray}
%
%
%
%
and
\begin{eqnarray}\label{H}
\hspace{-1.5em}\mathbf{H}_{t,k}=\left(
\begin{array}{cc}
H^{1(1)}_{t,k}&0\\0&\widehat{\mathcal H}^{2(2)}_{k,k}
\end{array}\right),~
\widetilde{\mathbf{H}}_{t,k}=\left(
\begin{array}{cc}
\mathcal H^{1(1)}_{t,k}&0\\0&\mathcal H^{2(2)}_{k,k}
\end{array}\right),~\mathbf{h}_{t,k}=\left(
\begin{array}{cc}
h_{t,k}^1
\\
h_{k,k}^2
\end{array}\right),~k\in \mathbb{T}_t.
\end{eqnarray}

\item[b)] The solutions of (\ref{U}) (\ref{V}) have the property $O_{t,k}\succeq 0$, $\mathcal{O}_{t,k}\succeq0$ and $\mathbb{O}_{k,k}\succeq 0$, $k\in  \mathbb{T}_t$.

\item[c)] For any $u\in l^2_\mathbb{F}(\mathbb{T}_t; \mathbb{R}^{m_1})$, the conditions
\begin{eqnarray}\label{condition-1}
&&\hspace{-4em}M_{t,k}({\alpha}^u_{k}-\mathbb{E}_t{\alpha}^u_{k})\in \mbox{Ran}(O_{t,k}),~~ a.s.,\\
&&\hspace{-4em} \label{condition-3}
\mathcal{M}_{t,k}\mathbb{E}_t{\alpha}_{k}^u\in \mbox{Ran}(\mathcal{O}_{t,k}),~~~k\in  \mathbb{T}_t
\end{eqnarray}
are satisfied, where ${\alpha}^u$ is given by
\begin{eqnarray}\label{system-7}
\left\{\begin{array}{l}
{\alpha}^u_{k+1}=\big{(}A_{k}{\alpha}^u_k+B^1_{k}\eta^u_k\big{)} +\sum_{i=1}^p\big{(}C^i_{k}{\alpha}^u_k+D^{1i}_{k}\eta^u_k\big{)}w^i_k, \\[1mm]
{\alpha}^u_t=0,~~k\in  \mathbb{T}_t
\end{array}
\right.
\end{eqnarray}
with
\begin{eqnarray}\label{system-8}
\eta^u_k=u_k-O_{t,k}^\dagger M_{t,k}({\alpha}^u_k-\mathbb{E}_t{\alpha}_k^u)
-\mathcal{O}_{t,k}^\dagger\mathcal{M}_{t,k}\mathbb{E}_t{\alpha}^u_k,~~k\in  \mathbb{T}_t.
\end{eqnarray}

\end{itemize}

\end{itemize}
Under i) or ii), the open-loop equilibrium of Problem (GLQ) can be selected as
\begin{eqnarray}\label{u*-v*}
&&\hspace{-3em}\left(
\begin{array}{c}
u_{k}^*\\v_{k}^*
\end{array}\right)=-\widetilde{\mathbf{W}}_{t,k}^\dagger\Bigg{[}\widetilde{\mathbf{H}}_{t,k}\left(
\begin{array}{c}
\mathbb{E}_tX_{k}^{*}\\\mathbb{E}_tX_{k}^{*}
\end{array}\right)+\mathbf{h}_{t,k}\Bigg{]}-\mathbf{W}_{t,k}^\dagger \mathbf{H}_{t,k}\left(
\begin{array}{c}
X _{k}^{*}-\mathbb{E}_tX_{k}^{*}\\X _{k}^{*}-\mathbb{E}_tX_{k}^{*}
\end{array}\right)
\end{eqnarray}
with $X^*$ given in (\ref{state-equilibrium}).

\end{theorem}

\begin{remark}
The condition a) in Theorem \ref{Theorem-Problem-GLQ} is characterizing the stationary conditions (\ref{stationary-condition}), and b)-c) is equivalent to the convex conditions (\ref{convex}). To the best of the authors, b)-c) is the first result about equivalently characterizing the convexity of mean-field LQ problems.

\end{remark}

The following result is straightforward by following Theorem \ref{Theorem-Problem-GLQ}.

\begin{theorem}\label{Theorem-Problem-GLQ-2}
If the conditions
\begin{eqnarray}\label{W-H-condition}
\left\{
\begin{array}{l}
\mathbf{W}_{t,k}\mathbf{W}_{t,k}^\dagger \mathbf{H}_{t,k}=\mathbf{H}_{t,k},~~~
\widetilde{\mathbf{W}}_{t,k}\widetilde{\mathbf{W}}_{t,k}^\dagger \widetilde{\mathbf{H}}_{t,k}=\widetilde{\mathbf{H}}_{t,k},\\[1mm]
\widetilde{\mathbf{W}}_{t,k}\widetilde{\mathbf{W}}_{t,k}^\dagger \mathbf{h}_{t,k}=\mathbf{h}_{t,k},~~~
{O}_{t,k}{O}_{t,k}^\dagger {M}_{t,k}={M}_{t,k},\\[1mm]
{\mathcal{O}}_{t,k}\mathcal{{O}}_{t,k}^\dagger{\mathcal{M}}_{t,k}={\mathcal{M}}_{t,k},~~O_{t,k},\mathcal{O}_{t,k}, \mathbb{O}_{t,k}\succeq 0, ~~k\in \mathbb{T}_t,~t\in \mathbb{T}
\end{array}
\right.
\end{eqnarray}
are satisfied, then for any initial pair $(t,y)\times \mathbb{R}^{\widetilde{n}}$, Problem (GLQ) admits an open-loop equilibrium that is given in (\ref{u*-v*}).

\end{theorem}

\begin{theorem}\label{Theorem-Problem-GLQ-3}
If conditions in (\ref{W-H-condition}) are satisfied and $\mathbf{W}_{t,k}, \widetilde{\mathbf{W}}_{t,k}$ are nonsigular $k\in \mathbb{T}_t$, then Problem (GLQ) admits unique open-loop equilibrium
\begin{eqnarray*}
&&\hspace{-3em}\left(
\begin{array}{c}
u_{k}^*\\v_{k}^*
\end{array}\right)=-\widetilde{\mathbf{W}}_{t,k}^{-1}\Bigg{[}\widetilde{\mathbf{H}}_{t,k}\left(
\begin{array}{c}
\mathbb{E}_tX_{k}^{*}\\\mathbb{E}_tX_{k}^{*}
\end{array}\right)+\mathbf{h}_{t,k}\Bigg{]}-\mathbf{W}_{t,k}^{-1} \mathbf{H}_{t,k}\left(
\begin{array}{c}
X _{k}^{*}-\mathbb{E}_tX_{k}^{*}\\X _{k}^{*}-\mathbb{E}_tX_{k}^{*}
\end{array}\right)
\end{eqnarray*}
with
\begin{eqnarray*}
\left\{
\begin{array}{l}
{X}^{*}_{k+1}=\big{(}A_{k}{X}^{*}_k+B^1_{k}u^*_k+{B}^2_{k}v^*_\ell\big{)} +\sum_{i=1}^p\big{(}C^i_{k}{X}^{*}_k+D^{1i}_{k}u^*_k+{D}^{2i}_{k}v^*_k\big{)}w^i_k,\\[1mm]
{X}^{*}_t=y,~~~ k\in \mathbb{T}_t.
\end{array}
\right.
\end{eqnarray*}

\end{theorem}

If all the weighting matrices in (\ref{J-1}) (\ref{J-2}) do not depend on the initial times, this corresponds to a special case of Problem (GLQ), which is denoted as {Problem (sGLQ)} below. For Problem (sGLQ), the corresponding $P_{t,k}, \mathcal{P}_{t,k}, U_{t,k}, T_{k,\ell}, \mathcal{T}_{k,\ell}, \widetilde{T}_{k,\ell}, V_{k,\ell}, \xi_{k,\ell}, k\in \mathbb{T}_t, \ell\in \mathbb{T}_k$ of (\ref{P}) (\ref{T}) (\ref{U}) (\ref{V}) are also independent of the initial times, and are denoted, respectively, by $P_{k}, \mathcal{P}_{k}, U_{k}, T_{k}, \mathcal{T}_{k}, \widetilde{T}_{k}, V_{k}, \xi_k, k\in \mathbb{T}_t$. Furthermore, matrices in (\ref{W-H}) (\ref{W-bf}) (\ref{h}) (\ref{M-O}) (\ref{H}) do not depend on the initial times too. For example, (\ref{T}) (\ref{W-bf}) become
\begin{eqnarray*}\label{T-2}
\left\{
\begin{array}{l}
T_{k}={Q}_{k}^2+(A_{k})^T{T}_{k+1}A_{k}+\sum_{i,j=1}^p\delta_{k}^{ij} (C^{i}_{k})^T{T}_{k+1}C_{k}^{j}\\[2mm]
\hphantom{T_{k}=}-\Big{[}(H_{k}^{2(1)})^T~\,(H_{k}^{2(2)})^T\Big{]}\mathbf{W}_{k}^\dagger
\left[
\begin{array}{c}
H^{1(1)}_{k}\\ \widehat{\mathcal H}^{2(2)}_{k}
\end{array}
\right],\\[2mm]
\mathcal{T}_{k}=\mathcal{Q}_{k}^2+(A_{k})^T\mathcal{T}_{k+1}A_{k}+\sum_{i,j=1}^p\delta_{k}^{ij} (C^{i}_{k})^T{T}_{k+1}C_{k}^{j}\\[2mm]
\hphantom{\mathcal{T}_{k}=}-\Big{[}(\widehat{\mathcal H}_{k}^{2(1)})^T~\,(\widehat{\mathcal H}_{k}^{2(2)})^T\Big{]}{\mathbf{W}}_{k}^\dagger
\left[
\begin{array}{c}
H^{1(1)}_{k}\\ \widehat{\mathcal H}^{2(2)}_{k}
\end{array}
\right],\\[2mm]
\widetilde{T}_{k}=(A_k)^T\widetilde{T}_{k+1}A_k -\Big{[}({\mathcal H}_{k}^{2(1)})^T~\,({\mathcal H}_{k}^{2(2)})^T\Big{]}\widetilde{\mathbf{W}}_{k}^\dagger\left[
\begin{array}{c}
\mathcal H^{1(1)}_{k}\\ \mathcal H^{2(2)}_{k}
\end{array}
\right]\\[2mm]
\hphantom{\mathcal{O}_{k}=}+\Big{[}(\widehat{\mathcal H}_{k}^{2(1)})^T~\,(\widehat{\mathcal H}_{k}^{2(2)})^T\Big{]}{\mathbf{W}}_{k}^\dagger \left[
\begin{array}{c}
H^{1(1)}_{k}\\ \widehat{\mathcal H}^{2(2)}_{k}
\end{array}
\right],\\[2mm]
%
%
%
\xi_{k}=-\Big{[}({\mathcal H}_{k}^{2(1)})^T~\,({\mathcal H}_{k}^{2(2)})^T\Big{]}\widetilde{\mathbf{W}}_{k}^\dagger
\left[
\begin{array}{c}
h_{k}^1\\h_{k}^2
\end{array}
\right]+A_{k}^T\xi_{k+1}+q_{k}^2,\\[2mm]
T_{N}=G^2,~~  \mathcal{T}_{N}=\mathcal{G}^2,~~   \widetilde{T}_{N}=0,~~  \xi_{N}=g^2, ~~
k\in \mathbb{T}_t,
\end{array}
\right.
\end{eqnarray*}
and
\begin{eqnarray*}\label{W-bf-2}
\hspace{-2em}\mathbf{W}_{k}=\left(
\begin{array}{ll}
W^{1(11)}_{k}&W^{1(12)}_{k}\\ \widehat{\mathcal W}^{2(21)}_{k}&\widehat{\mathcal W}^{2(22)}_{k}
\end{array}\right),~~~ \widetilde{\mathbf{W}}_{k}=\left(
\begin{array}{ll}
\mathcal W^{1(11)}_{k}&\mathcal W^{1(12)}_{k}\\ \mathcal W^{2(21)}_{k}&\mathcal W^{2(22)}_{k}
\end{array}\right),~~~k\in  \mathbb{T}_t.
\end{eqnarray*}

Note that looking for open-loop self-coordination control of Problem (LQ) is a special case of solving Problem (GLQ).  
Now consider Problem (LQ). Using the notations of (\ref{J-1}) (\ref{J-2}), the weighting matrices of (\ref{widehat-J-2}) (\ref{bar-J-2}) are %
\begin{equation*}\begin{array}{l}
Q_{t,k}^1=\left(
\begin{array}{cc}
{Q}^0_{t,k}&0\\0&0
\end{array}
\right),~
\bar{Q}_{t,k}^1=\left(
\begin{array}{cc}
\bar{{Q}}^0_{t,k}&0\\0&0
\end{array}
\right),~
R_{t,k}^1=\left(
\begin{array}{cc}
{R}^0_{t,k}+\mu_k\Psi_k&-\mu_k\Psi_k\\-\mu_k\Psi_k&\mu_k\Psi_k
\end{array}
\right),\\[3mm]
\bar{R}_{t,k}^1=\left(
\begin{array}{cc}
\bar{{R}}^0_{t,k}&0\\0&0
\end{array}
\right),~
G_{t}^1=\left(
\begin{array}{cc}
{G}^0_{t}&0\\0&0
\end{array}
\right),~
\bar{G}_{t}^1=\left(
\begin{array}{cc}
\bar{{G}}^0_{t}&0\\0&0
\end{array}
\right),\\[3mm]
S^1_{t,k}=0,~\bar{S}^1_{t,k}=0,~\rho_{t,k}^1=0,~q_{t,k}^1=g^1_t=0,
\end{array}
\end{equation*}
and
\begin{equation*}\begin{array}{l}
Q_{k,\ell}^2=\left(
\begin{array}{cc}
0&0\\0&{Q}^0_{k,\ell}
\end{array}
\right),~
\bar{Q}_{k,\ell}^2=\left(
\begin{array}{cc}
0&0\\0&\bar{{Q}}^0_{k,\ell}
\end{array}
\right),~\\[3mm]
R_{k,\ell}^2=\left\{
\begin{array}{l}
\left(
\begin{array}{cc}
\mu_k\Psi_k&-\mu_k\Psi_k\\-\mu_k\Psi_k&{R}^0_{k,k}+\mu_k\Psi_k
\end{array}
\right),~~~~~~~\ell=k,\\
\left(
\begin{array}{cc}
0&0\\0&{R}^0_{k,\ell}
\end{array}
\right),~~~~~~~~~~~~~~~~~~~~~~~~~~\ell\in \mathbb{T}_{k+1},
\end{array}
\right.\\[3mm]
\bar{R}_{k,\ell}^2=\left(
\begin{array}{cc}
0&0\\0&\bar{{R}}^0_{k,\ell}
\end{array}
\right),~
G_{k}^2=\left(
\begin{array}{cc}
0&0\\0&{G}^0_{k}
\end{array}
\right),~
\bar{G}_{k}^2=\left(
\begin{array}{cc}
0&0\\0&\bar{{G}}^0_{k}
\end{array}
\right),\\[3mm]
S^2_{k,\ell}=0,~\bar{S}^2_{k,\ell}=0,~\rho_{k,\ell}^2=0,~q_{k,\ell}^2=g^2_k=0.
\end{array}
\end{equation*}
Combining (\ref{widehat-J-2}) and (\ref{bar-J-2}), we can get results that are parallel to Theorem \ref{Theorem-Problem-GLQ}, Theorem \ref{Theorem-Problem-GLQ-2} and Theorem \ref{Theorem-Problem-GLQ-3} to obtain the open-loop self-coordination control of Problem (LQ). Due to space limitations, the results are not presented here.

\section{Multi-period mean-variance portfolio selection}\label{section-mv}

In this section, we find the open-loop self-coordination control of multi-period mean-variance portfolio selection, which is a special example of Problem (LQ). Consider a capital market consisting of one riskless asset and ${p_0}$ risky assets
over a finite time horizon $N$. Let $s_k(>1)$ be a given deterministic
return of the riskless asset at time period $k$ and $e_k =
(e^1_k,\cdots,e^{p_0}_k)^T$ the vector of random returns of the ${p_0}$
risky assets at period $k$. We assume that vectors $e_k, k = 0,
1,\cdots, N-1$, are statistically independent and the only
information known about the random return vector $e_k$ is its first
two moments: its mean $\mathbb{E}(e_k) = (\mathbb{E}e^1_k,
\mathbb{E}e^2_k, \cdots , \mathbb{E}e^{p_0}_k)^T$ and its  covariance
$\mbox{Cov}(e_k)=\mathbb{E}[(e_k-\mathbb{E}e_k)(e_k-\mathbb{E}e_k)^T]$.
Clearly, $\mbox{Cov}(e_k)$ is nonnegative definite, i.e.,
$\mbox{Cov}(e_k)\succeq 0$.

Let $X_k\in \mathbb{R}$ be the wealth of the investor at the beginning
of the $k$-th period, and let $u^i_k$ be the
amount invested in the $i$-th risky asset at period $k$, $i=1,2,\cdots,{p_0}$. Then,
$X_k-\sum_{i=1}^{p_0}u_k^i$ is the amount invested in the riskless asset
at period $k$, and the wealth at the beginning of the $(k+1)$-th
period \cite{Li-Duan} is given by
\begin{eqnarray}\label{system-mean-variance}
X_{k+1}=\sum_{i=1}^{p_0}e_k^iu_k^i+\Big{(}X_k-\sum_{i=1}^{p_0}u_k^i\Big{)}s_k=s_kX_k+\Theta_k^T u_k,
\end{eqnarray}
where $\Theta_k$ is the excess return vector of risky assets \cite{Li-Duan} defined as
$\Theta_k=(\Theta_k^1,\Theta_k^2,\cdots,\Theta_k^{p_0})^T=(e_k^1-s_k,e_k^2-s_k,\cdots,e^{p_0}_k-s_k)^T$.
In this section,  we consider the case where short-selling of
stocks is allowed, i.e., $u_k^i, i=1,...,k$, take values in
$\mathbb{R}$. This leads to a multi-period mean-variance
portfolio selection formulation. For this problem, we let
$\mathbb{F}^{m}_k=\sigma(e_\ell, \ell=0,1,\cdots,k-1)$, $k=0,...,N-1$.

To proceed, (\ref{system-mean-variance}) is transformed into a
linear system with multiplicative noises such that
the general theory of above section can work. Precisely,
define
\begin{eqnarray}\label{w-D}
\left\{
\begin{array}{l}
w^i_{k}=e^i_k-s_k-\mathbb{E}(e^i_k-s_k), \\
D^{{m}i}_{k}=(0,\cdots,0,1,0,\cdots,0), \\
~i=1,\cdots,{p_0},~k\in \mathbb{T},
\end{array}
\right.
\end{eqnarray}
where the $i$-th entry of $D^{{m}i}_k$ is 1.
Then, $\{w_k=(w^1_{k},...,w^{p_0}_{k})^T, k\in \mathbb{T}\}$ is a
martingale difference sequence as $e_k, k=0,..,N-1$, are statistically independent. Furthermore,
$$\mathbb{E}_k[w_{k}w_{k}^T]=\mathbb{E}[w_{k}w_{k}^T]=\mbox{Cov}(e_{k})\triangleq(\overline{{\delta}}^{ij}_k)_{{p_0}\times {p_0}},$$
and (\ref{system-mean-variance}) becomes
\begin{eqnarray}\label{system-mean-variance-2}
\left\{
\begin{array}{l}
X_{k+1}=(s_kX_k+(\mathbb{E}\Theta_k)^T
u_k)+\sum_{i=1}^{p_0}D^{{m}i}_{k}u_kw^i_{k},\\
X_t=z,~~k\in \mathbb{T}_t.
\end{array}
\right.
\end{eqnarray}
Then, a time-inconsistent version of multi-period mean-variance problem \cite{Li-Duan} is
formulated in the following.

\textbf{Problem (MV)}. \emph{For $t\in \mathbb{T}$ and $z\in l^2_{\mathbb{F}^m}(t; \mathbb{R})$, find a $u^*\in l^2_{\mathbb{F}^m}(\mathbb{T}_t; \mathbb{R}^{p_0})$ such that
\begin{eqnarray*}
{J}_m(t,z; u^*)=\inf_{u\in l^2_{\mathbb{F}^m}(\mathbb{T}_t; \mathbb{R}^{p_0})}{J}_m(t,z; u).
\end{eqnarray*}
Here,
\begin{eqnarray}\label{cost-mean-variance-1}
{J}_m(t,z;u)=\mathbb{E}_t\big{[}(X_N-\mathbb{E}_tX_N)^2\big{]}-\lambda\mathbb{E}_tX_N
\end{eqnarray}
with $\lambda>0$ the trade-off parameter between the mean and variance of the terminal wealth.}

In what follows, we look for the open-loop self-coordination control of Problem (MV). From the formulation of (\ref{widehat-J-2}) (\ref{bar-J-2}) and using similar notations, introduce the following objective functionals and system dynamics:
%
%
%
%
\begin{eqnarray}\label{widehat-J-mv-2}
&&\hspace{-4em}\widehat{J}_m(t,z;u,v)=\sum_{k=t}^{N-1}{\mathbb{E}}_t\big{[}
{\mathbf{u}}_k^T\Upsilon_k{\mathbf{u}}_k
\big{]}+
{\mathbb{E}}_t\big{[}
(X_N^a)^TG^1X^a_N\big{]}
+
({\mathbb{E}}_tX_N^a)^T\bar{G}^1{\mathbb{E}}_tX^a_N+2(g^1)^T{\mathbb{E}}_tX^a_N,\\
&&\label{bar-J-mv-2}
\hspace{-4em}\bar{J}_m(k, {X}^a_k;u|_{\mathbb{T}_k},v|_{\mathbb{T}_k})
={\mathbf{u}}_{k}^T\Upsilon_k{\mathbf{u}}_{k}+
{\mathbb{E}}_k\big{[}
(X_N^a)^TG^2X^a_N\big{]}
+
({\mathbb{E}}_kX^a_N)^T\bar{G}^2{\mathbb{E}}_kX^a_N+2(g^2)^T{\mathbb{E}}_kX^a_N,
\end{eqnarray}
and
\begin{eqnarray}\label{X-m}
\left\{
\begin{array}{l}
{X}^{a}_{k+1}=A_k{X}^a_k+B^1_ku_k+B^2_kv_k+\sum_{i=1}^{p_0}\big{(}D^{1i}_{k}u_k+{D}^{2i}_{k}v_k\big{)}{w}^i_k,\\
{X}^a_{t}=\left(
\begin{array}{c}
z\\ z
\end{array}
\right),~~~k\in \mathbb{T}_t
\end{array}
\right.
\end{eqnarray}
with
\begin{eqnarray}
&&\label{system-parameter-mv-1}
\hspace{-3.5em}\mathbf{u}_k=\left(
\begin{array}{c}
u_k\\v_k
\end{array}
\right),~
\Upsilon_k=\mu_k\left(
\begin{array}{cc}
\Phi_k&-\Phi_k\\[1mm] -\Phi_k&\Phi_k
\end{array}
\right)\succeq 0,~G^1=\left(
\begin{array}{cc}
1&0\\0& 0
\end{array}
\right),~\bar{G}^1=\left(
\begin{array}{cc}
-1&0\\0&0
\end{array}
\right),\\
&&\label{system-parameter-mv-2}
\hspace{-3.5em} g^1=\left(
\begin{array}{c}
-{\lambda}/{2}\\0
\end{array}
\right),~G^2=\left(
\begin{array}{cc}
0&0\\0& 1
\end{array}
\right),~\bar{G}^2=\left(
\begin{array}{cc}
0&0\\0&-1
\end{array}
\right),~g^2=\left(
\begin{array}{c}
0\\-{\lambda}/{2}
\end{array}
\right),\\
&&\hspace{-3.5em}\label{A-m} A_k=\left(
\begin{array}{cc}
s_k&0\\0&s_k
\end{array}
\right), B^1_k=\left(
\begin{array}{c}
(\mathbb{E}\Theta_k)^T\\0
\end{array}
\right), B^2_k=\left(
\begin{array}{c}
0\\(\mathbb{E}\Theta_k)^T
\end{array}
\right),\\
&&\hspace{-3.5em}\label{D-m}D_k^{1i}=\left(
\begin{array}{c}
{D}^{mi}_k\\0
\end{array}
\right),~~D^{2i}_k=\left(
\begin{array}{c}
0\\ {D}^{mi}_k
\end{array}
\right),~~~i=1,...,p_0.
\end{eqnarray}
Here, $\mu_k\geq 0, \Phi_k\succeq 0, k\in \mathbb{T}$. We then have the following two results, whose proofs are given in Section \ref{Section-proof}.

\begin{theorem}\label{Theorem-Problem-mv}
Given $\{\Upsilon_k, k\in \mathbb{T}\}$, let the conditions
\begin{eqnarray}\label{W-H-mv-2}
\overline{\mathbf{W}}_{k}\overline{\mathbf{W}}_{k}^\dagger \overline{\mathbf{H}}_{k}=\overline{\mathbf{H}}_{k},~~~
\widetilde{\overline{\mathbf{W}}}_{k}\widetilde{\overline{\mathbf{W}}}_{k}^\dagger \overline{\mathbf{h}}_{k}=\overline{\mathbf{h}}_{k}, ~~~k\in \mathbb{T}
\end{eqnarray}
be satisfied, where
\begin{eqnarray}\label{W-bf-mv-0}
\left\{
\begin{array}{l}
\overline{\mathbf{W}}_k=\Upsilon_k+\left(
\begin{array}{cc}
\overline{P}_{k+1}^{(11)}\mathbb{E}\big{(}\Theta_k\Theta_k^T \big{)}&0\\[2mm]\overline{T}_{k+1}^{(21)}\mbox{Cov}(\Theta_k)&\overline{T}_{k+1}^{(22)}\mbox{Cov}(\Theta_k)
\end{array}
\right),\\[2mm]
\widetilde{\overline{\mathbf{W}}}_k=\Upsilon_k+\left(
\begin{array}{cc}
\overline{P}_{k+1}^{(11)}\mbox{Cov}(\Theta_k)&0\\[2mm]\overline{T}_{k+1}^{(21)}\mbox{Cov}(\Theta_k)& \overline{T}_{k+1}^{(22)}\mbox{Cov}(\Theta_k)
\end{array}
\right), \\[2mm]
\overline{\mathbf{H}}_k=\left(
\begin{array}{cc}
s_k\overline{P}_{k+1}^{(11)}\mathbb{E}\Theta_k&0\\0&0
\end{array}
\right)\in \mathbb{R}^{2p_0\times 4},
 \\[2mm]
k\in \mathbb{T},
\end{array}
\right.
\end{eqnarray}
\begin{eqnarray}
%
\overline{\mathbf{h}}_k=\left\{
\begin{array}{ll}
\left(
\begin{array}{c}
-\frac{\lambda }{2}\mathbb{E}\Theta_{N-1}\\[1mm]
-\frac{\lambda }{2}\mathbb{E}\Theta_{N-1}
\end{array}
\right),&~~~k=N-1,\\[4mm]
\left(
\begin{array}{c}
-\frac{\lambda}{2}s_{N-1}\mathbb{E}\Theta_{N-2}\\[1mm]
-\frac{\lambda}{2}s_{N-1}\mathbb{E}\Theta_{N-2}
\end{array}
\right),&~~~k=N-2,\\[4mm]
\left(
\begin{array}{c}
-\frac{\lambda}{2}s_{k+1}\cdots s_{N-1}\mathbb{E}\Theta_{k}\\[1mm]
-\frac{\lambda}{2}s_{k+1}\cdots s_{N-1}\mathbb{E}\Theta_{k}
\end{array}
\right),&~~~k\in \{0,...,N-3\},
\end{array}
\right.
\end{eqnarray}
and
\begin{eqnarray}\label{P-mv-0-0}
\left\{
\begin{array}{l}
\overline{P}_{k}^{(11)}=s_k^2\overline{P}_{k+1}^{(11)}\Big{[}1-\overline{P}_{k+1}^{(11)}(\mathbb{E}\Theta_k)^T \big{(}\overline{\mathbf{W}}_k^\dagger\big{)}^{(11)} \mathbb{E}\Theta_k\Big{]},\\[2mm]
\overline{T}_{k}=(A_{k})^T\overline{T}_{k+1}A_{k}-\Big{[}(\overline{H}_{k}^{2(1)})^T~\,(\overline{H}_{k}^{2(2)})^T\Big{]} \overline{\mathbf{W}}_{k}^\dagger
\left[
\begin{array}{c}
\overline{H}^{1(1)}_{k}\\ 0
\end{array}
\right],\\[2mm]
\overline{H}_{k}^{1(s)}=(B_{k}^s)^T \left(
\begin{array}{cc}
\overline{P}_{k+1}^{(11)}&0\\0&0
\end{array}
\right)A_{k},~~\overline{H}_{k}^{2(s)}=(B_{k}^s)^T \overline{T}_{k+1}A_{k},~~~s=1,2,\\[2mm]
%
%
\overline{P}_{N}^{(11)}=1,~~~  \overline{T}_{N}=G^2, ~~~k\in  \mathbb{T}
\end{array}
\right.
\end{eqnarray}
with $\overline{T}^{21}_{k+1}$ and $\overline{T}_{k+1}^{(22)}$ being the $(2,1)$-th and $(2,2)$-th entries of $\overline{T}_{k+1}$, respectively. Then, for any $(t,z)\in \mathbb{T}\times \mathbb{R}$, Problem (MV) admits an open-loop self-coordination control for the initial pair $(t,z)$ and $\{\Upsilon_k, k\in \mathbb{T}_t\}$, which is selected as
\begin{eqnarray}\label{u*-v*-mv}
&&\hspace{-3em}
v_{k}^*=-\big{[}0~~I_{p_0} \big{]}\left[\overline{\mathbf{W}}_{k}^\dagger \overline{\mathbf{H}}^{(1)}_{k}(X _{k}^{a*}-{\mathbb{E}}_tX_{k}^{a*})+\widetilde{\overline{\mathbf{W}}}_{k}^\dagger\overline{\mathbf{h}}_{k}\right],~~~k\in \mathbb{T}_t.
\end{eqnarray}
Here, $ \overline{\mathbf{H}}^{(1)}_{k}$ is the first column block of $\overline{\mathbf{H}}_{k}$, i.e.,
\begin{eqnarray*}
\overline{\mathbf{H}}^{(1)}_k=\left(
\begin{array}{cc}
s_k\overline{P}_{k+1}^{(11)}\mathbb{E}\Theta_k&0\\0&0
\end{array}
\right)\in \mathbb{R}^{2p_0\times 2},~~~k\in \mathbb{T},
\end{eqnarray*}
and
\begin{eqnarray}\label{equlibrium-system-MV}
\left\{
\begin{array}{l}
{X}^{a*}_{k+1}=\big{(}A_{k}{X}^{a*}_k+B^1_{k}u^*_k+{B}^2_{k}v^*_\ell\big{)} +\sum_{i=1}^{p_0}\big{(}D^{1i}_{k}u^*_k+{D}^{2i}_{k}v^*_k\big{)}w^i_k,\\[1mm]
{X}^{a*}_t=\left(
\begin{array}{c}
z\\z
\end{array}
\right),~~~ k\in \mathbb{T}_t
\end{array}
\right.
\end{eqnarray}
with $A_k, B_k^1, B_k^2, D_k^{1i}, D_k^{2i}, k\in \mathbb{T}_t$ given in (\ref{A-m}) (\ref{D-m}).

\end{theorem}

When the punishment matrices become zero, the open-loop self-coordination control should be an open-loop time-consistent equilibrium control, and further the obtained Riccati-like equations (\ref{P-mv-0-0}) should reduce to and coincide with the ones for open-loop time-consistent equilibrium control. This is validated in Theorem \ref{Theorem-Problem-mv2} below.

\begin{theorem}\label{Theorem-Problem-mv2}
The following statements hold.
\begin{itemize}
\item[i)] Let $\mu_k=0, k\in \mathbb{T}$. Then, $\overline{T}^{(21)}_{k+1}=0, \overline{P}^{(11)}_k>0, \overline{T}^{(22)}_k>0, k\in \mathbb{T}$
with
\begin{eqnarray}\label{T-mv-0}
\left\{
\begin{array}{l}
%
\overline{T}^{(22)}_k=s_k^2\overline{T}^{(22)}_{k+1},\\[2mm]
\overline{T}^{(22)}_{N}=1, ~~~k\in  \mathbb{T}.
\end{array}
\right.
\end{eqnarray}

\item[ii)] Let $\mu_k=0, k\in \mathbb{T}$ and $\mathbb{E}\Theta_k\in \mbox{Ran}\big{[}\mbox{Cov}(\Theta_k)\big{]}, k\in \mathbb{T}$. Then, for any initial pair $(t,z)$ Problem (MV) admits an open-loop self-coordination control, which is an open-loop time-consistent equilibrium control.

\item[iii)] Let $\mathbb{E}\Theta_k\in \mbox{Ran}\big{[}\mbox{Cov}(\Theta_k)\big{]}, k\in \mathbb{T}$. Define $\Xi_k=\big{\{}\Phi|\Phi=a_1\mbox{Cov}(\Theta_k)+a_2\mathbb{E}\Theta_k(\mathbb{E}\Theta_k)^T, a_1, a_2\geq 0\big{\}}, k\in \mathbb{T}$ and let $\Phi_k\in \Xi_k, k\in \mathbb{T}$. Then, for any initial pair Problem (MV) admits an open-loop self-coordination control.


\end{itemize}

\end{theorem}

Given any $\Phi_k \succeq0, k\in \mathbb{T}$, denote
\begin{eqnarray*}
\left|\overline{\mathbf{W}}_k\right|=\left|
\begin{array}{cc}
\mu_k\Phi_k+\overline{P}_{k+1}^{(11)}\mathbb{E}\big{(}\Theta_k\Theta_k^T \big{)}&-\mu_k\Phi_k\\[2mm]-\mu_k\Phi_k+\overline{T}_{k+1}^{(21)}\mbox{Cov}(\Theta_k) &\mu_k\Phi_k+\overline{T}_{k+1}^{(22)}\mbox{Cov}(\Theta_k)
\end{array}
\right|\equiv{\mathbf{p}}(\mu_k).
\end{eqnarray*}
The polynomial ${\mathbf{p}}(\mu_k)$ is at most of order $2p_0$. Let $\Lambda_k^0=\{\mu_k\,|\, \mathbf{p}(\mu_k)=0\}$ and $\Lambda_k=\Lambda_k^0 \cap \mathbb{R}_+$ with $\mathbb{R}_+=(0,+\infty)$. Note that $\Lambda_k$ is related to the values of $\mu_{k+1}, ...,\mu_{N-1}$, $k\in \mathbb{T}$. Clearly, $\mathbf{p}(\mu_k)\neq 0$ for $\mu_k\in \mathbb{R}_+\setminus \Lambda_k$, and the Lebesgue measure $m(\Lambda_k)=0, k\in \mathbb{T}$. In this case, we call that $\overline{\mathbf{W}}_k$ is  nonsingular for sequently generic $\{\mu_k, k\in \mathbb{T}\}$.

To be precise, a property is parameterized by $\{a_k \in \mathbb{R}_+, k\in \mathbb{T}\}$, and this property is called to hold for sequently generic $\{a_k, k\in \mathbb{T}\}$ if this property is satisfied for any $\{a_k, k\in \mathbb{T}\}$ with
$a_{N-1}\in\mathbb{R}_+\setminus \overline{\Lambda}_{N-1}, a_{N-2}\in\mathbb{R}_+\setminus \overline{\Lambda}_{N-2}$,..., and $a_{0}\in\mathbb{R}_+\setminus \overline{\Lambda}_{0}$; here, for $k\in \mathbb{T}$, $m(\overline{\Lambda}_k)=0$ and $\overline{\Lambda}_k$ is related to the values of $a_{k+1}, ..., a_{N-1}$. 

\begin{theorem}
Give any $\Phi_k \succeq0, k\in \mathbb{T}$. Then, for sequently generic $\{\mu_k, k\in \mathbb{T}\}$ and any initial pair $(t,z)$, Problem (MV) admits unique open-loop self-coordination controls.

\end{theorem}

\begin{remark}
Clearly, $\mathbb{E}\Theta_k\in \mbox{Ran}\big{[}\mbox{Cov}(\Theta_k)\big{]}, k\in \mathbb{T}$ holds if $\mbox{Cov}(\Theta_k)\succ0, k\in \mathbb{T}$. Note that $\mbox{Cov}(\Theta_k)\succ 0, k\in \mathbb{T}$ is a common assumption in multi-period mean-variance portfolio selection \cite{Cui-2017} \cite{Cui-2017-2} \cite{Li-Duan}. 
Moreover, letting $\mu_k=0, k\in \mathbb{T}$, we recover the results on open-loop time-consistent control of multi-period mean-variance portfolio selection \cite{Cui-2017} \cite{Ni-Li-Zhang-Krstic-2020} \cite{Ni-Zhang-Li}.

\end{remark}

\section{Examples}\label{example}

In this section, two examples are presented to validate the theory developed above.

\subsection{Two examples}

\begin{example}\label{Exam-1}
 Consider a discrete-time stochastic LQ problem, whose system dynamics and cost functional are given, respectively, by
\begin{eqnarray*}
\left\{
\begin{array}{l}
X^0_{k+1}=(A^0_kX^0_k+B^0_ku_k)+D^0_ku_kw_k, \\[1mm]
X^0_t=x,~~t\in \{0,1,2,3\},~~k\in \{t,...,3\},
\end{array}
\right.
\end{eqnarray*}
and
\begin{eqnarray*}\label{example-A-equ}
&&J_e(t,x;u)=\sum_{k=t}^{3}\mathbb{E}_t\big{[}(X_k^0)^TQ^0_{k}X^0_k+ (\mathbb{E}_tX_k^0)^T\bar{Q}^0_{k}\mathbb{E}_tX^0_k+ R^0_{k}u^2_k\big{]}\nonumber \\
&&\hphantom{J(t,x;u)=}+\mathbb{E}_t\big{[}(X_4^0)^TG^0X^0_4\big{]}
+(\mathbb{E}_tX^0_4)^T \bar{G}^0\mathbb{E}_tX^0_{4}, 
\end{eqnarray*}
where
\begin{eqnarray*}
&&A^0_0=\left(
\begin{array}{cc}
1& 0.4\\0.3 &2
\end{array}
\right),~~A^0_1=\left(
\begin{array}{cc}
1.102& -0.24\\ 0.53& 1.89
\end{array}
\right),~~A^0_2=\left(
\begin{array}{cc}
1.89& 0.49\\0& 1.75
\end{array}
\right),~~\\[1mm]
&&A^0_3=\left(
\begin{array}{cc}
0.8& -0.4\\0.2& 0.7
\end{array}
\right),~~B^0_0=\left(
\begin{array}{c}
1.2\\ -0.5
\end{array}
\right),~~B^0_1=\left(
\begin{array}{c}
1\\ 1
\end{array}
\right),~~B^0_2=\left(
\begin{array}{c}
1.2\\ 0.2
\end{array}
\right),~~\\[1mm]
&&B^0_3=\left(
\begin{array}{c}
1\\ 0.3
\end{array}
\right),~~D^0_0=\left(
\begin{array}{c}
1\\ 0.3
\end{array}
\right),~~D^0_1=\left(
\begin{array}{c}
1\\ 0.4
\end{array}
\right),~D^0_2=\left(
\begin{array}{c}
0.45\\ 0.25
\end{array}
\right),~\\[1mm]
&&D^0_3=\left(
\begin{array}{c}
0.52\\ 0
\end{array}
\right),~
Q^0_0=\left(
\begin{array}{cc}
    0.55 &   0.25\\
    0.25 &   0.6
\end{array}
\right),~
Q^0_1=\left(
\begin{array}{cc}
    1  & -0.325\\
   -0.325  &  0.5
\end{array}
\right),
\\[1mm]
&&Q^0_2=\left(
\begin{array}{cc}
    1.25  &  0.25\\
    0.25  &  1.4
\end{array}
\right),~Q^0_3=\left(
\begin{array}{cc}
    0.5  &       0\\
         0  &  0.375
\end{array}
\right),~~\bar{Q}_0^0=\left(
\begin{array}{cc}
    1 &   0.325\\
    0.325&   1.15
\end{array}
\right),\\
&&\bar{Q}_1^0=\left(
\begin{array}{cc}
    1.265  &  0.175\\
    0.175  &  0.95
\end{array}
\right),~~
\bar{Q}_2^0=\left(
\begin{array}{cc}
    1.25  &  0.325\\
    0.325  &  0.9
\end{array}
\right),~~
\bar{Q}_3^0=\left(
\begin{array}{cc}
    1 &        0\\
         0 &   1.5
\end{array}
\right),\\
&&
R^0_0=1.5,~~R^0_1=1.4,~~R^0_2=1.6,~~R^0_3=2,~G^0=\left(
\begin{array}{cc}
1& -0.1\\-0.1& 1
\end{array}
\right),~~\\
&&\bar{G}^0=\left(
\begin{array}{cc}
0.5 &0\\0& 0.5
\end{array}
\right),
\end{eqnarray*}
and $\{w_k, k=0,1,2,3\}$ here is a scalar martingale difference with constant second-order conditional moment $\mathbb{E}_{k}(w_k^2)=1, k=0,1,2,3$. Letting the punishment matrices be
\begin{eqnarray*}
\Upsilon_k=\mu\left(
\begin{array}{cc}
1 & 1 \\ 1&1
\end{array}
\right), ~~ k\in \{0,...,3\},
\end{eqnarray*}
find the open-loop self-coordination control for the initial pair $(0,x)$ with $x=[0.5~\, 0.5]^T$.

\end{example}

\begin{example}\label{Exam-2}
Consider a multi-period mean-variance portfolio selection problem. A capital market
consists of one riskless asset and three risky assets over a finite time horizon
$N=4$, and the parameters of the model are as follows
\begin{eqnarray*}
&&z=10,~~s_k=1.04,~~ \mathbb{E}e_k^1=1.162,~~\mathbb{E}e_k^2=1.246,~~\lambda=1,\\
&&\mathbb{E}e_k^3=1.228,~~k=0,1,2,3,
\end{eqnarray*}
and the covariance of $e_k=(e_k^1, e_k^2,e^3_k)^T$ is
\begin{eqnarray*}
\mbox{Cov}(e_k)=\left(
\begin{array}{ccc}
0.0146& 0.0187 &0.0145\\
0.0187& 0.0854 &0.0104\\
0.0145& 0.0104 &0.0289
\end{array}
\right)\succ0,~k=0,1,2,3.
\end{eqnarray*}
%
%
For objective functional of the form (\ref{cost-mean-variance-1}), let the punishment matrices be
\begin{eqnarray*}
\Upsilon_k=\mu\left(
\begin{array}{cc}
I_3&-I_3\\[1mm] -I_3&I_3
\end{array}
\right), ~~~~k\in \{0,1,2,3\}
\end{eqnarray*}
with $I_3$ the identical matrix of order 3. Find the open-loop self-coordination control for the initial pair $(0,z)$.
\end{example}

\subsection{Findings }

For Example \ref{Exam-1} and Example \ref{Exam-2}, we have a basket of policy candidates to handle the time inconsistency, namely, precommitted optimal control, open-loop time-consistent equilibrium control, open-loop self-coordination control, and self-coordination policy of \cite{Cui-2017-2}. 
Among these candidates, a question arises naturally:
Which one should we select to handle the time inconsistency?
The answer depends on whether or not we have the discretion to reconsider the problems in the future. If we are not allowed to reconsider the problems in the future, precommitted optimal control is our unique selection. On the other hand, if we are given the discretion to reconsider the problems at any or some of intermediate time points, it might be reasonable to select the one as our policy, which outperforms the others at that instant (where we lastly reconsider the problem).

In this section, the expected objective functionals will be calculated for intermediate time instants and different policy candidates. 
For a specific $\mu$, let $v^*(\mu), v^{m*}(\mu)$ be the open-loop self-coordination controls of Example \ref{Exam-1} and Example \ref{Exam-2}, respectively, and
\begin{eqnarray}\label{cost-exam1-k-1}
%
{V}_k(\mu)=\mathbb{E}\big{[}{J}_e(k,{X}^0_k;{v}^*(\mu)|_{\mathbb{T}_k})\big{]}
\end{eqnarray}
and
\begin{eqnarray}
%
\label{cost-k-1}{V}^m_k(\mu)=\mathbb{E}\big{[}{{J}}_m(k,{X}^m_k;{v}^{m*}(\mu)|_{\mathbb{T}_k})\big{]} ={\mathbb{E}}\big{[}({X}^m_N -
{\mathbb{E}}_k{X}^m_N)^2\big{]} -\lambda{\mathbb{E}}{X}^m_N
\end{eqnarray}
are the corresponding expected objective functionals at $k\in \{0,1,2,3\}$ with $\mathbb{T}_k=\{k,...,3\}$ here. In (\ref{cost-exam1-k-1})-(\ref{cost-k-1}), the internal states are computed via
\begin{eqnarray*}
%
\left\{
\begin{array}{l}
{X}^0_{k+1}=(A^0_k{X}^0_k+B^0_ku^*_k)+D^0_kv^*_k(\mu)w_k, \\[1mm]
{X}^0_0=x,~~k\in \{0,...,3\},
\end{array}
\right.
\end{eqnarray*}
and
\begin{eqnarray*}
%
%
%
%
\left\{
\begin{array}{l}
{X}^m_{k+1}=s_k{X}^m_k+\Theta^T_kv^{m*}_k(\mu),\\[1mm]
{X}^m_{0}=z,~~~k\in \{0,1,2,3\}.
\end{array}
\right.
\end{eqnarray*}
Let $v^{pr}$ be the precommitted optimal control of Example \ref{Exam-1} and
\begin{eqnarray*}
V^{pr}_k=\mathbb{E}\big{[}{J}_e(k,{X}^0_k;{v}^{pr}|_{\mathbb{T}_k})\big{]},~~~~k\in \{0,1,2,3\}
\end{eqnarray*}
with
\begin{eqnarray*}
\left\{
\begin{array}{l}
{X}^0_{k+1}=(A^0_k{X}^0_k+B^0_kv^{pr}_k)+D^0_kv^{pr}_k(\mu)w_k, \\[1mm]
{X}^0_0=x,~~k\in \{0,...,3\}.
\end{array}
\right.
\end{eqnarray*}
Then, the following facts are straightforward:
\begin{itemize}

\item[1.] When $\mu=0$, the open-loop self-coordination control $v^{*}(\mu)$ becomes the open-loop time-consistent control, namely,
    \begin{eqnarray*}
    V_k^{tc}=V_k(0),~~~k\in \{0, 1, 2, 3\},
    \end{eqnarray*}
    where $V_k^{tc}$ denotes the expected objective functional under open-loop time-consistent control.

\item[2.] It holds that
$$
V^{pr}_k=\mathbb{E}\big{[}{J}_e(k,\widehat{X}^0_k;{u}^*(0)|_{\mathbb{T}_k})\big{]}.
$$ 
Here, $u^*(0)$ is a version of $u^*$ of fictitious game (\ref{Problem-LQ-0}) (\ref{Problem-LQ-0-2}) in Problem (LQ)$_g$ that corresponds to Example \ref{Exam-1} and $\mu=0$; and $\widehat{X}^0_k$ is computed via
\begin{eqnarray*}
\left\{
\begin{array}{l}
\widehat{{X}}^0_{k+1}=(A^0_k\widehat{{X}}^0_k+B^0_ku^*(0)_k)+D^0_ku^*_k(0)w_k, \\[1mm]
\widehat{{X}}^0_0=x,~~k\in \{0,...,3\}.
\end{array}
\right.
\end{eqnarray*}

\end{itemize}
Furthermore, similar facts hold for $V^{mpr}_k, V_k^{mtc}$ of Example \ref{Exam-2}, which are the expected objective functionals under precommitted optimal control and open-loop time-consistent equilibrium control, respectively; denote the expected objective functionals under self-coordination policy of \cite{Cui-2017-2} as
\begin{eqnarray*}
V^{mc}_k(\mu)=\mathbb{E}\big{[}{{J}}_m(k, {X}^{mc}_k;{u}^{mc}(\mu)|_{\mathbb{T}_k})\big{]} ={\mathbb{E}}\big{[}({X}^{mc}_N-{\mathbb{E}}_k {X}^{mc}_N)^2\big{]} -\lambda{\mathbb{E}}{X}^{mc}_N,~~~k\in \{0,1,2,3\},
\end{eqnarray*}
where
\begin{eqnarray*}
\left\{
\begin{array}{l}
{X}^{mc}_{k+1}=s_k{X}^{mc}_k+\Theta^T_kv^{mc}_k(\mu),\\[1mm]
{X}^{mc}_{0}=z,~~~k\in \{0,1,2,3\},
\end{array}
\right.
\end{eqnarray*}
and $v_k^{mc}(\mu)=K^{mc}_k(\mu){X}^{mc}_k+L^{mc}_k(\mu), k=0,1,2,3$ with $K^{mc}_k, L^{mc}_k$ given in Theorem 3.1 of \cite{Cui-2017-2}. 

For Example \ref{Exam-1} and Example \ref{Exam-2}, we have conducted simulations with $\mu$ valued in
\begin{eqnarray*}
\Lambda^e=\big{\{}\ell\times 10^{-5}, \ell\times 10^{-3}, \ell\,\big{|}\, \ell=0,1,2,...,10^5 \big{\}}.
\end{eqnarray*}
The following tables present the minima and minimizers of ${V}_k(\mu), {V}_k^{m}(\mu), {V}^{mc}_k(\mu), k\in \{0,1,2,3\}$ over $\mu\in \Lambda^e$.
{
\renewcommand{\arraystretch}{1.5}
\begin{center}
\begin{tabular}{|c|c|c|c|c|}
\hline
\multirow{2}*{$\uparrow$}&30.0160&29.0124&26.8679&12.2209\\ \cline{2-5}
&$\mbox{min}_{\mu}{V}_0(\mu)$&$\mbox{min}_{\mu}{V}_1(\mu)$&$\mbox{min}_{\mu}{V}_2(\mu)$& $\mbox{min}_{\mu}{V}_3(\mu)$\\
\hline
\multirow{2}*{$\downarrow$}&${\mu}_0^*$&${\mu}_1^*$& ${\mu}_2^*$&${\mu}_3^*$\\ \cline{2-5}
&0&0&0.38460&1.7760\\
\hline
\end{tabular}\\
[0.25em]
Table 1. values of $\mbox{min}_{\mu}{V}_k(\mu)$ and ${\mu}^*_k=\mbox{arg}\,\mbox{min}_{\mu}{V}_k(\mu)$, $k=0,1,2,3$.    \\[1.55em]
\begin{tabular}{|c|c|c|c|c|}
\hline
\multirow{2}*{$\uparrow$}&-14.8722&-22.1273&-27.0525&-34.3649\\ \cline{2-5}
&$\mbox{min}_{\mu}{V}^m_0(\mu)$&$\mbox{min}_{\mu}{V}^m_1(\mu)$&$\mbox{min}_{\mu}{V}^m_2(\mu)$& $\mbox{min}_{\mu}{V}^m_3(\mu)$\\
\hline
\multirow{2}*{$\downarrow$}&${\mu}_0^{m*}$&${\mu}_1^{m*}$& ${\mu}_2^{m*}$&${\mu}_3^{m*}$\\ \cline{2-5}
&0.06424&0.16591&0.19802&0.22226\\
\hline
\end{tabular}\\
[0.25em]
Table 2. values of $\mbox{min}_{\mu}{V}^m_k(\mu)$ and ${\mu}^{m*}_k=\mbox{arg}\,\mbox{min}_{\mu}{V}^m_k(\mu)$, $k=0,1,2,3$.
\end{center}
\begin{center}
%
%
%
\begin{tabular}{|c|c|c|c|c|}
\hline
\multirow{2}*{$\uparrow$}&-20.6331 &-20.9993& -21.9038&-24.1135\\ \cline{2-5}
&$\mbox{min}_{\mu}{V}^{mc}_0(\mu)$&$\mbox{min}_{\mu}{V}^{mc}_1(\mu)$&$\mbox{min}_{\mu}{V}^{mc}_2(\mu)$&$\mbox{min}_{\mu}{V}^{mc}_3(\mu)$\\
\hline
\multirow{2}*{$\downarrow$}&${\mu}_0^{mc*}$&${\mu}_1^{mc*}$& ${\mu}_2^{mc*}$&${\mu}_3^{mc*}$\\ \cline{2-5}
&99953&6.13&2.163&$10^5$\\
\hline
\end{tabular}\\[0.25em]
Table 3. values of $\mbox{min}_{\mu}{V}^{mc}_k(\mu)$ and ${\mu}^{mc*}_k=\mbox{arg}\,\mbox{min}_{\mu}{V}^{mc}_k(\mu)$, $k=0,1,2,3$. %
\end{center}}
In what follows, 16 pictures are presented to show the curves of expected objective functionals.  Several points need to be specialized:

\begin{itemize}

\item[1)]  Figures 1 is for Example \ref{Exam-1} and Figures 2-4 correspond to Example \ref{Exam-2}.

\item[2)] Figure 4 shows the curves of $\min_{\mu}{V}^m_k(\mu), \min_{\mu}{V}^{mc}_k(\mu), k\in \{0,1,2,3\}$.

\item[3)] In Figures 1-3, local curves are obtained by shrinking the interval scale of time.

\end{itemize}

\begin{center}
\begin{tabular}{cc}
\includegraphics[scale=0.55]{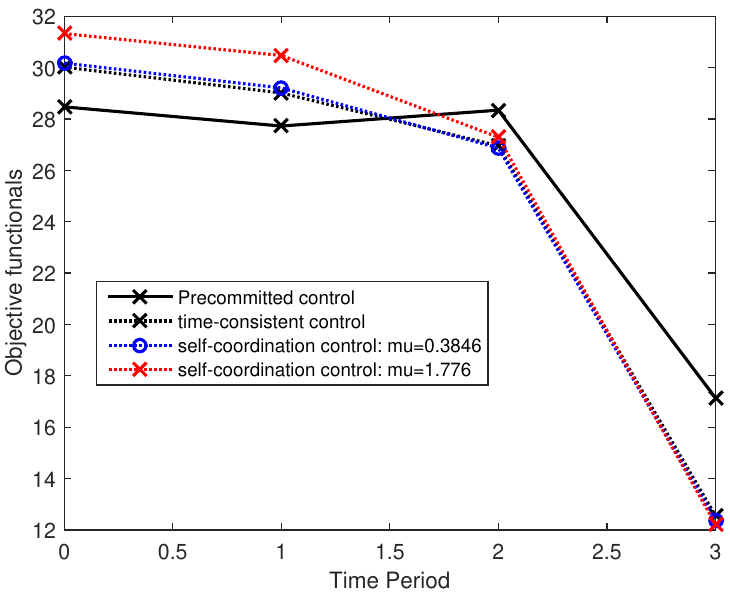}&\includegraphics[scale=0.55]{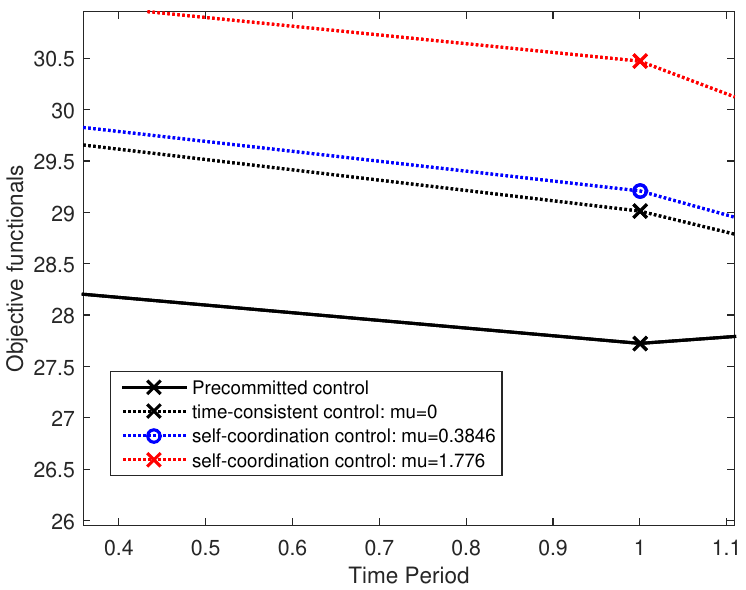}\\
\includegraphics[scale=0.55]{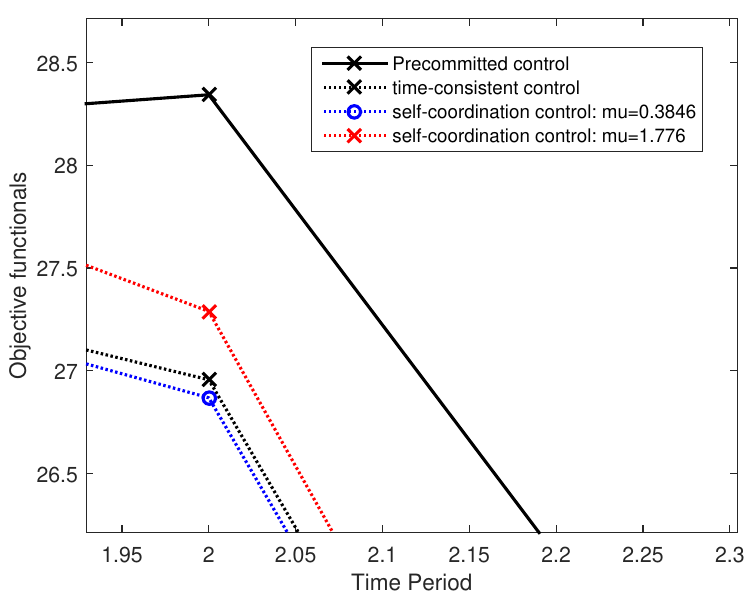}&\includegraphics[scale=0.55]{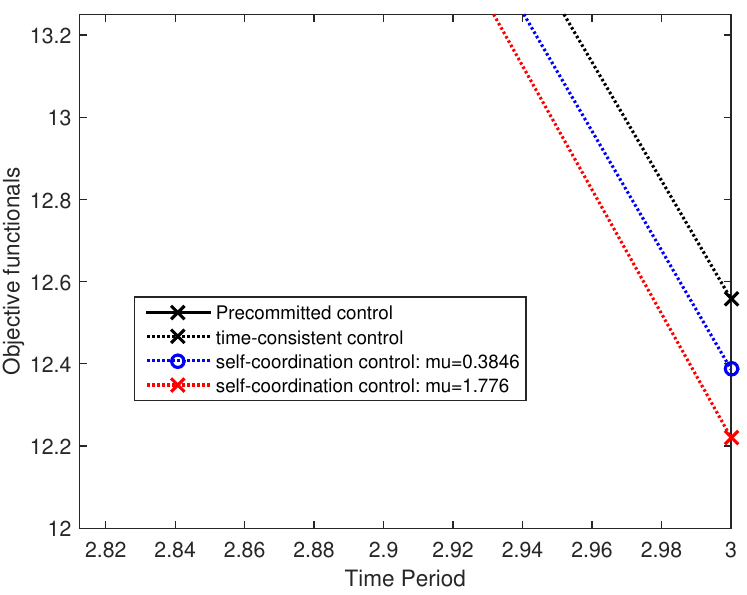}
\end{tabular}\\
{\small Figure 1: Curves and local curves of ${V}^{pr}_k, V_k^{tc}, {V}_k(\mu)$, $k=3$.}
\end{center}
%
%
%
\begin{center}
\begin{tabular}{cc}
\includegraphics[scale=0.55]{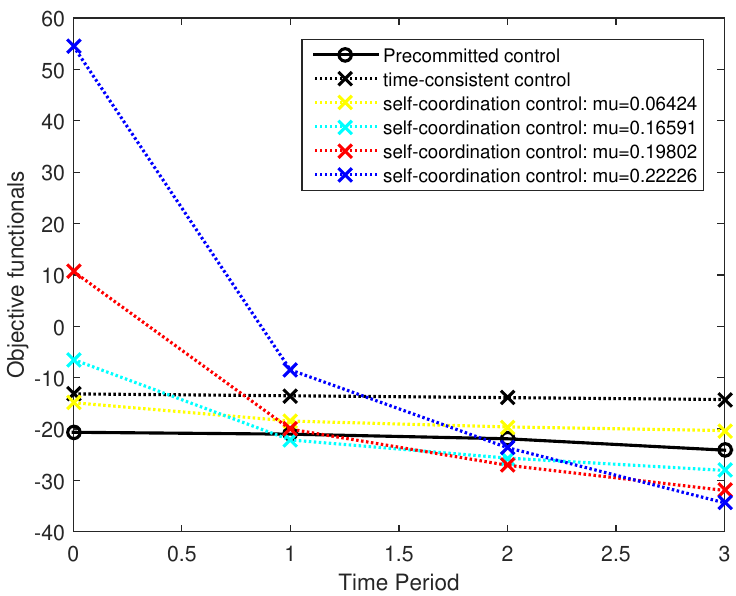}&\includegraphics[scale=0.55]{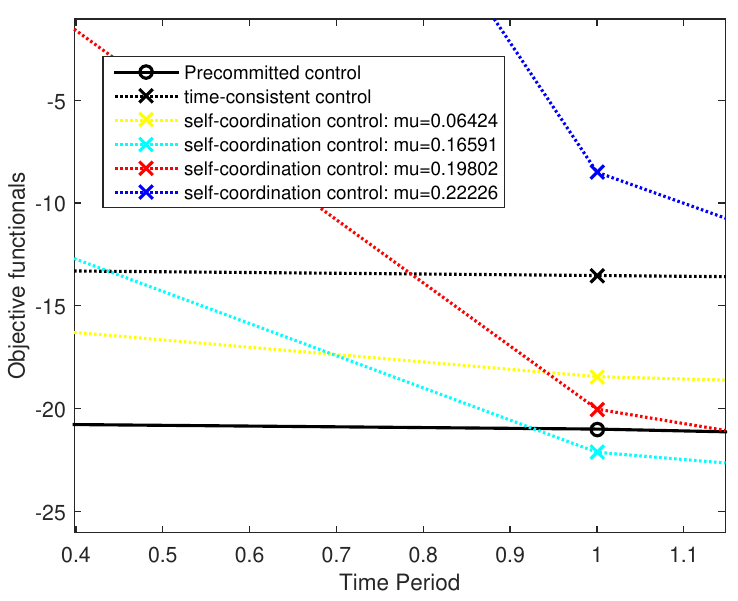}
\end{tabular}
%
\end{center}
\begin{center}
\begin{tabular}{cc}
\includegraphics[scale=0.55]{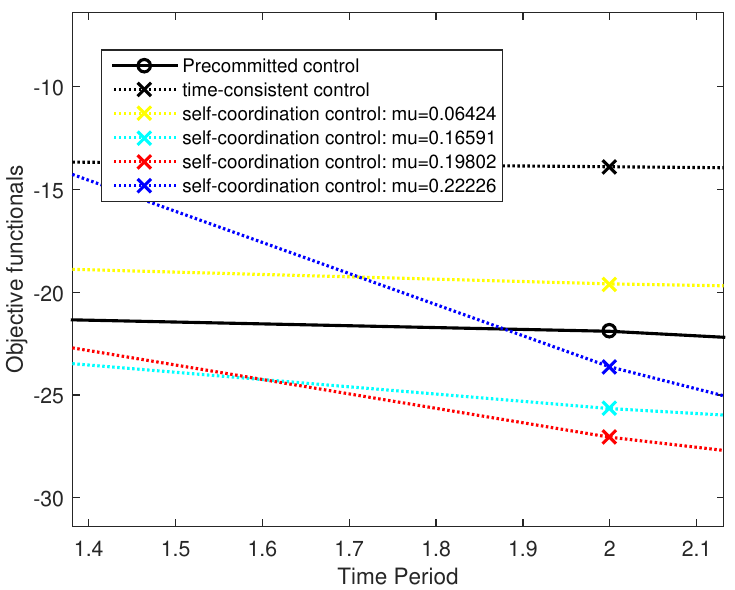}&\includegraphics[scale=0.55]{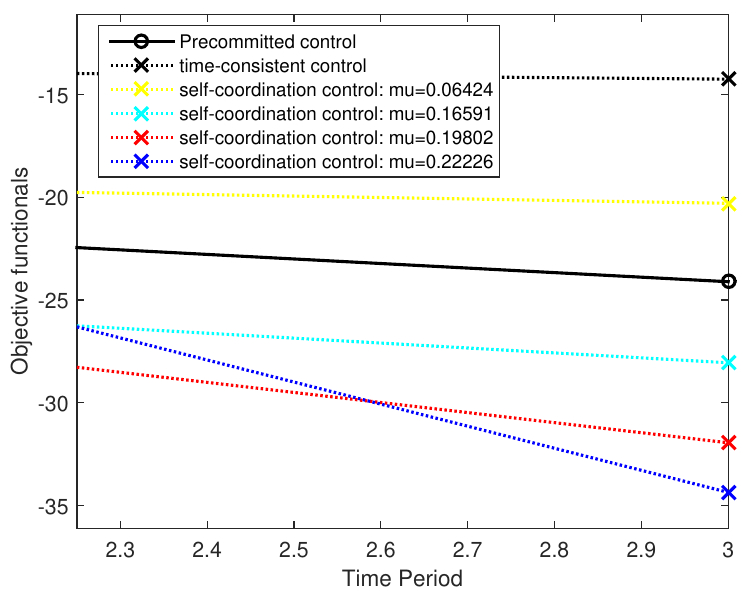}
\end{tabular}\\
{\small Figure 2: Curves and local curves of ${V}^{mpr}_k, V_k^{mtc}, {V}^m_k(\mu)$, $k=0,1,2,3$.}
\end{center}
\vspace{0.25em}
\begin{center}
\begin{tabular}{cc}
\includegraphics[scale=0.55]{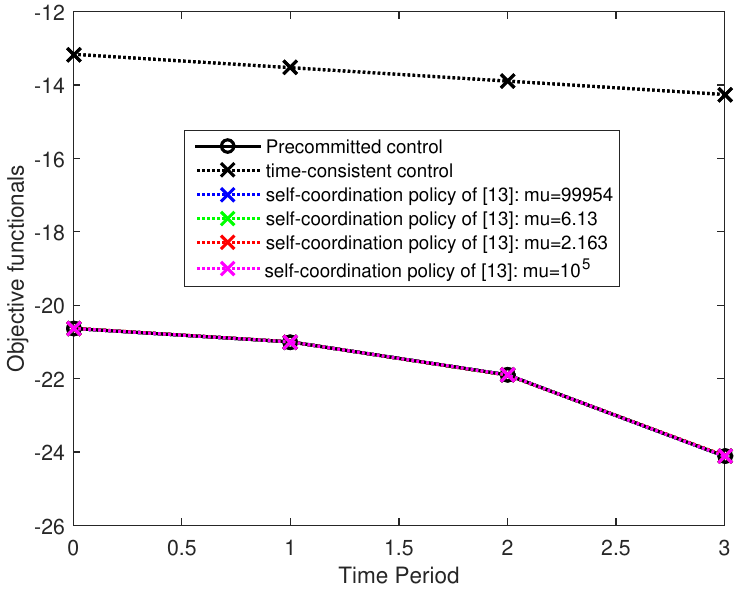}&\includegraphics[scale=0.55]{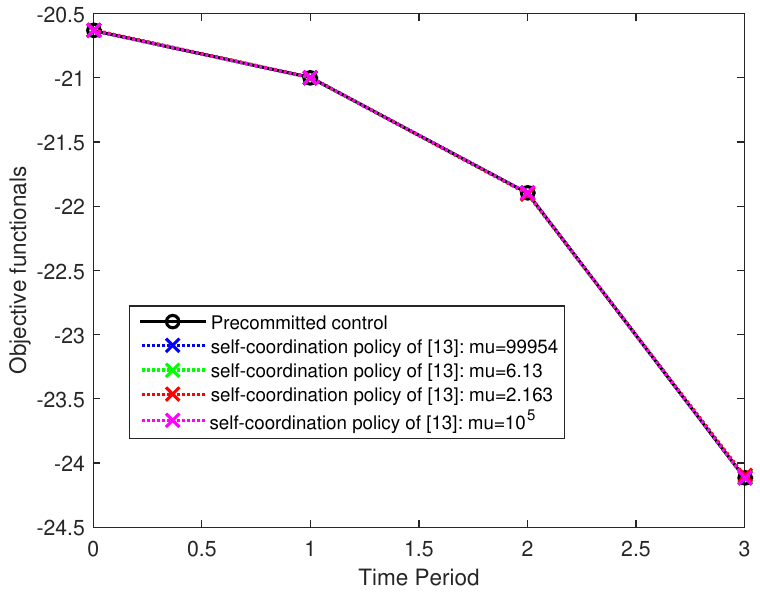}\\
\includegraphics[scale=0.55]{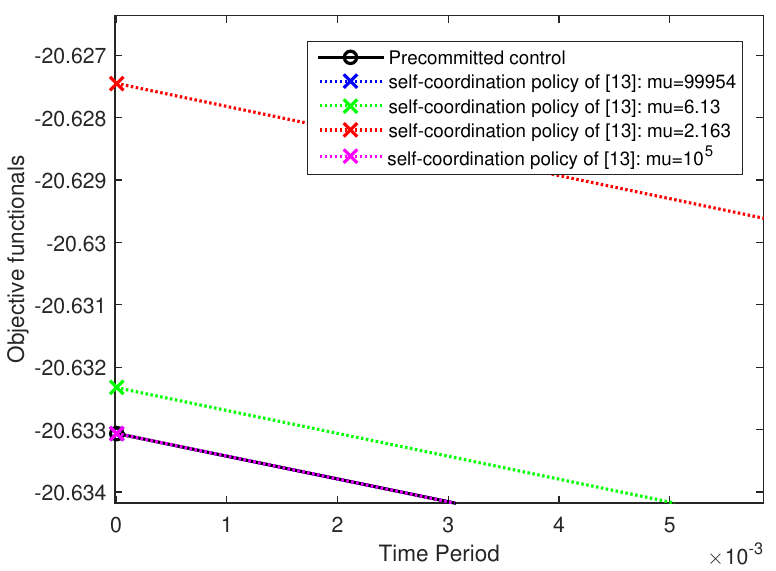}&\includegraphics[scale=0.55]{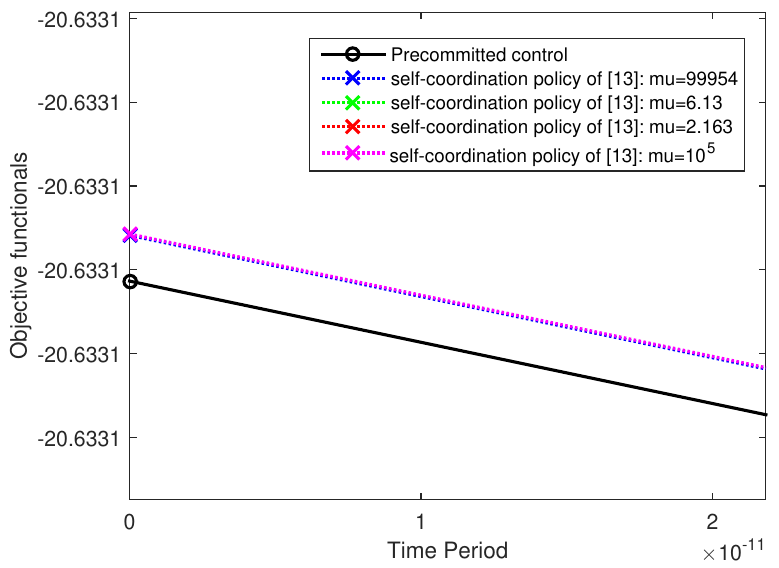}
%
\end{tabular}
%
\end{center}
\begin{center}
\begin{tabular}{cc}
%
\includegraphics[scale=0.55]{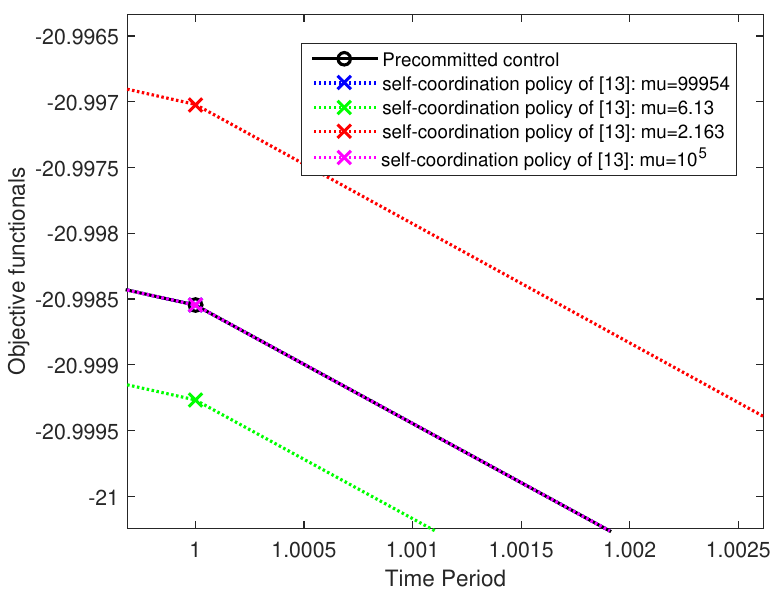}&\includegraphics[scale=0.55]{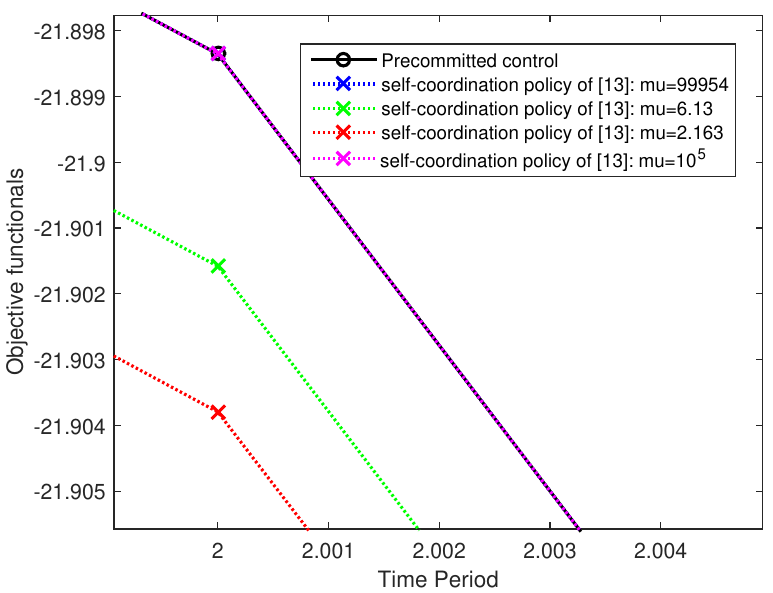}\\
\includegraphics[scale=0.55]{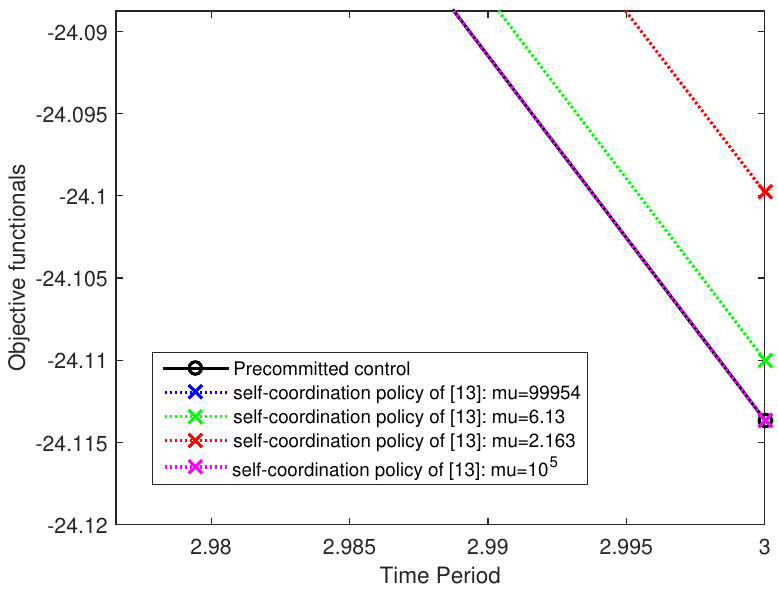}&\includegraphics[scale=0.55]{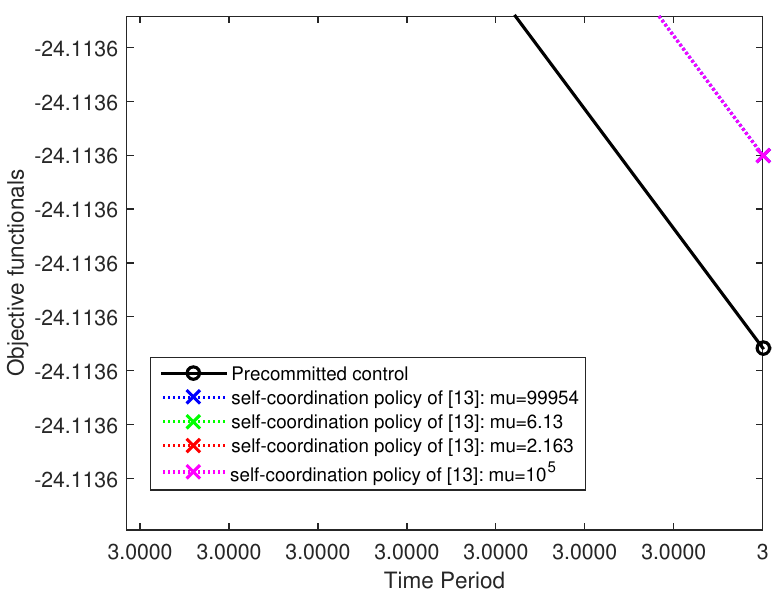}
\end{tabular}\\
{\small Figure 3: Curves and local curves of ${V}^{mpr}_k, V_k^{mtc}, {V}^{mc}_k(\mu)$, $k=1,2,3$.}
\end{center}
%
%
%
%
%
%
%
\vspace{0.25em}
\begin{center}
\includegraphics[scale=0.723]{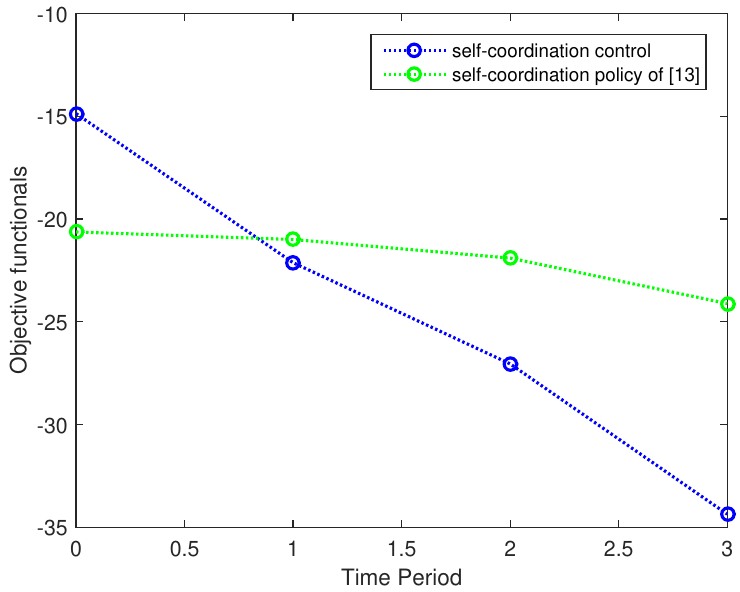}\\
{\small Figure 4: Curves of $\min_{\mu}{{V}}^{m}_k(\mu), \min_{\mu}{V}^{mc}_k(\mu)$, $k=0, 1,2,3$.}
\end{center}

%
%
%
%

%
%
Above figures indicate a large body of diversity  and the following cases are manifested:

\begin{itemize}

\item[1)] At some early instants, open-loop precommitted optimal control outperforms (with smaller expected objective functional) open-loop time-inconsistent equilibrium control; and at late instants, open-loop time-inconsistent equilibrium control outperforms open-loop precommitted optimal control: Figure 1.

\item[2)]At some late instants, open-loop self-coordination control outperforms open-loop precommitted optimal control and open-loop time-consistent equilibrium control: Figures 1 ($\mu=0.3846, 1.776$), Figure 2 ($\mu=0.16591, 0.19802, 0.22226$).

\item[3)] When $\mu$ is small enough, open-loop self-coordination control performs between open-loop precommitted optimal control and open-loop time-consistent equilibrium control: Figure 2 ($\mu=0.06424$).

\item[4)] Comparison with \cite{Cui-2017-2}: Example \ref{Exam-2}
\begin{itemize}

\item[a)] Self-coordination policy of \cite{Cui-2017-2} performs mostly between open-loop precommitted optimal control and open-loop time-consistent equilibrium control (Figure 3), though the minima of $V^{mc}_k(\mu), k=1, 2$ are slightly smaller than those of open-loop precommitted optimal control. 
    \\[-1em]


\item[] At $k=3$, open-loop precommitted optimal control outperforms all the self-coordination policies of \cite{Cui-2017-2} over $\Lambda^e$: subfigures 3, 4 of Figure 3.\\[-0.8em]

\item[b)] Concerned with the minima of expected objective functionals over $\mu\in \Lambda^e$, open-loop self-coordination control of this paper outperform self-coordination policy \cite{Cui-2017-2} at $k=1,2,3$: Figure 4.

    \vspace{0.45em}

    Namely, we can select open-loop self-coordination controls such that at $k=1,2,3$ they outperform self-coordination policy of \cite{Cui-2017-2}; see the line ($\mu=0.16591$) of Figure 2, and compare with those of Figure 3.

    \vspace{0.45em}

    This, yet, is at a price that self-coordination policy of \cite{Cui-2017-2} outperforms open-loop self-coordination control at $k=0$.\\[-0.8em]

\item[c)] To summarize and for the particular example (Example \ref{Exam-2}), self-coordination policy of \cite{Cui-2017-2} performs mostly between open-loop precommitted optimal control and open-loop time-inconsistent equilibrium control (this is also indicated in Page 102 of \cite{Cui-2017-2}); at a price at $k=0$, particular open-loop self-coordination controls can be selected such that at $k=1,2,3$ they outperform self-coordination policy of \cite{Cui-2017-2}, and also outperform open-loop precommitted optimal control and open-loop time-consistent equilibrium control.

    \vspace{0.45em}

    Namely, the scheme of \cite{Cui-2017-2} looks for self-coordination policy between open-loop precommitted optimal control and open-loop time-inconsistent equilibrium control, while to some extent this paper goes beyond open-loop precommitted optimal control and open-loop time-consistent equilibrium control. Therefore, our formulation might be viewed as some supplement to that of \cite{Cui-2017-2}, and adds a new dimension to handle time-inconsistent optimal control problems.

\end{itemize}

\end{itemize}


The last words of this section go to the necessity to study open-loop self-coordination control.

\section{Proofs}\label{Section-proof}

\subsection{Proof of Theorem \ref{Nece-suffic1}}

The proof is based on the method of discrete-time convex variation.

\emph{i)$\Rightarrow$ii)}. Let $({u}^*,v^*)\in l^2_\mathbb{F}(\mathbb{T}_t; \mathbb{R}^{m_1})\times l^2_\mathbb{F}(\mathbb{T}_t; \mathbb{R}^{m_2})$ be an open-loop equilibrium.  For $\varepsilon \in \mathbb{R}$ and $u\in l^2_{\mathbb{F}}(\mathbb{T}_t, \mathbb{R}^{m_1})$, let $X^{\varepsilon}$ satisfy the S$\Delta$E
\begin{eqnarray}\label{system-5}
\left\{\begin{array}{l}
X^{\varepsilon}_{k+1}=\big{(}A_{k}X^{\varepsilon}_k+B^1_{k}(u^*_k+\varepsilon u_k)+{B}^2_{k}v^*_k\big{)} \\
\hphantom{X^{t,\varepsilon}_{k+1}=}+\sum_{i=1}^p\big{(}C^i_{k}X^{\varepsilon}_k+D^{1i}_{k}(u^*_k+\varepsilon u_k)+{D}^{2i}_{k}v^*_k\big{)}w^i_k, \\[2mm]
X^{\varepsilon}_t=y, ~~k\in  \mathbb{T}_t.
\end{array}
\right.
\end{eqnarray}
From (\ref{System-equilibrium}) and (\ref{system-5}), we have
\begin{eqnarray*}
\left\{
\begin{array}{l}
\frac{X^{\varepsilon}_{k+1}-X^{*}_{k+1}}{\varepsilon}=A_{k}\frac{X^{\varepsilon}_k-X^{*}_k}{\varepsilon}+B^1_{k}u_k
+\sum_{i=1}^p\big{(}C^i_{k}\frac{X^{\varepsilon}_k-X^{*}_k}{\varepsilon}+D^{1i}_{k}u_k\big{)}w^i_k,\\[2mm]
\frac{X^{\varepsilon}_t-X^{*}_t}{\varepsilon}=0,~~k\in  \mathbb{T}_t.
\end{array}
\right.
\end{eqnarray*}
Denote $\frac{X^{\varepsilon}_k-X^{*}_k}{\varepsilon}$ by $\alpha_k$, then $\alpha=\{\alpha_k, ~k\in  \mathbb{T}_t\}$ satisfies (\ref{system-3}).
Obviously, $X^{\varepsilon}_k=X^{*}_k+\varepsilon\alpha_k, k\in \mathbb{T}_t$. Then, we obtain
\begin{eqnarray}\label{Nece-suffic1-1}
\hspace{-0.35em}
&&0\leq J_1(t,y;u^*+\varepsilon u, v^*)-J_1(t,y; u^*, v^*)\nonumber\\
&&\hphantom{0}=2\varepsilon\Big{\{}\mathbb{E}_t\big{[}(X^{*}_N)^TG^1_t\alpha_N\big{]}+(\mathbb{E}_t\alpha_N) ^T\bar{G}^1_t\mathbb{E}_tX^{*}_N+(g_t^1)^T\mathbb{E}_t\alpha_N\Big{\}}\nonumber\\
&&\hphantom{0\leq }\hphantom{=}+\varepsilon^2\Big{\{}\mathbb{E}_t[\alpha_N^TG^1_t\alpha_N]+(\mathbb{E}_t\alpha_N)^T \bar{G}^1_t\mathbb{E}_t\alpha_N\Big{\}}\nonumber\\
\ds
&&\hphantom{0\leq }
+\varepsilon\sum_{k=t}^{N-1}\mathbb{E}_t\Bigg{\{}\Bigg{[}\left(
\begin{array}{c}
\alpha_k\\{u}_k\\0
\end{array}\right)^T\left(
\begin{array}{ccc}
Q_{t,k}^1&(S^{1(1)}_{t,k})^T&(S^{1(2)}_{t,k})^T\\S^{1(1)}_{t,k}&R^{1(11)}_{t,k}&R^{1(12)}_{t,k}\\
S^{1(2)}_{t,k}&R^{1(21)}_{t,k}&R^{1(22)}_{t,k}
\end{array}\right)\left(
\begin{array}{c}
2X^{*}_k+\varepsilon\alpha_k\\2u^*_k+\varepsilon u_k\\2v^*_k
\end{array}\right)\nonumber\\
&&\hphantom{0\leq }+
\left(
\begin{array}{c}
\mathbb{E}_t\alpha_k\\\mathbb{E}_t{u}_k\\0
\end{array}\right)^T\left(
\begin{array}{ccc}
\bar Q_{t,k}^1&(\bar S^{1(1)}_{t,k})^T&(\bar S^{1(2)}_{t,k})^T\\ \bar S^{1(1)}_{t,k}&\bar R^{1(11)}_{t,k}&\bar R^{1(12)}_{t,k}\\
\bar S^{1(2)}_{t,k}&\bar R^{1(21)}_{t,k}&\bar R^{1(22)}_{t,k}
\end{array}\right)\left(
\begin{array}{c}
2\mathbb{E}_tX^{*}_k+\varepsilon\mathbb{E}_t\alpha_k\\2\mathbb{E}_tu^*_k+\varepsilon \mathbb{E}_tu_k\\2\mathbb{E}_tv^*_k
\end{array}\right)
\Bigg{]}\nonumber\\
&&\hphantom{0\leq }+2(q_{t,k}^1)^T\alpha_k +2(\rho_{t,k}^{1(1)})^Tu_k
\Bigg{\}} \nonumber\\
\ds
&&\hphantom{0}=2\varepsilon\Bigg{\{}\mathbb{E}_t\left[(G^1_tX^{*}_N+g^1_t)^T\alpha_N\right] +(\mathbb{E}_tX^{*}_N)^T\bar{G}^1_t\mathbb{E}_t\alpha_N\nonumber\\
\ds
&&\hphantom{0\leq }
+ \sum_{k=t}^{N-1}\mathbb{E}_t\Big{[}\Big{(}Q_{t,k}^1X^{*}_k+\big{(}S^{1(1)}_{t,k}\big{)}^Tu^*_k +\big{(}S^{1(2)}_{t,k}\big{)}^Tv^*_k+q_{t,k}^1\Big{)}^T \alpha_k\nonumber\\[1mm]
&&\hphantom{0\leq }+\Big{(}S^{1(1)}_{t,k}X^{*}_k+R^{1(11)}_{t,k}u^*_k+ R^{1(12)}_{t,k}v^*_k+\rho^{1(1)}_{t,k}\Big{)}^Tu_k\nonumber\\
&&\hphantom{0\leq }+\left(\bar Q_{t,k}^1\mathbb{E}_tX^{*}_k+\big{(}\bar S^{1(1)}_{t,k}\big{)}^T\mathbb{E}_tu^*_k+\big{(}\bar{S}^{1(2)}_{t,k}\big{)}^T\mathbb{E}_tv^*_k\right)^T\mathbb{E}_t\alpha_k \nonumber\\
&&\hphantom{0\leq }+\left(\bar S^{1(1)}_{t,k}\mathbb{E}_tX^{*}_k+\bar R^{1(11)}_{t,k}\mathbb{E}_tu^*_k+\bar R^{1(12)}_{t,k}\mathbb{E}_tv^*_k\right)^T\mathbb{E}_tu_k\Big{]}\Bigg{\}}\nonumber\\
\ds
&&\hphantom{0\leq }+\varepsilon^2\Bigg{\{}\mathbb{E}_t\left[\alpha_N^TG^1_t\alpha^t_N\right]+(\mathbb{E}_t\alpha_N)^T \bar{G}^1_t\mathbb{E}_t\alpha_N\nonumber\\
\ds
&&\hphantom{0\leq }+\sum_{k=t}^{N-1}\mathbb{E}_t\Big{[}\alpha_k^T{Q}^1_{t,k}\alpha_k +2u_k^TS^{1(1)}_{t,k}\alpha_k+u_k^TR^{1(11)}_{t,k}u_k\nonumber\\
&&\hphantom{0\leq }+(\mathbb{E}_t\alpha_k)^T\bar Q_{t,k}^1\mathbb{E}_t\alpha_k+2(\mathbb{E}_tu_k)^T\bar S^{1(1)}_{t,k}\mathbb{E}_t\alpha_k+(\mathbb{E}_tu_k)^T\bar R^{1(11)}_{t,k}\mathbb{E}_tu_k\Big{]}\Bigg{\}}.
\end{eqnarray}
Noting (\ref{system-y-k}) and (\ref{system-3}), we have
\begin{eqnarray*}
&&\hspace{-2.5em}\mathbb{E}_t\left[(G^1_tX^{*}_N+g^1_t)^T\alpha_N\right] +(\mathbb{E}_tX^{*}_N)^T\bar{G}^1_t\mathbb{E}_t\alpha_N
+ \sum_{k=t}^{N-1}\mathbb{E}_t\Big{\{}\Big{(}Q_{t,k}^1X^{*}_k+\big{(}S^{1(1)}_{t,k}\big{)}^Tu^*_k +\big{(}S^{1(2)}_{t,k}\big{)}^Tv^*_k+q_{t,k}^1\Big{)}^T \alpha_k\nonumber\\
&&\hspace{-2.5em} +\Big{(}S^{1(1)}_{t,k}X^{*}_k+R^{1(11)}_{t,k}u^*_k+ R^{1(12)}_{t,k}v^*_k+\rho^{1(1)}_{t,k}\Big{)}^Tu_k+\left(\bar Q_{t,k}^1\mathbb{E}_tX^{*}_k+\big{(}\bar S^{1(1)}_{t,k}\big{)}^T\mathbb{E}_tu^*_k+\big{(}\bar S^{1(2)}_{t,k}\big{)}^T\mathbb{E}_tv^*_k\right)^T\mathbb{E}_t\alpha_k\nonumber\\
&&\hspace{-2.5em}+\left(\bar S^{1(1)}_{t,k}\mathbb{E}_tX^{*}_k+\bar R^{1(11)}_{t,k}\mathbb{E}_tu^*_k+\bar R^{1(12)}_{t,k}\mathbb{E}_tv^*_k\right)^T\mathbb{E}_tu_k\Big{\}}\nonumber\\
&&\hspace{-2.5em}=\mathbb{E}_t\sum_{k=t}^{N-1}\Bigg{\{}\Big{[}Q_{t,k}^1(X^{*}_k-\mathbb{E}_tX^{*}_k) +\big{(}S^{1(1)}_{t,k}\big{)}^T(u^*_k-\mathbb{E}_tu^*_k)+\big{(}S^{1(2)}_{t,k}\big{)}^T(v^*_k-\mathbb{E}_tv^*_k)+
A_{k}^T(\mathbb{E}_kY^{*}_{k+1}-\mathbb{E}_tY^{*}_{k+1})\nonumber\\
&&\hspace{-2.5em}\hphantom{=}+\sum_{i=1}^p(C_{k}^i)^T\Big{(}\mathbb{E}_k(Y^{*}_{k+1}w^i_k)-\mathbb{E}_t(Y^{*}_{k+1}w^i_k)\Big{)} -(Y^{*}_k-\mathbb{E}_tY^{*}_k)\Big{]}^T(\alpha_k-\mathbb{E}_t\alpha_k)+\Big{[}\mathcal{Q}_{t,k}^1\mathbb{E}_tX^{*}_k\nonumber\\
&&\hspace{-2.5em}\hphantom{=}+\big{(}\mathcal{S}^{1(1)}_{t,k}\big{)}^T\mathbb{E}_tu^*_k +\big{(}\mathcal{S}^{1(2)}_{t,k}\big{)}^T\mathbb{E}_tv^*_k +q_{t,k}^1+A_{k}^T\mathbb{E}_tY^{*}_{k+1}+\sum_{i=1}^p(C_{k}^i)^T\mathbb{E}_t(Y^{*}_{k+1}w^i_k)-\mathbb{E}_tY^{*}_k\Big{]}^T\mathbb{E}_t\alpha_k \nonumber\\
&&\hspace{-2.5em}\hphantom{=}+\Big{[}S^{1(1)}_{t,k}(X^{*}_k-\mathbb{E}_tX^{*}_k)+R^{1(11)}_{t,k}(u^*_k-\mathbb{E}_tu^*_k) +R^{1(12)}_{t,k}(v^*_k-\mathbb{E}_tv^*_k) +(B^1_{k})^T(Y^{*}_{k+1}-\mathbb{E}_tY^{*}_{k+1})\nonumber\\
&&\hspace{-2.5em}\hphantom{=}+\sum_{i=1}^p(D_{k}^{1i})^T\Big{(}\mathbb{E}_k(Y^{*}_{k+1}w^i_k) -\mathbb{E}_t(Y^{*}_{k+1}w^i_k)\Big{)}\Big{]}^T (u_k-\mathbb{E}_tu_k)
+\Big{[}\mathcal{S}_{t,k}^{1(1)}\mathbb{E}_tX^{*}_k+\mathcal{R}^{1(11)}_{t,k}\mathbb{E}_tu^*_k\nonumber\\
&&\hspace{-2.5em}\hphantom{=}+\mathcal{R}^{1(12)}_{t,k}\mathbb{E}_tv^*_k+\rho_{t,k}^{1(11)}+(B^1_{k})^T\mathbb{E}_tY^{*}_{k+1}
+\sum_{i=1}^p(D_{k}^{1i})^T\mathbb{E}_t(Y^{*}_{k+1}w^i_k)\Big{]}^T\mathbb{E}_tu_k\Bigg{\}}\nonumber \\
&&\hspace{-2.5em}=\mathbb{E}_t\sum_{t=k}^{N-1}\Bigg{\{}\Big{[}S^{1(1)}_{t,k}(X^{*}_k-\mathbb{E}_tX^{*}_k)+R^{1(11)}_{t,k} (u^*_k-\mathbb{E}_tu^*_k)+R^{1(12)}_{t,k}(v^*_k-\mathbb{E}_tv^*_k) +(B^1_{k})^T(Y^{*}_{k+1}-\mathbb{E}_tY^{*}_{k+1})\nonumber\\
&&\hspace{-2.5em}\hphantom{=}+\sum_{i=1}^p(D_{k}^{1i})^T\Big{(}\mathbb{E}_k(Y^{*}_{k+1}w^i_k) -\mathbb{E}_t(Y^{*}_{k+1}w^i_k)\Big{)}\Big{]}^T (u_k-\mathbb{E}_tu_k)
+\Big{[}\mathcal{S}_{t,k}^{1(1)}\mathbb{E}_tX^{*}_k+\mathcal{R}^{1(11)}_{t,k}\mathbb{E}_tu^*_k\nonumber\\
&&\hspace{-2.5em}\hphantom{=}+\mathcal{R}^{1(12)}_{t,k}\mathbb{E}_tv^*_k+\rho_{t,k}^{1(11)}+(B^1_{k})^T\mathbb{E}_tY^{*}_{k+1}
+\sum_{i=1}^p(D_{k}^{1i})^T\mathbb{E}_t(Y^{*}_{k+1}w^i_k)\Big{]}^T\mathbb{E}_tu_k\Bigg{\}}.
\end{eqnarray*}
Then, (\ref{Nece-suffic1-1}) becomes
\begin{eqnarray}\label{Nece-suffic1-2}
&&\hspace{-2em}J_1(t,y;u^*+\varepsilon u, v^*)-J_1(t,y; u^*, v^*)\nonumber\\[1mm]
&&\hspace{-2em}=2\varepsilon \mathbb{E}_t\sum_{k=t}^{N-1}\Bigg{\{}\Big{[}S^{1(1)}_{t,k}(X^{*}_k-\mathbb{E}_tX^{*}_k) +R^{1(11)}_{t,k}(u^*_k-\mathbb{E}_tu^*_k)+R^{1(12)}_{t,k}(v^*_k-\mathbb{E}_tv^*_k) +(B^1_{k})^T(Y^{*}_{k+1}-\mathbb{E}_tY^{*}_{k+1})\nonumber\\
&&\hspace{-2em}\hphantom{=}+\sum_{i=1}^p(D_{k}^{1i})^T\Big{(}\mathbb{E}_k(Y^{*}_{k+1}w^i_k)-\mathbb{E}_t(Y^{*}_{k+1}w^i_k)\Big{)}\Big{]}^T (u_k-\mathbb{E}_tu_k)
+\Big{[}\mathcal{S}_{t,k}^{1(1)}\mathbb{E}_tX^{*}_k+\mathcal{R}^{1(11)}_{t,k}\mathbb{E}_tu^*_k\nonumber\\
&&\hspace{-2em}\hphantom{=}+\mathcal{R}^{1(12)}_{t,k}\mathbb{E}_tv^*_k+\rho_{t,k}^{1(11)}+(B^1_{k})^T\mathbb{E}_tY^{*}_{k+1}
+\sum_{i=1}^p(D_{k}^{1i})^T\mathbb{E}_t(Y^{*}_{k+1}w^i_k)\Big{]}^T\mathbb{E}_tu_k\Bigg{\}}+\varepsilon^2 \widetilde{J}_1(t,0;u)\nonumber\\
&&\hspace{-2em}\geq 0.
\end{eqnarray}
As (\ref{Nece-suffic1-2}) holds for any $\varepsilon \in \mathbb{R}$ and any $u\in l^2_{\mathbb{F}}(\mathbb{T}_t, \mathbb{R}^{m_1})$, we must have%
\begin{eqnarray*}
\inf_{u\in l^2_\mathbb{F}(\mathbb{T}_t; \mathbb{R}^{m_1})} \widetilde{J}_1(t,0;u)\geq 0,~~~a.s.,
\end{eqnarray*}
and
\begin{eqnarray*}
\left\{
\begin{array}{l}
0=S^{1(1)}_{t,k}(X^{*}_k-\mathbb{E}_tX^{*}_k)+R^{1(11)}_{t,k}(u^*_k-\mathbb{E}_tu^*_k)+R^{1(12)}_{t,k}(v^*_k-\mathbb{E}_tv^*_k) +(B^1_{k})^T(\mathbb{E}_kY^{*}_{k+1}-\mathbb{E}_tY^{*}_{k+1})\nonumber\\[1mm]
\hphantom{0=}+\sum_{i=1}^p(D_{k}^{1i})^T\Big{(}\mathbb{E}_k(Y^{*}_{k+1}w^i_k)-\mathbb{E}_t(Y^{*}_{k+1}w^i_k)\Big{)},\\[1mm]
0=\mathcal{S}_{t,k}^{1(1)}\mathbb{E}_tX^{*}_k+\mathcal{R}^{1(11)}_{t,k}\mathbb{E}_tu^*_k+\mathcal{R}^{1(12)}_{t,k}\mathbb{E}_tv^*_k +\rho_{t,k}^{1(11)}+(B^1_{k})^T\mathbb{E}_tY^{*}_{k+1}+\sum_{i=1}^p(D_{k}^{1i})^T\mathbb{E}_t(Y^{*}_{k+1}w^i_k),\\[1mm]
k\in \mathbb{T}_t,
\end{array}
\right.
\end{eqnarray*}
which implies the first equation of (\ref{stationary-condition}).

On the other hand, for any $\lambda\in \mathbb{R}$ and $v_k\in l^2_{\mathbb{F}}(k, \mathbb{R}^{m_2})$, let $X^{\lambda}$ satisfy the S$\Delta$E,
\begin{eqnarray}\label{system-6}
\left\{\begin{array}{l}
X^{\lambda}_{\ell+1}=\big{(}A_{\ell}X^{\lambda}_\ell+B^1_{\ell}u^*_\ell+{B}^2_{\ell}v^*_\ell\big{)}
+\sum_{i=1}^p\big{(}C^i_{\ell}X^{\lambda}_\ell+D^{1i}_{\ell}u^*_\ell+{D}^{2i}_{\ell}v^*_\ell\big{)}w^i_\ell, \\[2mm]
X^{\lambda}_{k+1}=\big{(}A_{k}X^{\lambda}_k+B^1_{k}u^*_k+{B}^2_{k}(v^*_k+\lambda v_k)\big{)}\\[2mm]
\hphantom{X^{\lambda}_{k+1}=}+\sum_{i=1}^p\big{(}C^i_{k}X^{\lambda}_k+D^{1i}_{k}u^*_k+{D}^{2i}_{k}(v^*_k+\lambda v_k)\big{)}w^i_k, \\[2mm]
X^{\lambda}_k=y,~~\ell\in  \mathbb{T}_{k+1}.
\end{array}
\right.
\end{eqnarray}
From (\ref{System-equilibrium}) and (\ref{system-6}), we have
\begin{eqnarray*}
\left\{
\begin{array}{l}
\frac{X^{\lambda}_{\ell+1}-X^{*}_{\ell+1}}{\lambda}=A_{\ell}\frac{X^{\lambda}_\ell-X^{*}_\ell}{\lambda}
+\sum_{i=1}^pC^i_{\ell}\frac{X^{\lambda}_\ell-X^{*}_\ell}{\lambda}w^i_\ell,\\
\frac{X^{\lambda}_{k+1}-X^{*}_{k+1}}{\lambda}=B_{k}^2v_k+\sum_{i=1}^pD_{k}^{2i}v_kw^i_k,\\
\frac{X^{\lambda}_k-X^{*}_k}{\lambda}=0,~~\ell\in  \mathbb{T}_{k+1}.
\end{array}
\right.
\end{eqnarray*}
Denote $\frac{X^{\lambda}_\ell-X^{*}_\ell}{\lambda}$ by $\beta_\ell$, then $\beta=\{\beta_\ell, ~~\ell\in  \mathbb{T}_k\}$ satisfies (\ref{system-4}), and $X^{\lambda}_\ell=X^{*}_\ell+\lambda\beta_\ell, \ell\in  \mathbb{T}_k$. Therefore, it holds that
\begin{eqnarray}\label{Nece-suffic1-3}
&&\hspace{-2.15em}J_2(k,X^{*}_k; u^*|_{\mathbb{T}_k},( v^*_k+\lambda v_k,v^*|_{\mathbb{T}_{k+1}}))-J_2(k,X^{*}_k; u^*|_{\mathbb{T}_k}, v^*|_{\mathbb{T}_k})\nonumber\\[1mm]
&&\hspace{-2.15em}=2\lambda\mathbb{E}_k\Bigg{\{}\sum_{\ell=k}^{N-1}\Big{[}\Big{(}Q_{k,\ell}^2X^{*}_\ell +\big{(}S^{2(1)}_{k,\ell}\big{)}^Tu^*_\ell +\big{(}S^{2(2)}_{k,\ell}\big{)}^Tv^*_\ell+q_{k,\ell}^2\Big{)}^T\beta_\ell\nonumber\\
&&\hspace{-2.15em}\hphantom{=}
+\Big{(}\bar Q_{k,\ell}^2\mathbb{E}_kX^{*}_\ell+\big{(}\bar S^{2(1)}_{k,\ell}\big{)}^T\mathbb{E}_ku^*_\ell+\big{(}\bar S^{2(2)}_{k,\ell}\big{)}^T\mathbb{E}_kv^*_\ell\Big{)}^T\mathbb{E}_k\beta_\ell\Big{]}\nonumber\\
&&\hspace{-2.15em}\hphantom{=}+\Big{(}\mathcal{S}^{2(2)}_{k,k}X^{*}_k+\mathcal{R}^{2(21)}_{k,k}u^*_k+\mathcal{R}^{2(22)}_{k,k}v^*_k +\rho^{2(2)}_{k,k}\Big{)}^Tv_k\nonumber\\
&&\hspace{-2.15em}\hphantom{=}+(G^2_kX^{*}_N+g^2_k)^T\beta_N +(\mathbb{E}_kX^{*}_N)^T
\bar{G}^2_k(\mathbb{E}_k\beta_N)\nonumber\Bigg{\}}\nonumber\\
&&\hspace{-2.15em}\hphantom{=}+\lambda^2\Bigg{\{}\sum_{\ell=k}^{N-1}\mathbb{E}_k\Big{[}\beta_\ell^TQ_{k,\ell}^2\beta_\ell +(\mathbb{E}_k\beta_\ell)^T\bar Q_{k,\ell}^2\mathbb{E}_k\beta_\ell\Big{]}\nonumber\\
&&\hspace{-2.15em}\hphantom{=}+\mathbb{E}_k[\beta_N^TG^2_k\beta_N]+(\mathbb{E}_k\beta_N)^T\bar{G}^2_k\mathbb{E}_k\beta_N +v_k^T\mathcal{R}^{2(22)}_{k,k}v_k\Bigg{\}}\nonumber\\
&&\hspace{-2.15em}\geq 0.
\end{eqnarray}
From (\ref{system-z-k}) and (\ref{system-4}), we have
\begin{eqnarray*}
&&\hspace{-2em}\mathbb{E}_k\Bigg{\{}\sum_{\ell=k}^{N-1}\Big{[}\Big{(}Q_{k,\ell}^2X^{*}_\ell +\big{(}S^{2(1)}_{k,\ell}\big{)}^Tu^*_\ell +\big{(}S^{2(2)}_{k,\ell}\big{)}^Tv^*_\ell+q_{k,\ell}^2\Big{)}^T\beta_\ell+\Big{(}\bar Q_{k,\ell}^2\mathbb{E}_kX^{*}_\ell+\big{(}\bar S^{2(1)}_{k,\ell}\big{)}^T\mathbb{E}_ku^*_\ell\nonumber\\
&&\hspace{-2em} \hphantom{=}+\big{(}\bar S^{2(2)}_{k,\ell}\big{)}^T\mathbb{E}_kv^*_\ell\Big{)}^T\mathbb{E}_k\beta_\ell\Big{]} +\Big{(}\mathcal{S}^{2(2)}_{k,k}X^{*}_k+\mathcal{R}^{2(21)}_{k,k}u^*_k+\mathcal{R}^{2(22)}_{k,k}v^*_k +\rho^{2(2)}_{k,k}\Big{)}^Tv_k\nonumber\\
&&\hspace{-2em} \hphantom{=}+(G^2_kX^{*}_N+g^2_k)^T\beta_N+(\mathbb{E}_kX^{*}_N)^T
\bar{G}^2_k\mathbb{E}_k\beta_N\nonumber\Bigg{\}}\nonumber\\
&&\hspace{-2em}=\mathbb{E}_k\Bigg{\{}\sum_{\ell=k}^{N-1}\Big{[}Q_{k,\ell}^2(X^{*}_\ell-\mathbb{E}_kX^{*}_\ell) +\big{(}S^{2(1)}_{k,\ell}\big{)}^T(u^*_\ell-\mathbb{E}_ku^*_\ell)+\big{(}S^{2(2)}_{k,\ell}\big{)}^T(v^*_\ell-\mathbb{E}_kv^*_\ell)+A_{\ell}^T(\mathbb{E}_\ell Z^{k,*}_{\ell+1}-\mathbb{E}_kZ^{k,*}_{\ell+1})\\
&&\hspace{-2em} \hphantom{=}+\sum_{i=1}^p(C_{\ell}^i)^T\Big{(}\mathbb{E}_\ell(Z^{k,*}_{\ell+1}w^i_\ell)-\mathbb{E}_k(Z^{k,*}_{\ell+1} w^i_\ell)\Big{)} -(Z^{k,*}_\ell-\mathbb{E}_kZ^{k,*}_\ell)\Big{]}^T(\beta_\ell-\mathbb{E}_k\beta_\ell)+\Big{[}\mathcal{Q}_{k,\ell}^2\mathbb{E}_kX^{*}_\ell\nonumber\\
&&\hspace{-2em} \hphantom{=} +\big{(}\mathcal{S}^{2(1)}_{k,\ell}\big{)}^T\mathbb{E}_ku^*_\ell+q_{k,\ell}^2 +\big{(}\mathcal{S}^{2(2)}_{k,\ell}\big{)}^T\mathbb{E}_kv^*_\ell +A_{\ell}^T\mathbb{E}_kZ^{k,*}_{\ell+1}+\sum_{i=1}^p(C_{\ell}^i)^T\mathbb{E}_k(Z^{k,*}_{\ell+1}w^i_\ell)-\mathbb{E}_kZ^{k,*}_\ell\Big{]}^T \mathbb{E}_k\beta_\ell\Bigg{\}}\\
&&\hspace{-2em} \hphantom{=}
+\Big{(}\mathcal{S}^{2(2)}_{k,k}X^{*}_k+\mathcal{R}^{2(21)}_{k,k}u^*_k+\mathcal{R}^{2(22)}_{k,k}v^*_k +(B_{k}^2)^T\mathbb{E}_kZ^{k,*}_{k+1}+\sum_{i=1}^p(D_{k}^{2i})^T\mathbb{E}_k(Z^{k,*}_{k+1}w^i_k)+\rho^{2(2)}_{k,k}\Big{)}^Tv_k\\
&&\hspace{-2em}=\Big{(}\mathcal{S}^{2(2)}_{k,k}X^{*}_k+\mathcal{R}^{2(21)}_{k,k}u^*_k+\mathcal{R}^{2(22)}_{k,k}v^*_k +(B_{k}^2)^T\mathbb{E}_kZ^{k,*}_{k+1}+\sum_{i=1}^p(D_{k}^{2i})^T\mathbb{E}_k(Z^{k,*}_{k+1}w^i_k)+\rho^{2(2)}_{k,k}\Big{)}^Tv_k.
\end{eqnarray*}
Hence, (\ref{Nece-suffic1-3}) becomes
\begin{eqnarray*}
&&\hspace{-3em}J_2(k,X^{*}_k; u^*|_{\mathbb{T}_k},( v^*_k+\lambda v_k,v^*|_{\mathbb{T}_{k+1}}))-J_2(k,X^{*}_k; u^*|_{\mathbb{T}_k}, v^*|_{\mathbb{T}_k})\nonumber\\[1mm]
&&\hspace{-3em}= 2\lambda \Big{(}\mathcal{S}^{2(2)}_{k,k}X^{*}_k+\mathcal{R}^{2(21)}_{k,k}u^*_k+\mathcal{R}^{2(22)}_{k,k}v^*_k+(B_{k}^2)^T\mathbb{E}_kZ^{k,*}_{k+1}\nonumber\\
&&\hspace{-3em} \hphantom{=}+\sum_{i=1}^p(D_{k}^{2i})^T\mathbb{E}_k(Z^{k,*}_{k+1}w^i_k)+\rho^{2(2)}_{k,k}\Big{)}^Tv_k+\lambda^2 \widetilde{J}_2(k,0;v_k)\\
&&\hspace{-3em}\geq 0,
\end{eqnarray*}
which holds for any $\lambda\in \mathbb{R}$ and any $v_k\in l^2_{\mathbb{F}}(k, \mathbb{R}^{m_2})$. Therefore,
\begin{eqnarray*}
\inf_{v\in l^2_\mathbb{F}(\mathbb{T}_k; \mathbb{R}^{m_2})} \widetilde{J}_2(k,0;v_k)\geq 0,~~~a.s.,
\end{eqnarray*}
and
the second equation of (\ref{stationary-condition}) holds.

\emph{ii)$\Rightarrow$i)}. By reversing the proof of i)$\Rightarrow$ii), we can obtain the conclusion. \hfill $\square$

\subsection{Proof of Theorem \ref{Theorem-Problem-GLQ}}

\begin{proposition}\label{Thm-stationary-condition}
The following statements are equivalent.

i). There exists a $({u}^*,v^*)\in l^2_\mathbb{F}(\mathbb{T}_t; \mathbb{R}^{m_1})\times l^2_\mathbb{F}(\mathbb{T}_t; \mathbb{R}^{m_2})$ such that the stationary conditions of (\ref{stationary-condition}) hold.

ii). a) of Theorem \ref{Theorem-Problem-GLQ} is satisfied.

Under the condition ii), the backward states $Y^{*}, Z^{k,*}$ of (\ref{system-y-k}) and (\ref{system-z-k}) have the following expressions
\begin{eqnarray}
\left\{
\begin{array}{l}
Y_{k}^{*}=P_{t,k}(X_{k}^*-\mathbb{E}_tX_{k}^*)+{\mathcal{P}}_{t,k}\mathbb{E}_tX_{k}^*+\sigma_{t,k},~~~~~k\in \mathbb{T}_t,\\[1mm]
Z_{\ell}^{k,*}=T_{k,\ell}(X_{\ell}^{*}-\mathbb{E}_{k}X_{\ell}^{*})+{\mathcal{T}}_{k,\ell}\mathbb{E}_{k}X_{\ell}^{*}+\widetilde{T}_{k,\ell}\mathbb{E}_tX_{\ell}^{*}+\xi_{k,\ell},~~~k\in \mathbb{T}_t,~~~\ell\in \mathbb{T}_k.
\end{array}
\right.
\end{eqnarray}

\end{proposition}

\emph{Proof.} \emph{i)$\Rightarrow$ii)}. From (\ref{stationary-condition}), it holds that
\begin{eqnarray}\label{stationary-condition-(N-1)-1}
\left\{
\begin{array}{l}
0=\mathcal{{S}}^{1(1)}_{t,k}\mathbb{E}_tX^{*}_k+\mathcal{R}^{1(11)}_{t,k}\mathbb{E}_tu^*_k +\mathcal{R}^{1(12)}_{t,k}\mathbb{E}_tv^*_k+(B^1_{k})^T\mathbb{E}_tY^{*}_{k+1}\\[1mm]
\hphantom{0=}+\sum_{i=1}^p(D_{k}^{1i})^T\mathbb{E}_t(Y^{*}_{k+1}w^i_k) +\rho_{t,k}^{1(1)},~~~~~~k\in  \mathbb{T}_t,\\[2mm]
0=\mathcal{S}^{2(2)}_{k,k}\mathbb{E}_tX^{*}_k+\mathcal{R}^{2(21)}_{k,k}\mathbb{E}_tu^*_k +\mathcal{R}^{2(22)}_{k,k}\mathbb{E}_tv^*_k +(B_{k}^2)^T\mathbb{E}_tZ^{k,*}_{k+1}\\[1mm]
\hphantom{0=}+\sum_{i=1}^p(D_{k}^{2i})^T\mathbb{E}_t(Z^{k,*}_{k+1}w^i_k)+\rho^{2(2)}_{k,k},~~~~~~k\in  \mathbb{T}_t,
\end{array}
\right.
\end{eqnarray}
and
\begin{eqnarray}\label{stationary-condition-(N-1)-2}
\left\{
\begin{array}{l}
0=S^{1(1)}_{t,k}\big{(}X_k^{*}-\mathbb{E}_tX_k^{*}\big{)}+R^{1(11)}_{t,k}\big{(}u_k^*-\mathbb{E}_tu_k^*\big{)} +R^{1(12)}_{t,k}\big{(}v_k^*-\mathbb{E}_tv_k^*\big{)}\\[1mm]
\hphantom{0=}+(B^1_{k})^T\big{(}\mathbb{E}_kY_{k+1}^{*}-\mathbb{E}_tY_{k+1}^{*}\big{)}+\sum_{i=1}^p(D^{1i}_{k})^T\big{(} \mathbb{E}_k(Y_{k+1}^{*}w_{k}^i))-\mathbb{E}_t(Y_{k+1}^{*}w_{k}^i)\big{)},~~~~k\in  \mathbb{T}_t,\\[2mm]
0=\mathcal{S}^{2(2)}_{k,k}(X^{*}_k-\mathbb{E}_tX^{*}_k)+\mathcal{R}^{2(21)}_{k,k}(u^*_k-\mathbb{E}_tu^*_k) +\mathcal{R}^{2(22)}_{k,k}(v^*_k-\mathbb{E}_tv^*_k) \\[1mm]
\hphantom{0=}+(B_{k}^2)^T(\mathbb{E}_kZ^{k,*}_{k+1}-\mathbb{E}_tZ^{k,*}_{k+1})+\sum_{i=1}^p(D_{k}^{2i})^T \big{(}\mathbb{E}_k{(}Z^{k,*}_{k+1}w^i_k)-\mathbb{E}_t(Z^{k,*}_{k+1}w^i_k)\big{)},~~~~k\in  \mathbb{T}_t.
\end{array}
\right.
\end{eqnarray}

Let us first consider the case $k=N-1$. We have
%
%
%
%
%
%
%
\begin{eqnarray}
\label{stationary-condition-(N-1)-0}
\begin{array}{l}
\mathbb{E}_{t}Y_{N}^{*}=\mathcal{G}_t^1A_{N-1}\mathbb{E}_tX_{N-1}^{*}+\mathcal{G}_t^1B^1_{N-1}\mathbb{E}_tu^*_{N-1} +\mathcal{G}_t^1B^2_{N-1}\mathbb{E}_tv^*_{N-1}+g_t^1,\\
\mathbb{E}_{t}{(}Y_{N}^{*}w_{N-1}^i{)}=G_t^1\sum_{j=1}^p\delta_{N-1}^{ij}\big{(}C^j_{N-1}\mathbb{E}_tX_{N-1}^{*} +D^{1j}_{N-1}\mathbb{E}_tu_{N-1}^*+D^{2j}_{N-1}\mathbb{E}_tv_{N-1}^*\big{)}.
\end{array}
\end{eqnarray}
Then, the first equations of (\ref{stationary-condition-(N-1)-1})and (\ref{stationary-condition-(N-1)-2}) become
\begin{eqnarray*}\label{stationary-condition-(N-1)-1-1}
\begin{array}{l}
0=\big{(}\mathcal{S}_{t,N-1}^{1(1)}+(B_{N-1}^1)^T\mathcal{G}_t^1A_{N-1}+\sum_{i,j=1}^p\delta_{N-1}^{ij} (D^{1i}_{N-1})^TG_t^1C_{N-1}^j \big{)}\mathbb{E}_tX_{N-1}^{*}\\[1mm]
\hphantom{0=}+\big{(}{\mathcal{R}}_{t,N-1}^{1(11)}+(B_{N-1}^1)^T\mathcal{G}_t^1B_{N-1}^1+\sum_{i,j=1}^p\delta_{N-1}^{ij} (D^{1i}_{N-1})^TG_t^1D_{N-1}^{1j} \big{)} \mathbb{E}_tu_{N-1}^*\\[1mm]
\hphantom{0=}+\big{(}{\mathcal{R}}_{t,N-1}^{1(12)}+(B_{N-1}^1)^T\mathcal{G}_t^1B_{N-1}^2+\sum_{i,j=1}^p\delta_{N-1}^{ij} (D^{1i}_{N-1})^TG_t^1D_{N-1}^{2j} \big{)} \mathbb{E}_tv_{N-1}^*\\[1mm]
\hphantom{0=}+(B_{N-1}^1)^Tg_t^1+\rho_{t,N-1}^{1(1)},
\end{array}
\end{eqnarray*}
and
\begin{eqnarray*}\label{stationary-condition-(N-1)-2-1}
\begin{array}{l}
\hspace{-2em}0=\big{(}{S}_{t,N-1}^{1(1)}+(B_{N-1}^1)^T{G}_t^1A_{N-1}+\sum_{i,j=1}^p\delta_{N-1}^{ij} (D^{1i}_{N-1})^TG_t^1C_{N-1}^j \Big{)}\big{(}X^{*}_{N-1}-\mathbb{E}_tX_{N-1}^{*}\big{)}\\
\hspace{-2em}\hphantom{0=}+\Big{(}{{R}}_{t,N-1}^{1(11)}+(B_{N-1}^1)^T{G}_t^1B_{N-1}^1+\sum_{i,j=1}^p\delta_{N-1}^{ij} (D^{1i}_{N-1})^TG_t^1D_{N-1}^{1j} \Big{)}\big{(}u_{N-1}^*-\mathbb{E}_tu_{N-1}^*\big{)}\\
\hspace{-2em}\hphantom{0=}+\Big{(}{{R}}_{t,N-1}^{1(12)}+(B_{N-1}^1)^T{G}_t^1B_{N-1}^2+\sum_{i,j=1}^p\delta_{N-1}^{ij} (D^{1i}_{N-1})^TG_t^1D_{N-1}^{2j} \Big{)}\big{(}v_{N-1}^*- \mathbb{E}_tv_{N-1}^*\big{)}.
\end{array}
\end{eqnarray*}
Furthermore,
%
%
%
%
%
%
\begin{eqnarray*}
&&\mathbb{E}_{t}Z_{N}^{N-1,*}=\mathcal{G}_{N-1}^2A_{N-1}\mathbb{E}_tX_{N-1}^{*} +\mathcal{G}_{N-1}^2B^1_{N-1}\mathbb{E}_tu^*_{N-1} +\mathcal{G}_{N-1}^2B^2_{N-1}\mathbb{E}_tv^*_{N-1}+g_{N-1}^2,\\
&&\mathbb{E}_{t}{(}Z_{N}^{N-1,*}w_{N-1}^i{)}=G_{N-1}^2\sum_{j=1}^p\delta_{N-1}^{ij}\big{(}C^j_{N-1}\mathbb{E}_tX_{N-1}^{*} +D^{1j}_{N-1}\mathbb{E}_tu_{N-1}^*+D^{2j}_{N-1}\mathbb{E}_tv_{N-1}^*\big{)}.
\end{eqnarray*}
Therefore, the second equations of (\ref{stationary-condition-(N-1)-1}) and (\ref{stationary-condition-(N-1)-2}) become
\begin{eqnarray}\label{stationary-condition-(N-1)-1-2}
\begin{array}{l}
\hspace{-1.75em}0=\big{(}\mathcal{S}_{N-1,N-1}^{2(2)}+(B_{N-1}^2)^T\mathcal{G}_{N-1}^2A_{N-1}+\sum_{i,j=1}^p\delta_{N-1}^{ij} (D^{2i}_{N-1})^TG_{N-1}^2C_{N-1}^j \big{)}\mathbb{E}_tX_{N-1}^{*}\\
\hspace{-1.75em}\hphantom{0=}+\big{(}{\mathcal{R}}_{{N-1},N-1}^{2(21)}+(B_{N-1}^2)^T\mathcal{G}_{N-1}^2B_{N-1}^1 +\sum_{i,j=1}^p\delta_{N-1}^{ij} (D^{2i}_{N-1})^TG_{N-1}^2D_{N-1}^{1j} \big{)} \mathbb{E}_tu_{N-1}^*\\
\hspace{-1.75em}\hphantom{0=}+\big{(}{\mathcal{R}}_{{N-1},N-1}^{2(22)}+(B_{N-1}^2)^T\mathcal{G}_{N-1}^2B_{N-1}^2 +\sum_{i,j=1}^p\delta_{N-1}^{ij} (D^{2i}_{N-1})^TG_{N-1}^2D_{N-1}^{2j} \big{)} \mathbb{E}_tv_{N-1}^* \\
\hspace{-1.75em}\hphantom{0=}+(B_{N-1}^2)^Tg_{N-1}^2+\rho_{{N-1},N-1}^{2(2)}.
\end{array}
\end{eqnarray}
and
\begin{eqnarray}\label{stationary-condition-(N-1)-2-2}
\begin{array}{l}
\hspace{-0.8em}0=\Big{(}\mathcal{S}_{N-1,N-1}^{2(2)}+(B_{N-1}^2)^T\mathcal{G}_{N-1}^2A_{N-1}+\sum_{i,j=1}^p\delta_{N-1}^{ij} (D^{2i}_{N-1})^TG_{N-1}^2C_{N-1}^j \Big{)}(X _{N-1}^{*}-\mathbb{E}_tX_{N-1}^{*})\\
\hspace{-0.8em}\hphantom{0=}+\Big{(}{\mathcal{R}}_{{N-1},N-1}^{2(21)}+(B_{N-1}^2)^T\mathcal{G}_{N-1}^2B_{N-1}^1 +\sum_{i,j=1}^p\delta_{N-1}^{ij} (D^{2i}_{N-1})^TG_{N-1}^2D_{N-1}^{1j} \Big{)} (u_{N-1}^*-\mathbb{E}_tu_{N-1}^*)\\
\hspace{-0.8em}\hphantom{0=}+\Big{(}{\mathcal{R}}_{{N-1},N-1}^{2(22)}+(B_{N-1}^2)^T\mathcal{G}_{N-1}^2B_{N-1}^2 +\sum_{i,j=1}^p\delta_{N-1}^{ij} (D^{2i}_{N-1})^TG_{N-1}^2D_{N-1}^{2j} \Big{)} (v_{N-1}^*-\mathbb{E}_tv_{N-1}^*).
\end{array}
\end{eqnarray}
With the notations of (\ref{W-H}) (\ref{W-bf}), we have from above equations
\begin{eqnarray}\label{Propo-4-1-1}
&&\hspace{-3em}0=\widetilde{\mathbf{H}}_{t,N-1}\left(
\begin{array}{c}
\mathbb{E}_tX_{N-1}^{*}\\\mathbb{E}_tX_{N-1}^{*}
\end{array}\right)+\widetilde{\mathbf{W}}_{t,N-1}\left(
\begin{array}{c}
\mathbb{E}_tu_{N-1}^*\\\mathbb{E}_tv_{N-1}^*
\end{array}\right)
+\mathbf{h}_{t,N-1}
\end{eqnarray}
and
\begin{eqnarray}\label{Propo-4-1-2}
&&0=\mathbf{H}_{t,N-1}\left(
\begin{array}{c}
X _{N-1}^{*}-\mathbb{E}_tX_{N-1}^{*}\\X _{N-1}^{*}-\mathbb{E}_tX_{N-1}^{*}
\end{array}\right)+\mathbf{W}_{t,N-1}\left(
\begin{array}{c}
u_{N-1}^{*}-\mathbb{E}_tu_{N-1}^{*}\\v_{N-1}^{*}-\mathbb{E}_tv_{N-1}^{*}
\end{array}\right).
\end{eqnarray}
Therefore, by a property of Moore-Penrose inverse (Lemma 3.1 of \cite{Ait-Chen-Zhou-2002}), (\ref{stationary-condition-(1)}) (\ref{stationary-condition-(2)}) hold for $k=N-1$ and we can select
\begin{eqnarray*}
\left(
\begin{array}{c}
\mathbb{E}_tu_{N-1}^*\\\mathbb{E}_tv_{N-1}^*
\end{array}\right)=-\widetilde{\mathbf{W}}_{t,N-1}^\dagger\Bigg{[}\widetilde{\mathbf{H}}_{t,N-1}\left(
\begin{array}{c}
\mathbb{E}_tX_{N-1}^{*}\\\mathbb{E}_tX_{N-1}^{*}
\end{array}\right)+\mathbf{h}_{t,N-1}\Bigg{]},
\end{eqnarray*}
and
\begin{eqnarray*}
\left(
\begin{array}{l}
u_{N-1}^{*}-\mathbb{E}_tu_{N-1}^{*}\\v_{N-1}^{*}-\mathbb{E}_tv_{N-1}^{*}
\end{array}\right)=-\mathbf{W}_{t,N-1}^\dagger \mathbf{H}_{t,N-1}\left(
\begin{array}{c}
X _{N-1}^{*}-\mathbb{E}_tX_{N-1}^{*}\\X _{N-1}^{*}-\mathbb{E}_tX_{N-1}^{*}
\end{array}\right).
\end{eqnarray*}
Hence, (\ref{u*-v*}) holds for $k=N-1$.
Substituting $(u_{N-1}^*, v^*_{N-1})$ into (\ref{system-y-k}) (\ref{system-z-k}), we have
\begin{eqnarray*}
&&\hspace{-1.5em}Y_{N-1}^{*}=P_{t,N-1}(X_{N-1}^*-\mathbb{E}_tX_{N-1}^*)+{\mathcal{P}}_{t,N-1}\mathbb{E}_tX_{N-1}^*+\sigma_{t,N-1},\\
&&\hspace{-1.5em}Z_{N-1}^{r,*}=T_{r,N-1}(X_{N-1}^{*}-\mathbb{E}_{r}X_{N-1}^{*})+{\mathcal{T}}_{r,N-1}\mathbb{E}_{r}X_{N-1}^{*} +\widetilde{T}_{r,N-1}\mathbb{E}_tX_{N-1}^{*}+\xi_{r,N-1},~~\forall r \in \{t,...,N-2\}.
\end{eqnarray*}

For $k=N-2$ and by mimic the derivations between (\ref{stationary-condition-(N-1)-0}) and (\ref{stationary-condition-(N-1)-2-2}), we have
\begin{eqnarray}\label{Propo-4-1-3}
0=\widetilde{\mathbf{H}}_{t,N-2}\left(
\begin{array}{c}
\mathbb{E}_tX_{N-2}^{*}\\\mathbb{E}_tX_{N-2}^{*}
\end{array}\right)+\widetilde{\mathbf{W}}_{t,N-2}\left(
\begin{array}{c}
\mathbb{E}_tu_{N-2}^*\\\mathbb{E}_tv_{N-2}^*
\end{array}\right)+\mathbf{h}_{t,N-2}
\end{eqnarray}

and

\begin{eqnarray}\label{Propo-4-1-4}
0=\mathbf{H}_{t,N-2}\left(
\begin{array}{c}
X _{N-2}^{*}-\mathbb{E}_tX_{N-2}^{*}\\X _{N-2}^{*}-\mathbb{E}_tX_{N-2}^{*}
\end{array}\right)+\mathbf{W}_{t,N-2}\left(
\begin{array}{c}
u_{N-2}^{*}-\mathbb{E}_tu_{N-2}^{*}\\v_{N-2}^{*}-\mathbb{E}_tv_{N-2}^{*}
\end{array}\right).
\end{eqnarray}
Therefore, (\ref{stationary-condition-(1)}) (\ref{stationary-condition-(2)}) hold for $k=N-2$ and we can select
\begin{eqnarray*}
&&\hspace{-6em}\left(
\begin{array}{c}
\mathbb{E}_tu_{N-2}^*\\\mathbb{E}_tv_{N-2}^*
\end{array}\right)=-\widetilde{\mathbf{W}}_{t,N-2}^\dagger\Bigg{[}\widetilde{\mathbf{H}}_{t,N-2}\left(
\begin{array}{c}
\mathbb{E}_tX_{N-2}^{*}\\\mathbb{E}_tX_{N-2}^{*}
\end{array}\right)+\mathbf{h}_{t,N-2}\Bigg{]},
\end{eqnarray*}
and
\begin{eqnarray*}
&&\left(
\begin{array}{c}
u_{N-2}^{*}-\mathbb{E}_tu_{N-2}^{*}\\v_{N-2}^{*}-\mathbb{E}_tv_{N-2}^{*}
\end{array}\right)=-\mathbf{W}_{t,N-2}^\dagger \mathbf{H}_{t,N-2}\left(
\begin{array}{c}
X _{N-2}^{*}-\mathbb{E}_tX_{N-2}^{*}\\X _{N-2}^{*}-\mathbb{E}_tX_{N-2}^{*}
\end{array}\right).
\end{eqnarray*}
Hence, Hence, (\ref{u*-v*}) holds for $k=N-2$. Substituting $(u^*_{N-2}, v^*_{N-2})$ into (\ref{system-y-k}) (\ref{system-z-k}), we have
\begin{eqnarray*}
&&Y_{N-2}^{*}=P_{t,N-2}(X_{N-2}^*-\mathbb{E}_tX_{N-2}^*)+{\mathcal{P}}_{t,N-2}\mathbb{E}_tX_{N-2}^*+\sigma_{t,N-2},\\
&&Z_{N-2}^{r,*}=T_{r,N-2}(X_{N-2}^{*}-\mathbb{E}_{r}X_{N-2}^{*})+{\mathcal{T}}_{r,N-2}\mathbb{E}_{r}X_{N-2}^{*} \\
&&\hphantom{Z_{N-2}^{r,*}=}+\widetilde{T}_{r,N-1}\mathbb{E}_tX_{N-2}^{*}+\xi_{r,N-2},~~~~~\forall r\in\{t,...,N-3\}.
\end{eqnarray*}
By repeating the above procedure and the arguments of induction, we can get the expressions of the backward states $Y^{*}, Z^{k,*}$ and $u^*, v^*$.

\emph{ii)$\Rightarrow$i)}. Due to the property of Moore-Penrose inverse and by reversing the proof of i)$\Rightarrow$ii), we can obtain the conclusion. \hfill $\square$

We now study the convex condition (\ref{convex}). By adding to and subtracting
\begin{eqnarray*}
\sum_{k=t}^{N-1}\mathbb{E}_t\Big{[}(\alpha_{k+1})^TU_{t,k+1}\alpha_{k+1}-(\alpha_k)^TU_{t,k}\alpha_k +(\mathbb{E}_t\alpha_{k+1})^T\bar{U}_{t,k+1}\mathbb{E}_t\alpha_{k+1}
-(\mathbb{E}_t\alpha_k)^T\bar{U}_{t,k}\mathbb{E}_t\alpha_k
\Big{]}
\end{eqnarray*}
from $\widetilde{J}_1(t,0;u)$ (with $\bar{U}_{t,k}=\mathcal{U}_{t,k}-U_{t,k}, k\in \mathbb{T}_t$), we have
\begin{eqnarray}\label{convex-u}
&&\hspace{-4em}\widetilde{J}_1(t,0;u)=\sum_{k=t}^{N-1}\mathbb{E}_t\Big{[}(\alpha_k-\mathbb{E}_t\alpha_k)^T M_{t,k}^TO_{t,k}^\dagger M_{t,k}(\alpha_k-\mathbb{E}_t\alpha_k)+2(u_k-\mathbb{E}_tu_k)^T M_{t,k}(\alpha_k-\mathbb{E}_t\alpha_k)\nonumber\\
&&\hspace{-4em}\hphantom{\widetilde{J}_1(t,0;u)=}+(u_k-\mathbb{E}_tu_k)^TO_{t,k}(u_k-\mathbb{E}_tu_k)
+(\mathbb{E}_t\alpha_k)^T\mathcal{M}_{t,k}^T\mathcal{O}_{t,k}^\dagger \mathcal{M}_{t,k}\mathbb{E}_t\alpha_k \nonumber\\
&&\hspace{-4em}\hphantom{\widetilde{J}_1(t,0;u)=}+2(\mathbb{E}_tu_k)^T\mathcal{M}_{t,k}\mathbb{E}_t\alpha_k +(\mathbb{E}_tu_k)^T\mathcal{O}_{t,k}\mathbb{E}_tu_k\Big{]}.
\end{eqnarray}
Similarly,
\begin{eqnarray}\label{convex-v}
&&\widetilde{J}_2(k,0;v_k)=v_k^T\mathbb{O}_{k,k}v_k.
\end{eqnarray}
As $O_{t,k}$, $\mathcal{O}_{t,k}$ are symmetric, there exist orthogonal matrices $F_{t,k}$, $\mathcal{F}_{t,k}$ such that
\begin{eqnarray*}
O_{t,k}=(F_{t,k})^T\left(
\begin{array}{cc}
\Sigma_{t,k}&0\\0&0
\end{array}\right)F_{t,k},\\
\mathcal{O}_{t,k}=(\mathcal{F}_{t,k})^T\left(
\begin{array}{cc}
\Gamma_{t,k}&0\\0&0
\end{array}\right)\mathcal{F}_{t,k}.
\end{eqnarray*}
In the above, $\Sigma_{t,k}$, $\Gamma_{t,k}$, are diagonal matrices, whose diagonal elements are the nonzero eigenvalues of $O_{t,k}$,$\mathcal{O}_{t,k}$, respectively. Let $\mbox{rank}(O_{t,k})=r_k^1$, $\mbox{rank}(\mathcal{O}_{t,k})=r_k^2$. Then, we have
\begin{eqnarray*}
O_{t,k}^\dagger =(F_{t,k})^T\left(
\begin{array}{cc}
\Sigma_{t,k}^{-1}&0\\0&0
\end{array}\right) F_{t,k},\\
\mathcal{O}_{t,k}^\dagger =(\mathcal{F}_{t,k})^T\left(
\begin{array}{cc}
\Gamma_{t,k}^{-1}&0\\0&0
\end{array}\right)\mathcal{F}_{t,k}.
\end{eqnarray*}
Moreover, $F_{t,k}$, $\mathcal{F}_{t,k}$ can be decomposed as $F_{t,k}=[(F_{t,k}^{(1)})^T, (F_{t,k}^{(2)})^T]^T$, $\mathcal{F}_{t,k}=[(\mathcal{F}_{t,k}^{(1)})^T, (\mathcal{F}_{t,k}^{(2)})^T]^T$, respectively, where the lines of $F_{t,k}^{(2)}$, $\mathcal{F}_{t,k}^{(2)}$ form the bases of $\mbox{Ker}(O_{t,k})$ and $\mbox{Ker}(\mathcal{O}_{t,k})$, respectively. Let
\begin{eqnarray*}
&&O_{t,k}u_k=\left(
\begin{array}{c}
F_{t,k}^{(1)}u_k\\F_{t,k}^{(2)}u_k
\end{array}\right),~~~~
\mathcal{O}_{t,k}u_k=\left(
\begin{array}{c}
\mathcal{F}_{t,k}^{(1)}u_k\\ \mathcal{F}_{t,k}^{(2)}u_k
\end{array}\right).
%
%
%
%
\end{eqnarray*}
Hence, we have
\begin{eqnarray}\label{convex-u-1}
&&\hspace{-3.8em}\widetilde{J}_1(t,0;u)=\sum_{k=t}^{N-1}\mathbb{E}_t\Big{\{}\Big{[}{F}_{t,k}^{(1)} (u_k-\mathbb{E}_tu_k)+{\Sigma}_{t,k}^{-1}{F}_{t,k}^{(1)}M_{t,k}(\alpha_k-\mathbb{E}_t\alpha_k)
\Big{]}^T{\Sigma}_{t,k}\Big{[}{F}_{t,k}^{(1)}(u_k-\mathbb{E}_tu_k)\nonumber\\
&&\hspace{-3.8em}\hphantom{\widetilde{J}_1(t,0;u)=}+{\Sigma}_{t,k}^{-1}{F}_{t,k}^{(1)}M_{t,k} (\alpha_k-\mathbb{E}_t\alpha_k)\Big{]}+\Big{[}\mathcal{F}_{t,k}^{(1)}\mathbb{E}_tu_k+\Gamma_{t,k}^{-1}\mathcal{F}_{t,k}^{(1)}
\mathcal{M}_{t,k}\mathbb{E}_t\alpha_k\Big{]}^T
\Gamma_{t,k}\Big{[}\mathcal{F}_{t,k}^{(1)}\mathbb{E}_tu_k\nonumber\\
&&\hspace{-3.8em}\hphantom{\widetilde{J}_1(t,0;u)=}+\Gamma_{t,k}^{-1}\mathcal{F}_{t,k}^{(1)} \mathcal{M}_{t,k}\mathbb{E}_t\alpha_k\Big{]}\Big{\}}+2\sum_{k=t}^{N-1}\mathbb{E}_t\Big{[}\Big{(}{F}_{t,k}^{(2)}M_{t,k} (\alpha_k-\mathbb{E}_t\alpha_k)\Big{)}^T
{F}_{t,k}^{(2)}(u_k-\mathbb{E}_tu_k)\Big{]}\nonumber\\
&&\hspace{-3.8em}\hphantom{\widetilde{J}_1(t,0;u)=}
+2\sum_{k=t}^{N-1}\mathbb{E}_t\Big{[}\Big{(}\mathcal{F}_{t,k}^{(2)}\mathcal{M}_{t,k}\mathbb{E}_t\alpha_k\Big{)}^T \mathcal{F}_{t,k}^{(2)}\mathbb{E}_tu_k\Big{]}.
\end{eqnarray}
Note that the space spanned by lines of ${F}_{t,k}^{(1)}$ is $\mbox{Ran}(O_{t,k})$. Let
\begin{eqnarray*}
&&{\mathbb{U}}^{1}(\mbox{Ran})=\Big{\{}u\mid u\in l^2_\mathbb{F}(\mathbb{T}_t; \mathbb{R}^{m_1}), u_k-\mathbb{E}_tu_k\in \mbox{Ran}({O}_{t,k}),\mbox{ and }\mathbb{E}_tu_k=0,~k\in  \mathbb{T}_t\Big{\}},\\
&&{\mathbb{U}}^{1}(\mbox{Ker})=\Big{\{}u\mid u\in l^2_\mathbb{F}(\mathbb{T}_t; \mathbb{R}^{m_1}),  u_k-\mathbb{E}_tu_k\in \mbox{Ker}({O}_{t,k}),\mbox{ and }\mathbb{E}_tu_k=0,~k\in  \mathbb{T}_t\Big{\}}.
\end{eqnarray*}

\begin{proposition}\label{Thm-convexity-condition}
The following statements are equivalent.
\begin{itemize}
\item[i)] The convex conditions of (\ref{convex}) hold.

\item[ii)] b) and c) of Theorem \ref{Theorem-Problem-GLQ} are satisfied.
\end{itemize}
\end{proposition}

\emph{Proof.} \emph{i)$\Rightarrow$ii)}. Note that $u\mapsto \widetilde{J}_1(t,0; u)$ is convex. If $u\in{\mathbb{U}}^{1}(\mbox{Ran})$, we have $\mathbb{E}_t\alpha_k=0, k\in \mathbb{T}_t$, and
\begin{eqnarray}\label{convex-u-1-1}
\widetilde{J}_1(t,0;u)=\sum_{k=t}^{N-1}\mathbb{E}_t\Big{\{}\Big{[}{F}_{t,k}^{(1)}u_k+\Sigma_{t,k}^{-1} {F}_{t,k}^{(1)}M_{t,k}\alpha_k\Big{]}^T{\Sigma}_{t,k}[{F}_{t,k}^{(1)}u_k +{\Sigma}_{t,k}^{-1}{F}_{t,k}^{(1)}M_{t,k}\alpha_k\Big{]}\Big{\}}\geq0.
\end{eqnarray}
Introduce a set
\begin{eqnarray*}
\widetilde{\mathbb{U}}^{1}(\mbox{Ran})=\Big{\{}({F}_{t,t}^{(1)}u_t,\cdots,{F}_{t,N-1}^{(1)}u_{N-1})\big{|} u\in{\mathbb{U}}^{1}(\mbox{Ran})\Big{\}}.
\end{eqnarray*}
For $k\in  \mathbb{T}_t$, let $f_k^1, \cdots, f_k^{r_k^1}$ be the lines of ${F}_{t,k}^{(1)}$, then $(f_k^1)^T, \cdots, (f_k^{r_k^1})^T$ form a basis of $\mbox{Ran}(O_{t,k})$. For any $u\in{\mathbb{U}}^{1}(\mbox{Ran})$ and $k\in  \mathbb{T}_t$, there exit $\lambda_k^1, \cdots, \lambda_k^{r_k^1}\in{\mathbb{R}}$ such that $u_k-\mathbb{E}_tu_k=u_k=\sum_{i=1}^{r_k^1}\lambda_k^i(f_k^i)^T$. Then,
\begin{eqnarray*}
&&{F}_{t,k}^{(1)}u_k=\sum_{i=1}^{r_k^1}\lambda_k^i\left(
\begin{array}{c}
f_k^1\\\vdots\\f_k^{r_k^1}
\end{array}\right)(f_k^i)^T=\left(
\begin{array}{c}
\lambda_k^1\\\vdots\\\lambda_k^{r_k^1}
\end{array}\right)\triangleq \lambda_k.
\end{eqnarray*}
For ${k}\in  \mathbb{T}_t$, $u_k$ is $\mathbb{F}_k$\mbox{-measurable} and $\mathbb{E}|u_k|^2<\infty$, this implies that$\lambda_k$ is $\mathbb{F}_k$\mbox{-measurable} and $\mathbb{E}|\lambda_k|^2<\infty$. Therefore, $\widetilde{\mathbb{U}}^{1}(\mbox{Ran})=l^2_\mathbb{F}(t; \mathbb{R}^{r_t^1})\times\cdots\times l^2_\mathbb{F}(N-1; \mathbb{R}^{r_{N-1}^1})$.

Introduce a bounded linear operator $\tau$ from ${\mathbb{U}}^{1}(\mbox{Ran})$ to $\widetilde{\mathbb{U}}^{1}(\mbox{Ran})$:
\begin{eqnarray*}
({\tau}u)_k={F}_{t,k}^{(1)}u_k+{\Sigma}_{t,k}^{-1}{F}_{t,k}^{(1)}M_{t,k}(\alpha_k-\mathbb{E}_t\alpha_k),~~{k}\in  \mathbb{T}_t.
\end{eqnarray*}
We now prove that $\tau$ is a surjection. In fact, for any $\theta\in\widetilde{\mathbb{U}}^{1}(\mbox{Ran})$, we have $\mathbb{E}_t\theta_k=0, k\in \mathbb{T}_t$ and let
\begin{eqnarray*}
\left\{\begin{array}{l}
\bar{\alpha}_{k+1}=\big{(}A_{k}\bar{\alpha}_k+B^1_{k}({F}_{t,k}^{(1)})^T[\theta_k -{\Sigma}_{t,k}^{-1}{F}_{t,k}^{(1)}M_{t,k}(\bar{\alpha}_k-\mathbb{E}_t\bar{\alpha}_k)]
\big{)}\\[1mm]
\hphantom{\bar{\alpha}_{k+1}=}+\sum_{i=1}^p\big{(}C^i_{k}\bar{\alpha}_k+D^{1i}_{k} ({F}_{t,k}^{(1)})^T[\theta_k-{\Sigma}_{t,k}^{-1}{F}_{t,k}^{(1)}M_{t,k}(\bar{\alpha}_k
-\mathbb{E}_t\bar{\alpha}_k)]\big{)}w^i_k, \\[1mm]
\bar{\alpha}_t=0,~~k\in  \mathbb{T}_t,
\end{array}
\right.
\end{eqnarray*}
and
\begin{eqnarray}\label{u_k_1}
u_k=({F}_{t,k}^{(1)})^T[\theta_k-{\Sigma}_{t,k}^{-1}{F}_{t,k}^{(1)}M_{t,k}(\bar{\alpha}_k-\mathbb{E}_t\bar{\alpha}_k)],~~k\in  \mathbb{T}_t.
\end{eqnarray}
Note that $u$ in (\ref{u_k_1}) is in ${\mathbb{U}}^{1}(\mbox{Ran})$. As ${F}_{t,k}^{(1)}({F}_{t,k}^{(1)})^T=I_{r_k^1}$, from (\ref{u_k_1}) we have
\begin{eqnarray*}
\theta_k=({\tau}u)_k,~~k\in  \mathbb{T}_t.
\end{eqnarray*}
Hence, $\tau$ is a surjection defined from ${\mathbb{U}}^{1}(\mbox{Ran})$ to $\widetilde{\mathbb{U}}^{1}(\mbox{Ran})$. From this, (\ref{convex-u-1-1}) and the procedure of contradiction, we have $\Sigma_{t,k}\succ0,k\in  \mathbb{T}_t$. This further implies $O_{t,k} \succeq0, k\in  \mathbb{T}_t$.

Let
\begin{eqnarray*}
{\mathbb{U}}^{2}(\mbox{Ran})=\left\{u \Bigg{|}\begin{array}
{l}u\in l^2_{\mathbb{F}}(\mathbb{T}_t; \mathbb{R}^{m_1}), \mathbb{E}_tu_k\in \mbox{Ran}(\mathcal{O}_{t,k}) \mbox{ and } \\
u_k-\mathbb{E}_tu_k=-({F}_{t,k}^{(1)})^T{\Sigma}_{t,k}^{-1}{F}_{t,k}^{(1)}M_{t,k}(\alpha_k-\mathbb{E}_t\alpha_k), k\in\mathbb{T}_t
\end{array} \right\}.
\end{eqnarray*}
Note further that $u\mapsto \widetilde{J}_1(t,0; u)$ is convex. If $u\in{\mathbb{U}}^{2}(\mbox{Ran})$, from (\ref{convex-u-1}) we have
\begin{eqnarray}\label{convex-u-2-1}
\hspace{-1em}\widetilde{J}_1(t,0;u)=\sum_{k=t}^{N-1}\mathbb{E}_t\Big{\{}\Big{[}\mathcal{F}_{t,k}^{(1)}\mathbb{E}_tu_k +\Gamma_{t,k}^{-1}\mathcal{F}_{t,k}^{(1)}\mathcal{M}_{t,k}\mathbb{E}_t\alpha_k\Big{]}^T
\Gamma_{t,k}\Big{[}\mathcal{F}_{t,k}^{(1)}\mathbb{E}_tu_k+\Gamma_{t,k}^{-1}\mathcal{F}_{t,k}^{(1)}\mathcal{M}_{t,k}\mathbb{E}_t\alpha_k\Big{]} \Big{\}}\geq0.
\end{eqnarray}
Introduce a set
\begin{eqnarray*}
\widetilde{\mathbb{U}}^{2}(\mbox{Ran})=\Big{\{}(\mathcal{F}_{t,t}^{(1)}\mathbb{E}_tu_t,\cdots,\mathcal{F}_{t,N-1}^{(1)} \mathbb{E}_tu_{N-1})\big{|} u\in{\mathbb{U}}^{2}(\mbox{Ran})\Big{\}}.
\end{eqnarray*}
For $k\in  \mathbb{T}_t$, let $\nu_k^1,\cdots,\nu_k^{r_k^2}$ denote the lines of $\mathcal{F}_{t,k}^{(1)}$. For any $u\in{\mathbb{U}}^{2}(\mbox{Ran})$ and $k\in  \mathbb{T}_t$, there exit $\chi_k^1,\cdots,\chi_k^{r_k^2}$ $\in{\mathbb{R}}$ such that $\mathbb{E}_tu_k=\sum_{i=1}^{r_k^2}\chi_k^i(\nu_k^i)^T$. Then,

\begin{eqnarray*}
&&\mathcal{F}_{t,k}^{(1)}\mathbb{E}_tu_k=\sum_{i=1}^{r_k^2}\chi_k^i\left(
\begin{array}{c}
\nu_k^1\\\vdots\\\nu_k^{r_k^2}
\end{array}\right)(\nu_k^i)^T=\left(
\begin{array}{c}
\chi_k^1\\\vdots\\\chi_k^{r_k^2}
\end{array}\right)\triangleq \chi_k
\end{eqnarray*}
Therefore, $\widetilde{\mathbb{U}}^{2}(\mbox{Ran})=l^2(t; \mathbb{R}^{r_t^2})\times\cdots\times l^2(N-1; \mathbb{R}^{r_{N-1}^2})$ is a deterministic space.
Furthermore, introduce a bounded linear operator $\phi$ from ${\mathbb{U}}^{2}(\mbox{Ran})$ to $\widetilde{\mathbb{U}}^{2}(\mbox{Ran})$:
\begin{eqnarray*}
({\phi}u)_k=\mathcal{F}_{t,k}^{(1)}\mathbb{E}_tu_k+\Gamma_{t,k}^{-1}\mathcal{F}_{t,k}^{(1)}
\mathcal{M}_{t,k}\mathbb{E}_t{\alpha}_k,~~{k}\in  \mathbb{T}_t.
\end{eqnarray*}
We now prove that $\phi$ is a surjection. In fact, for any $\varsigma\in\widetilde{\mathbb{U}}^{2}(\mbox{Ran})$, let
\begin{eqnarray*}
\left\{\begin{array}{l}
\widetilde{\alpha}_{k+1}=\Big{\{}A_{k}\widetilde{\alpha}_k+B^1_{k}\Big{[}(\mathcal{F}_{t,k}^{(1)})^T(\varsigma_k-\Gamma_{t,k}^{-1}\mathcal{F}_{t,k}^{(1)}\mathcal{M}_{t,k} \mathbb{E}_t\widetilde{\alpha}_k)\\[2mm]
\hphantom{\widetilde{\alpha}_{k+1}=}-({F}_{t,k}^{(1)})^T{\Sigma}_{t,k}^{-1}{F}_{t,k}^{(1)}M_{t,k}(\alpha_k-\mathbb{E}_t\alpha_k) \Big{]}\Big{\}}\\[2mm]
\hphantom{\widetilde{\alpha}_{k+1}=}+\sum_{i=1}^p\Big{\{}C^i_{k}\widetilde{\alpha}_k +D^{1i}_{k}\Big{[}(\mathcal{F}_{t,k}^{(1)})^T(\varsigma_k-\Gamma_{t,k}^{-1}\mathcal{F}_{t,k}^{(1)}\mathcal{M}_{t,k} \mathbb{E}_t\widetilde{\alpha}_k)\\[2mm]
\hphantom{\widetilde{\alpha}_{k+1}=}-({F}_{t,k}^{(1)})^T{\Sigma}_{t,k}^{-1}{F}_{t,k}^{(1)}M_{t,k}(\alpha_k-\mathbb{E}_t\alpha_k) \Big{]}\Big{\}}w^i_k, \\[2mm]
\widetilde{\alpha}_t=0,~~k\in  \mathbb{T}_t,
\end{array}
\right.
\end{eqnarray*}
and
\begin{eqnarray}\label{u_k_2}
u_k=(\mathcal{F}_{t,k}^{(1)})^T[\varsigma_k-\Gamma_{t,k}^{-1}\mathcal{F}_{t,k}^{(1)}\mathcal{M}_{t,k} \mathbb{E}_t\widetilde{\alpha}_k]-({F}_{t,k}^{(1)})^T{\Sigma}_{t,k}^{-1}{F}_{t,k}^{(1)}M_{t,k}(\alpha_k-\mathbb{E}_t\alpha_k),~~k\in  \mathbb{T}_t.
\end{eqnarray}
Note that $u$ in (\ref{u_k_2}) is in ${\mathbb{U}}^{2}(\mbox{Ran})$. As $\mathcal{F}_{t,k}^{(1)}(\mathcal{F}_{t,k}^{(1)})^T=I_{r_k^2}$, from (\ref{u_k_2}) we have
\begin{eqnarray*}
\varsigma_k=({\phi}u)_k,~~k\in  \mathbb{T}_t.
\end{eqnarray*}
Hence $\phi$ is a surjection defined from ${\mathbb{U}}^{2}(\mbox{Ran})$ to $\widetilde{\mathbb{U}}^{2}(\mbox{Ran})$. From this, (\ref{convex-u-2-1}) and the procedure of contradiction, we have $\Gamma_{t,k}\succ0, k\in  \mathbb{T}_t$. This further implies $\mathcal{O}_{t,k}\succeq0,k\in  \mathbb{T}_t$. Furthermore, from (\ref{convex-v}), it is easy to get $\mathbb{O}_{k,k}\succeq0,k\in  \mathbb{T}_t$.

We now prove c) of Theorem \ref{Theorem-Problem-GLQ}. Note that
\begin{eqnarray*}
\eta^u_k=u_k-({F}_{t,k}^{(1)})^T{\Sigma}_{t,k}^{-1}{F}_{t,k}^{(1)}M_{t,k}^{1(1)}({\alpha}^u_k-\mathbb{E}_t{\alpha}^u_k)
-(\mathcal{F}_{t,k}^{(1)})^T\Gamma_{t,k}^{-1}\mathcal{F}_{t,k}^{(1)}\mathcal{M}_{t,k}\mathbb{E}_t{\alpha}^u_k,~~k\in  \mathbb{T}_t.
\end{eqnarray*}
Then,
\begin{eqnarray}\label{convex-u-3}
&&\widetilde{J}_1(t,0;\eta^u)=\sum_{k=t}^{N-1}\mathbb{E}_t\Big{\{}\Big{[}{F}_{t,k}^{(1)} (u_k-\mathbb{E}_tu_k)\Big{]}^T\Sigma_{t,k}\Big{[}{F}_{t,k}^{(1)}(u_k-\mathbb{E}_tu_k)\Big{]}
+\Big{[}\mathcal{F}_{t,k}^{(1)}\mathbb{E}_tu_k\Big{]}^T\Gamma_{t,k}\Big{[}\mathcal{F}_{t,k}^{(1)}\mathbb{E}_tu_k\Big{]}\Big{\}}\nonumber\\
&&\hphantom{\widetilde{J}_1(t,0;v^u)=}+2\sum_{k=t}^{N-1}\mathbb{E}_t\Big{[}\Big{(}{F}_{t,k}^{(2)}M_{t,k} ({\alpha}_k^u-\mathbb{E}_t{\alpha}_k^u)\Big{)}^T{F}_{t,k}^{(2)}
(u_k-\mathbb{E}_tu_k)\Big{]}\nonumber\\
&&\hphantom{\widetilde{J}_1(t,0;v^u)=}
+2\sum_{k=t}^{N-1}\mathbb{E}_t\Big{[}\Big{(}\mathcal{F}_{t,k}^{(2)}\mathcal{M}_{t,k} \mathbb{E}_t{\alpha}_k^u\Big{)}^T\mathcal{F}_{t,k}^{(2)}\mathbb{E}_tu_k\Big{]}\geq0.
\end{eqnarray}
In the above, we must have
\begin{eqnarray}\label{condition_4}
\left(
\begin{array}{c}
({F}_{t,k}^{(2)})^T{F}_{t,k}^{(2)}M_{t,k} ({\alpha}_{k}^u-\mathbb{E}_t{\alpha}_{k}^u)\\[1mm]
({\mathcal{F}}_{t,k}^{(2)})^T\mathcal{F}_{t,k}^{(2)}\mathcal{M}_{t,k} \mathbb{E}_t{\alpha}_{k}^u
\end{array}
\right)=0,~~a.s.,~~k\in  \mathbb{T}_t.
\end{eqnarray}
Otherwise, assume there exist $k_1\in \mathbb{T}_t$ and $\widehat{u}$ such that
\begin{eqnarray*}
\left(
\begin{array}{c}
a_1\\a_2
\end{array}
\right)=
\left(
\begin{array}{c}
a_1\\a_2
\end{array}
\right)(\omega)
\equiv
\left(
\begin{array}{c}
({F}_{t,k_1}^{(2)})^T{F}_{t,k_1}^{(2)}M_{t,k_1}({\alpha}_{k_1}^{\widehat{u}}-\mathbb{E}_t{\alpha}_{k_1}^{\widehat{u}})\\[1mm]
({\mathcal{F}}_{t,k_1}^{(2)})^T\mathcal{F}_{t,k_1}^{(2)}\mathcal{M}_{t,k_1}\mathbb{E}_t{\alpha}_{k_1}^{\widehat{u}}
\end{array}
\right)(w)\neq 0,~~~\mbox{for }\omega\in \Lambda_1
\end{eqnarray*}
with $\Lambda_1\in \mathbb{F}_{k_1}$ and its probability $\mathbb{P}(\Lambda_1)>0$. If  $\mathbb{E}_t\big{(}|a_1|^2+|a_2|^2\big{)}(\omega)=0$, a.s., we must have
\begin{eqnarray*}
0=\int_{\Omega}\mathbb{E}_t\big{(}|a_1|^2+|a_2|^2\big{)}\mathbb{P}(d\omega)
=\int_{\Omega}\big{(}|a_1|^2+|a_2|^2\big{)}\mathbb{P}(d\omega)\geq \int_{\Lambda_1}\big{(}|a_1|^2+|a_2|^2\big{)}\mathbb{P}(d\omega)>0,
\end{eqnarray*}
which is impossible. Hence, there exists $\Lambda\in \mathbb{F}_t$ with $\mathbb{P}(\Lambda)>0$ such that $\mathbb{E}_t\big{(}|a_1|^2+|a_2|^2\big{)}(\omega)>0$, for $\omega\in \Lambda$. Introduce a new control
\begin{eqnarray}\label{control_u}
u_k=\left\{
\begin{array}{ll}
\widehat{u}_k,&~~~k=t,...,k_1-1,\\
\widehat{u}_{k_1}+c_1a_1+c_2a_2,&~~~k=k_1,\\
0,&~~~k=k_1+1,...,N-1,
\end{array}
\right.
\end{eqnarray}
where
\begin{eqnarray*}
&&c_1=c_1(\omega)=
\left\{
\begin{array}{ll}
\left\{
\begin{array}{l}
-\frac{\ds b_1+b_3
}
{2\mathbb{E}_t\left|a_1\right|^2}(\omega)\cdot I_{\big{\{}\mathbb{E}_t\left|a_1\right|^2>0,\, \mathbb{E}_t\left|a_2\right|^2>0\big{\}}}(w)\\[1.5mm]
+0\cdot I_{\big{\{}\mathbb{E}_t\left|a_1\right|^2=0,\, \mathbb{E}_t\left|a_2\right|^2>0\big{\}}}(w)\\[1.5mm]
-\frac{\ds 1+b_1+b_2+b_3+b_4}
{2\mathbb{E}_t\left|a_1\right|^2}(\omega)\cdot I_{\big{\{}\mathbb{E}_t\left|a_1\right|^2>0, \, \mathbb{E}_t\left|a_2\right|^2=0\big{\}}}(\omega)\end{array}
\right.,&~~~\omega\in \Lambda,\\[1.5mm]
0,&~~~\omega\in \Omega- \Lambda,
\end{array}
\right.\\[2mm]
&&c_2=c_2(\omega)=\left\{
\begin{array}{ll}
\left\{
\begin{array}{l}
-\frac{\ds 1+b_2+b_4}
{2\mathbb{E}_t\left|a_2\right|^2}(\omega)\cdot I_{\big{\{}\mathbb{E}_t\left|a_1\right|^2>0, \, \mathbb{E}_t\left|a_2\right|^2>0\big{\}}}(\omega)\\[1.5mm]
-\frac{\ds 1+b_1+b_2+b_3+b_4}
{2\mathbb{E}_t\left|a_2\right|^2}(\omega)\cdot I_{\big{\{}\mathbb{E}_t\left|a_1\right|^2=0, \, \mathbb{E}_t\left|a_2\right|^2>0\big{\}}}(\omega)\\[1.5mm]
+0\cdot I_{\big{\{}\mathbb{E}_t\left|a_1\right|^2>0,\, \mathbb{E}_t\left|a_2\right|^2=0\big{\}}}(w)
\end{array}
\right.,&~~~\omega\in \Lambda,\\[1.5mm]
0,&~~~\omega\in \Omega- \Lambda
\end{array}
\right.
\end{eqnarray*}
with $I_{\{\cdot\}}(\omega)$ be the indicator function. Then, under (\ref{control_u}) and for  $\omega\in \Lambda$, we have
\begin{eqnarray}\label{convex-u-4}
&&\hspace{-2em}\widetilde{J}_1(t,0;\eta^u)(\omega)=\big{(}b_1+b_2+b_3+b_4+2c_1\mathbb{E}_t\left|a_1\right|^2
+2c_2\mathbb{E}_t\left|a_2\right|^2\big{)}(\omega)\nonumber\\
&&\hspace{-2em}\hphantom{\widetilde{J}_1(t,0;\eta^u)(\omega)}=-1,
\end{eqnarray}
where
\begin{eqnarray*}
%
&&b_1=\sum_{k=t}^{k_1}\mathbb{E}_t\Big{\{}\Big{[}{F}_{t,k}^{(1)} (\widehat{u}_k-\mathbb{E}_t\widehat{u}_k)\Big{]}^T{\Sigma}_{t,k}
\Big{[}{F}_{t,k}^{(1)}(\widehat{u}_k-\mathbb{E}_t\widehat{u}_k)\Big{]}\Big{\}},\\
&&
b_2=\sum_{k=t}^{k_1}\mathbb{E}_t\Big{\{}\Big{[}\mathcal{F}_{t,k}^{(1)}\mathbb{E}_t\widehat{u}_k\Big{]}^T \Gamma_{t,k}\Big{[}\mathcal{F}_{t,k}^{(1)}\mathbb{E}_t\widehat{u}_k\Big{]}
\Big{\}},\\
&&
b_3=2\sum_{k=t}^{k_1}\mathbb{E}_t\Big{[}\Big{(}{F}_{t,k}^{(2)}M_{t,k} ({\alpha}_k^{\widehat{u}}-\mathbb{E}_t{\alpha}_k^{\widehat{u}})\Big{)}^T{F}_{t,k}^{(2)}
(\widehat{u}_k-\mathbb{E}_t\widehat{u}_k)\Big{]},\\
&&
b_4=2\sum_{k=t}^{k_1}\mathbb{E}_t\Big{[}\Big{(}\mathcal{F}_{t,k}^{(2)}\mathcal{M}_{t,k} \mathbb{E}_t{\alpha}_k^{\widehat{u}}\Big{)}^T\mathcal{F}_{t,k}^{(2)}\mathbb{E}_t\widehat{u}_k\Big{]}.
\end{eqnarray*}
As $\mathbb{P}(\Lambda)>0$, (\ref{convex-u-4}) contradicts the convex condition (\ref{convex-u-3}). Hence, we have (\ref{condition_4}), and (\ref{condition-1}) (\ref{condition-3}) follow.

\emph{ii)$\Rightarrow$i)}. From (\ref{convex-u-3}) and  b), c) of Theorem \ref{Theorem-Problem-GLQ}, we have  $\widetilde{J}_1(t,0;\eta^u)\geq0$ and $\widetilde{J}_2(k, 0; v_k)\geq 0$ for any $u\in l^2_{\mathbb{F}}(\mathbb{T}_t; \mathbb{R}^{m_1}), v_k\in l^2_{\mathbb{F}}(k; \mathbb{R}^{m_2}), k\in \mathbb{T}_t$. Then, we need only to show
\begin{eqnarray}\label{condition_2}
\Big{\{}\eta^u\mid u\in l^2_{\mathbb{F}}(\mathbb{T}_t; \mathbb{R}^{m_1})\Big{\}}=l^2_{\mathbb{F}}(\mathbb{T}_t; \mathbb{R}^{m_1}).
\end{eqnarray}
In fact, for any $\widetilde{\eta}\in l^2_{\mathbb{F}}(\mathbb{T}_t; \mathbb{R}^{m_1})$, let
\begin{eqnarray*}
u_k=\widetilde{\eta}_k+({F}_{t,k}^{(1)})^T{\Sigma}_{t,k}^{-1}{F}_{t,k}^{(1)}M_{t,k}^{1(1)}(\widetilde{\alpha}_k-\mathbb{E}_t\widetilde{\alpha}_k)
+(\mathcal{F}_{t,k}^{(1)})^T\Gamma_{t,k}^{-1}\mathcal{F}_{t,k}^{(1)}\mathcal{M}_{t,k}^{1(1)}\mathbb{E}_t\widetilde{\alpha}_k,~~k\in  \mathbb{T}_t,
\end{eqnarray*}
where
\begin{eqnarray*}
\left\{\begin{array}{l}
\widetilde{\alpha}_{k+1}=\big{(}A_{k}\widetilde{\alpha}_k+B^1_{k}\widetilde{\eta}_k\big{)}
+\sum_{i=1}^p\big{(}C^i_{k}\widetilde{\alpha}_k+D^{1i}_{k}\widetilde{\eta}_k\big{)}w^i_k, \\[1mm]
\widetilde{\alpha}_t=0,~~k\in  \mathbb{T}_t.
\end{array}
\right.
\end{eqnarray*}
Thus $\eta^u=\widetilde{\eta}$. Hence, (\ref{condition_2}) holds, which together with $\widetilde{J}_1(t,0; \eta^u)\geq0$ implies
\begin{eqnarray*}
\inf_{u\in l^2_\mathbb{F}(\mathbb{T}_t; \mathbb{R}^m)} {\widetilde{J}_1(t,0;u)}=\inf_{u\in l^2_\mathbb{F}(\mathbb{T}_t; \mathbb{R}^m)} {\widetilde{J}_1(t,0; \eta^u)}\geq0.
\end{eqnarray*}
This completes the proof. \hfill $\square$

\textbf{Proof of Theorem \ref{Theorem-Problem-GLQ}}. This follows from Theorem \ref{Nece-suffic1}, Proposition \ref{Thm-stationary-condition} and Proposition \ref{Thm-convexity-condition}. \hfill $\square$

\textbf{Proof of Theorem \ref{Theorem-Problem-GLQ-3}}. Following (\ref{Propo-4-1-1})-(\ref{Propo-4-1-4}) and by deduction, we have that Problem (GLQ) admits a unique open-loop equilibrium, as $\mathbf{W}_{t,k}, \widetilde{\mathbf{W}}_{t,k}$ are nonsigular $k\in \mathbb{T}_t$. The expression of open-loop equilibrium follows easily form Theorem \ref{Theorem-Problem-GLQ}.  This completes the proof. \hfill $\square$

\subsection{Proof of Theorem \ref{Theorem-Problem-mv}}

%

%
%
%
Applying the general theory of Section \ref{Section-results} to Problem (MV), (\ref{P})-(\ref{h}) (with $t=0$) becomes to
\begin{eqnarray*}\label{P-mv}
\left\{
\begin{array}{l}
\overline{P}_{k}=A_k^T\overline{P}_{k+1}A_k-\Big{[}(\overline{H}^{1(1)}_{k})^T~\,(\overline{H}^{1(2)}_{k})^T\Big{]} \overline{\mathbf{W}}_{k}^\dagger\left[
\begin{array}{c}
\overline{H}^{1(1)}_{k}\\ \widehat{\overline{\mathcal H}}^{2(2)}_{k}
\end{array}
\right],\\[3mm]
\overline{\mathcal{P}}_{k}=A_{k}^T\overline{\mathcal{P}}_{k+1}A_{k}-\Big{[}(\overline{\mathcal H}^{1(1)}_{k})^T~\,(\overline{\mathcal H}^{1(2)}_{k})^T\Big{]}
\widetilde{\overline{\mathbf{W}}}_{k}^\dagger
\left[
\begin{array}{c}
\overline{\mathcal H}^{1(1)}_{k}\\\overline{\mathcal H}^{2(2)}_{k}
\end{array}
\right], \\[2mm]
\overline{\sigma}_{k}=-\Big{[}(\overline{\mathcal H}^{1(1)}_{k})^T~\, (\overline{\mathcal H}_{k}^{1(2)})^T\Big{]}\widetilde{\overline{\mathbf{W}}}_{k}^\dagger
\left[
\begin{array}{c}
\overline{h}_{k}^1\\h_{k}^2
\end{array}
\right]
+A_{k}^T\overline{\sigma}_{k+1},\\[2mm]
\overline{P}_{N}=G^1,  \overline{\mathcal{P}}_{N}=0,  \overline{\sigma}_{N}=g^1,   ~~~~k\in  \mathbb{T},
\end{array}
\right.
\end{eqnarray*}
\begin{eqnarray*}\label{T-mv}
\left\{
\begin{array}{l}
\overline{T}_{k}=(A_{k})^T\overline{T}_{k+1}A_{k}-\Big{[}(\overline{H}_{k}^{2(1)})^T~\,(\overline{H}_{k}^{2(2)})^T\Big{]} \overline{\mathbf{W}}_{k}^\dagger
\left[
\begin{array}{c}
\overline{H}^{1(1)}_{k}\\ \widehat{\overline{\mathcal H}}^{2(2)}_{k}
\end{array}
\right],\\[2mm]
\overline{\mathcal{T}}_{k}=(A_{k})^T\overline{\mathcal{T}}_{k+1}A_{k}-\Big{[}(\widehat{\overline{\mathcal H}}_{k}^{2(1)})^T~\,(\widehat{\overline{\mathcal H}}_{k}^{2(2)})^T\Big{]}{\overline{\mathbf{W}}}_{k}^\dagger
\left[
\begin{array}{c}
\overline{H}^{1(1)}_{k}\\ \widehat{\overline{\mathcal H}}^{2(2)}_{k}
\end{array}
\right],\\[2mm]
\widetilde{\overline{T}}_{k}=(A_k)^T\widetilde{\overline{T}}_{k+1}A_k -\Big{[}({\overline{\mathcal H}}_{k}^{2(1)})^T~\,({\overline{\mathcal H}}_{k}^{2(2)})^T\Big{]}\widetilde{\overline{\mathbf{W}}}_{k}^\dagger\left[
\begin{array}{c}
\overline{\mathcal H}^{1(1)}_{k}\\ \overline{\mathcal H}^{2(2)}_{k}
\end{array}
\right]\\[2mm]
\hphantom{\mathcal{T}_{k}=}+\Big{[}(\widehat{\overline{\mathcal H}}_{k}^{2(1)})^T~\,(\widehat{\overline{\mathcal H}}_{k}^{2(2)})^T\Big{]}{\overline{\mathbf{W}}}_{k}^\dagger \left[
\begin{array}{c}
\overline{H}^{1(1)}_{k}\\ \widehat{\overline{\mathcal H}}^{2(2)}_{k}
\end{array}
\right],\\[2mm]
%
%
%
\overline{\xi}_{k}=-
\Big{[}({\overline{\mathcal H}}_{k}^{2(1)})^T~\,({\overline{\mathcal H}}_{k}^{2(2)})^T\Big{]}\widetilde{\overline{\mathbf{W}}}_{k}^\dagger
\left[
\begin{array}{c}
\overline{h}_{k}^1\\\overline{h}_{k}^2
\end{array}
\right]+A_{k}^T\overline{\xi}_{k+1},\\[2mm]
\overline{T}_{k,N}=G^2,~~  \overline{\mathcal{T}}_{k,N}=0,~~
\widetilde{\overline{T}}_{k,N}=0,~~  \overline{\xi}_{N}=g^2, ~~~k\in \mathbb{T},
\end{array}
\right.
\end{eqnarray*}
where
\begin{eqnarray*}\label{W-H-mv}
\left\{
\begin{array}{l}
\overline{\mathcal W}^{1(1s)}_{k}=\Upsilon_{k}^{(1s)}+(B_{k}^{1})^T\overline{\mathcal{P}}_{k+1}B_{k}^s+\sum_{i,j=1}^{p_0}\overline{\delta}_{k}^{ij} (D^{1 i}_{k})^T\overline{P}_{k+1}D_{k}^{sj},\\[2mm]
\overline{\mathcal W}^{2(2s)}_{k}=\Upsilon_{k}^{(2s)}+(B_{k}^2)^T(\overline{\mathcal{T}}_{k+1}+\widetilde{\overline{T}}_{k+1})B_{k}^s +\sum_{i,j=1}^{p_0}\overline{\delta}_{k}^{ij} (D^{2 i}_{k})^T\overline{T}_{k+1}D_{k}^{sj},\\[2mm]
\overline{W}^{1(1s)}_{k}=\Upsilon_{k}^{(1s)}+(B_{k}^1)^T\overline{P}_{k+1}B_{k}^s+\sum_{i,j=1}^{p_0}\overline{\delta}_{k}^{ij} (D^{1 i}_{k})^T\overline{P}_{k+1}D_{k}^{sj},\\[2mm]
\widehat{\overline{\mathcal{W}}}^{2(2s)}_{k}=\Upsilon_{k}^{(2s)}+(B_{k}^2)^T\overline{\mathcal{T}}_{k+1}B_{k}^s +\sum_{i,j=1}^{p_0}\overline{\delta}_{k}^{ij} (D^{2 i}_{k})^T\overline{T}_{k+1}D_{k}^{sj},\\[2mm]
\overline{\mathcal{H}}_{k}^{1(s)}=(B_{k}^s)^T\overline{\mathcal{P}}_{k+1}A_{k},~~~~
\overline{\mathcal{H}}_{k}^{2(s)}=(B_{k}^s)^T (\overline{\mathcal{T}}_{k+1}+\widetilde{\overline{T}}_{k+1})A_{k},\\[2mm]
\widehat{\overline{\mathcal{H}}}_{k}^{2(s)}=(B_{k}^s)^T \overline{\mathcal{T}}_{k+1}A_{k},~~~~
\overline{H}_{k}^{1(s)}=(B_{k}^s)^T \overline{P}_{k+1}A_{k},\\[2mm]
\overline{H}_{k}^{2(s)}=(B_{k}^s)^T \overline{T}_{k+1}A_{k},
~~~k\in  \mathbb{T},~~~s=1,2,
\end{array}
\right.
\end{eqnarray*}
and
\begin{eqnarray*}\label{W-bf-mv}
&&\hspace{-4em}\overline{\mathbf{W}}_{k}=\left(
\begin{array}{ll}
\overline{W}^{1(11)}_{k}&\overline{W}^{1(12)}_{k}\\\widehat{\overline{\mathcal W}}^{2(21)}_{k}&\widehat{\overline{\mathcal W}}^{2(22)}_{k}
\end{array}\right),~~~ \widetilde{\overline{\mathbf{W}}}_{k}=\left(
\begin{array}{ll}
\overline{\mathcal W}^{1(11)}_{k}&\overline{\mathcal W}^{1(12)}_{k}\\ \overline{\mathcal W}^{2(21)}_{k}&\overline{\mathcal W}^{2(22)}_{k}
\end{array}\right),\\
&&\hspace{-4em}
\overline{h}_{k}^1=(B_{k}^1)^T\overline{\sigma}_{k+1}, ~~\overline{h}_{k}^2=(B_{k}^2)^T\overline{\xi}_{k+1},~~~~~k\in \mathbb{T}.
\end{eqnarray*}
with $\Upsilon^{(11)}_k=\Phi_k, \Upsilon^{(12)}_k=-\Phi_k, \Upsilon_k^{(21)}=-\Phi_k, \Upsilon_k^{(22)}=\Phi_k$. 

Noting $\overline{\mathcal{P}}_N=\overline{\mathcal{T}}_N=\widetilde{\overline{T}}_N=0$, simple calculations show
\begin{eqnarray*}
\overline{\mathcal{P}}_k=\overline{\mathcal{T}}_k=\widetilde{\overline{T}}_k=0,~~\overline{\mathcal{H}}_{k}^{1(s)} =\overline{\mathcal{H}}_{k}^{2(s)}=\widehat{\overline{\mathcal{H}}}_{k}^{2(s)}=0, ~~k\in \mathbb{T},
\end{eqnarray*}
and
\begin{eqnarray*}
\overline{P}_{N-1}=\left(
\begin{array}{cc}
s_{N-1}^2&0\\0&0
\end{array}
\right)-s_{N-1}^2\left(
\begin{array}{cc}
(\mathbb{E}\Theta_{N-1})^T&0\\0&0
\end{array}
\right)\overline{\mathbf{W}}_{N-1}^\dagger \left(
\begin{array}{cc}
\mathbb{E}\Theta_{N-1}&0\\0&0
\end{array}
\right)=\left(
\begin{array}{cc}
\overline{P}_{N-1}^{(11)}&0\\0&0
\end{array}
\right)
\end{eqnarray*}
for some $\overline{P}_{N-1}^{(11)}$, which implies the form of $\overline{P}_k$:
\begin{eqnarray*}
\overline{P}_k=\left(
\begin{array}{cc}
\overline{P}_{k}^{(11)}&0\\0&0
\end{array}
\right),~~~k\in \mathbb{T}.
\end{eqnarray*}
Therefore,
\begin{eqnarray}\label{barP}
\overline{P}_{k}^{(11)}=s_k^2\overline{P}_{k+1}^{(11)}\Big{[}1-\overline{P}_{k+1}^{(11)}(\mathbb{E}\Theta_k)^T \big{(}\overline{\textbf{W}}_k^\dagger\big{)}^{(11)} \mathbb{E}\Theta_k\Big{]},~~~k\in \mathbb{T}.
\end{eqnarray}
Hence,
\begin{eqnarray*}
&&\overline{\mathbf{W}}_k=\Upsilon_k+\left(
\begin{array}{cc}
\overline{P}_{k+1}^{(11)}\mathbb{E}\big{(}\Theta_k\Theta_k^T \big{)}&0\\[2mm]\overline{T}_{k+1}^{(21)}\mbox{Cov}(\Theta_k)&\overline{T}_{k+1}^{(22)}\mbox{Cov}(\Theta_k)
\end{array}
\right),\\[2mm]
&&\widetilde{\overline{\mathbf{W}}}_k=\Upsilon_k+\left(
\begin{array}{cc}
\overline{P}_{k+1}^{(11)}\mbox{Cov}(\Theta_k)&0\\[2mm]\overline{T}_{k+1}^{(21)}\mbox{Cov}(\Theta_k)& \overline{T}_{k+1}^{(22)}\mbox{Cov}(\Theta_k)
\end{array}
\right).
\end{eqnarray*}
Moreover,
\begin{eqnarray*}
&&\sigma_k=\left\{
\begin{array}{ll}
\left(
\begin{array}{c}
-\frac{\lambda }{2}\\0
\end{array}
\right),&~~~k=N,\\[4mm]
\left(
\begin{array}{c}
-\frac{\lambda}{2}s_{N-1}\\0
\end{array}
\right),&~~~k=N-1,\\[4mm]
\left(
\begin{array}{c}
-\frac{\lambda}{2}s_k\cdots s_{N-1}\\0
\end{array}
\right),&~~~k\in \{0,...,N-2\},
\end{array}
\right.\\
&&\xi_k=\left\{
\begin{array}{ll}
\left(
\begin{array}{c}
0\\-\frac{\lambda }{2}
\end{array}
\right),&~~~k=N,\\[4mm]
\left(
\begin{array}{c}
0\\
-\frac{\lambda}{2}s_{N-1}
\end{array}
\right),&~~~k=N-1,\\[4mm]
\left(
\begin{array}{c}
0\\-\frac{\lambda}{2}s_k\cdots s_{N-1}
\end{array}
\right),&~~~k\in \{0,...,N-2\}.
\end{array}
\right.
\end{eqnarray*}
Hence, we have (\ref{W-bf-mv-0}) and (\ref{P-mv-0-0}).  Furthermore, for Problem (MV) and under the parameters (\ref{system-parameter-mv-1})-(\ref{D-m}), (\ref{U})-(\ref{M-O}) becomes to
\begin{eqnarray*}\label{U-mv}
\left\{
\begin{array}{l}
\overline{U}_{k}=A_k^T\overline{U}_{k+1}A_k-\overline{M}_{k}^T\overline{O}_{k}^\dagger \overline{M}_{k},\\[2mm]
\overline{\mathcal{U}}_{k}=A_k^T\overline{\mathcal{U}}_{k+1}A_k -\overline{\mathcal{M}}_{k}^T\overline{\mathcal{O}}_{k}^\dagger \overline{\mathcal{M}}_{k},\\[2mm]
\overline{U}_{N}=G^1,~~\overline{\mathcal{U}}_{N}=0,~~~~k\in  \mathbb{T},
\end{array}
\right.
\end{eqnarray*}
and
\begin{eqnarray*}\label{V-mv}
\left\{
\begin{array}{l}
\overline{V}_{k}=A_k^T\overline{V}_{k+1}A_k,\\[1mm]
\overline{\mathcal{V}}_{k}=A_k^T\overline{\mathcal{V}}_{k+1}A_k\equiv 0,\\[1mm]
\overline{V}_{N}=G^2,~~\overline{\mathcal{V}}_{N}=0,~~~~\ell\in  \mathbb{T}
\end{array}
\right.
\end{eqnarray*}
with
\begin{eqnarray*}\label{M-O-mv}
\left\{
\begin{array}{l}
\overline{M}_{k}=(B_k^1)^T\overline{U}_{k+1}A_k,\\[2mm]
\overline{\mathcal{M}}_{k}=(B_k^1)^T\overline{\mathcal{U}}_{k+1}A_k,\\[2mm]
\overline{O}_{k}=\Upsilon^{(11)}_k+(B_k^1)^T\overline{U}_{k+1}B_k^1+\sum_{i,j=1}^{p_0}\overline{\delta}_k^{ij}(D_k^{1i})^T\overline{U}_{k+1}D_k^{1j},\\[2mm]
\overline{\mathcal{O}}_{t,k}=\Upsilon^{(11)}_k+(B_k^1)^T\overline{\mathcal{U}}_{k+1}B_k^1+\sum_{i,j=1}^{p_0}\overline{\delta}_k^{ij} (D_k^{1i})^T\overline{U}_{k+1}D_k^{1j},\\[2mm]
\overline{\mathbb{O}}_{k}=\Upsilon^{(22)}_k+\sum_{i,j=1}^{p_0}\overline{\delta}_k^{ij}(D_k^{2i})^T\overline{V}_{k+1}D_k^{2j}.
\end{array}
\right.
\end{eqnarray*}
Introduce a new optimal control problem with the system dynamics
\begin{eqnarray}\label{delta}
\left\{
\begin{array}{l}
\theta_{k+1}=A_k\theta_k+B_k^1\nu_k+\sum_{i=1}^{p_0}D_k^{1i}\nu_k{w}_k^i,\\[1mm]
\theta_{t}=\bar{\theta},~~~k\in \mathbb{T},
\end{array}
\right.
\end{eqnarray}
and the objective functional
\begin{eqnarray}\label{J-delta}
&&\hspace{-3em}{J}^\theta(t,\bar{\theta};\nu)=\sum_{k=t}^{N-1}\mathbb{E}\big{[}\nu_k^TL^{1}_{k}\nu_k\big{]} +\mathbb{E}[\nu_N^TG^1\nu_N]+(\mathbb{E}\nu_N)^T \bar{G}^1\mathbb{E}\nu_N
\end{eqnarray}
that is to be minimized within $l^2_{\mathbb{F}^m}(\mathbb{T}_t; \mathbb{R}^{p_0})$.  Here, the parameters in (\ref{delta})-(\ref{J-delta}) are from (\ref{system-parameter-mv-1})-(\ref{D-m}). Clearly, this is a special example of the static mean-field LQ optimal control problem that is considered in \cite{Ni-Zhang-Li}.  As $\Upsilon_k^{(11)}\succeq 0, G^1\succeq 0, G^1+\bar{G}^1\succeq 0, k\in \mathbb{T}$, we have
from Theorem 4.3 of \cite{Ni-Zhang-Li} that
\begin{eqnarray*}
\overline{{O}}_{k}\overline{{O}}_{k}^\dagger\overline{{M}}_{k}=\overline{{M}}_{k},~~ \overline{\mathcal{{O}}}_{k}{\overline{\mathcal{O}}}_{k}^\dagger{\overline{\mathcal{M}}}_{k}={\overline{\mathcal{M}}}_{k},~~
\overline{O}_k\succeq 0,~~ \overline{\mathcal{O}}_k\succeq 0,~~
k\in \mathbb{T}.
\end{eqnarray*}
As  $\Upsilon^{(22)}_k, G^2\succeq 0$, we have $\overline{\mathbb{O}}_k\succeq 0, k\in \mathbb{T}$. This completes the proof by following Theorem \ref{Theorem-Problem-GLQ} and using the notations of (\ref{W-bf-mv-0})-(\ref{P-mv-0-0}). \hfill $\square$

\subsection{Proof of Theorem \ref{Theorem-Problem-mv2}}

\emph{i).}
Let $\overline{T}^{(ij)}_{k}$ be the $(i,j)$-th entry of $\overline{T}_{k}, i,j=1,2,~ k\in \mathbb{T}$, then
\begin{eqnarray*}
\left(
\begin{array}{cc}
\overline{T}^{(11)}_k&\overline{T}^{(12)}_k\\ \overline{T}^{(21)}_k&\overline{T}^{(22)}_k
\end{array}
\right)=s_k^2\left(
\begin{array}{cc}
\overline{T}^{(11)}_{k+1}&\overline{T}^{(12)}_{k+1}\\ \overline{T}^{(21)}_{k+1}&\overline{T}^{(22)}_{k+1}
\end{array}
\right)-s_k\Big{[}(\overline{H}_{k}^{2(1)})^T~\,(\overline{H}_{k}^{2(2)})^T\Big{]} \overline{\mathbf{W}}_{k}^{\dagger} \left(
\begin{array}{cc}
\overline{P}_{k+1}^{(11)}\mathbb{E}\Theta_{k}&0\\0&0
\end{array}
\right)
\end{eqnarray*}
with
\begin{eqnarray*}
\Big{[}(\overline{H}_{k}^{2(1)})^T~\,(\overline{H}_{k}^{2(2)})^T\Big{]}=s_k\left[
\begin{array}{cc}
\overline{T}_{k+1}^{(11)}(\mathbb{E}\Theta_k)^T&\overline{T}_{k+1}^{(21)}(\mathbb{E}\Theta_k)^T\\
\overline{T}_{k+1}^{(12)}(\mathbb{E}\Theta_k)^T&\overline{T}_{k+1}^{(22)}(\mathbb{E}\Theta_k)^T
\end{array}
\right].
\end{eqnarray*}
Therefore, for $k\in \mathbb{T}$,
\begin{eqnarray}\label{T-mv-1}
&&\hspace{-2em}\left(
\begin{array}{cc}
\overline{T}^{(11)}_k&\overline{T}^{(12)}_k\\ \overline{T}^{(21)}_k&\overline{T}^{(22)}_k
\end{array}
\right)=-s_k^2 \left(
\begin{array}{cc}
\overline{P}_{k+1}^{(11)}\overline{T}_{k+1}^{(11)}(\mathbb{E}\Theta_{k})^T\big{(} \overline{\mathbf{W}}^\dagger_k\big{)}^{(11)}\mathbb{E}\Theta_k+\overline{P}_{k+1}^{(11)}\overline{T}_{k+1}^{(21)} (\mathbb{E}\Theta_{k})^T\big{(}\overline{\mathbf{W}}^\dagger_k\big{)}^{(21)}\mathbb{E}\Theta_k&0\\ \overline{P}_{k+1}^{(11)}\overline{T}_{k+1}^{(12)}(\mathbb{E}\Theta_{k})^T\big{(} \overline{\mathbf{W}}^\dagger_k\big{)}^{(11)}\mathbb{E}\Theta_k+\overline{P}_{k+1}^{(11)}\overline{T}_{k+1}^{(22)} (\mathbb{E}\Theta_{k})^T\big{(}\overline{\mathbf{W}}^\dagger_k\big{)}^{(21)}\mathbb{E}\Theta_k&0
\end{array}
\right)\nonumber\\
&&\hspace{-2em}\hphantom{\left(
\begin{array}{cc}
\overline{T}^{(11)}_k&\overline{T}^{(12)}_k\\ \overline{T}^{(21)}_k&\overline{T}^{(22)}_k
\end{array}
\right)=}+s_k^2\left(
\begin{array}{cc}
\overline{T}^{(11)}_{k+1}&\overline{T}^{(12)}_{k+1}\\ \overline{T}^{(21)}_{k+1}&\overline{T}^{(22)}_{k+1}
\end{array}
\right)
\end{eqnarray}
with $\big{(}\overline{\mathbf{W}}_k^\dagger\big{)}^{(11)}, \big{(}\overline{\mathbf{W}}_k^\dagger\big{)}^{(21)}$ being the $(1,1)$-th and $(2,1)$-th blocks of $\overline{\mathbf{W}}_k^\dagger$. Therefore, (\ref{T-mv-0}) holds with the property $\overline{T}^{(22)}_{k}\geq 1, k\in {\mathbb{T}}$. Also, as $\overline{T}^{(12)}_N=0$, we have $\overline{T}^{(12)}_k=0, k\in \mathbb{T}$.

Note that $\mu_k=0, k\in \mathbb{T}$. In this case,
\begin{eqnarray*}
\overline{\mathbf{W}}_k=\left(
\begin{array}{cc}
\overline{P}_{k+1}^{(11)}\mathbb{E}\big{(}\Theta_k\Theta_k^T \big{)}&0\\[2mm]\overline{T}_{k+1}^{(21)}\mbox{Cov}(\Theta_k)&\overline{T}_{k+1}^{(22)}\mbox{Cov}(\Theta_k)
\end{array}
\right),~~~~k\in \mathbb{T}.
\end{eqnarray*}
Checking the definition of Moore-Penrose inverse, we have
\begin{eqnarray}\label{W-0-inv-2}
\overline{\mathbf{W}}^{\dagger}_{k}=\left(
\begin{array}{cc}
\big{[}\overline{P}_{k+1}^{(11)}\mathbb{E}\big{(}\Theta_{k}\Theta_{k}^T \big{)}\big{]}^{\dagger}&0\\[2mm]\big{(}\overline{\mathbf{W}}^{\dagger}_{k} \big{)}^{(21)}& \big{[}\overline{T}_{k+1}^{(22)}\mbox{Cov}(\Theta_{k})\big{]}^{\dagger}
\end{array}
\right),
\end{eqnarray}
where
\begin{eqnarray}\label{W-0-inv-(21)}
\big{(}\overline{\mathbf{W}}^{\dagger}_{k} \big{)}^{(21)}=-\big{[}\overline{T}_{k+1}^{(22)}\mbox{Cov}(\Theta_{k})\big{]}^{\dagger}\overline{T}_{k+1}^{(21)} \mbox{Cov}(\Theta_k) \big{[}\overline{P}_{k+1}^{(11)}\mathbb{E}\big{(}\Theta_{k}\Theta_{k}^T \big{)}\big{]}^{\dagger}.
\end{eqnarray}
From (\ref{T-mv-1}), (\ref{W-0-inv-(21)}), $\overline{T}_{k}^{(12)}=0$ and $\overline{T}_N^{(21)}=0$, it holds that $\overline{T}_k^{(21)}=0, k\in \mathbb{T}$. Hence, we have proved that $\overline{\mathbf{W}}_k$ has the following form:
\begin{eqnarray}\label{W-0}
\overline{\mathbf{W}}_k=\left(
\begin{array}{cc}
\overline{P}_{k+1}^{(11)}\mathbb{E}\big{(}\Theta_k\Theta_k^T \big{)}&0\\[2mm]0&\overline{T}_{k+1}^{(22)}\mbox{Cov}(\Theta_k)
\end{array}
\right),~~~~k\in \mathbb{T}.
\end{eqnarray}
From (\ref{barP}), we have
\begin{eqnarray}\label{matrix-2}
\overline{P}_{k}^{(11)}=s_{k}^2 \overline{P}_{k+1}^{(11)} \Big{[}1-(\mathbb{E}\Theta_{k})^T \big{[}\mathbb{E}\big{(}\Theta_k\Theta_k^T \big{)}\big{]}^\dagger \mathbb{E}\Theta_{k}\Big{]}.
\end{eqnarray}
By (50) of \cite{Ni-Zhang-Li}, we know $1-(\mathbb{E}\Theta_{k})^T \big{[}\mathbb{E}\big{(}\Theta_k\Theta_k^T \big{)}\big{]}^\dagger \mathbb{E}\Theta_{k}>0$. Therefore, $\overline{P}_{k}^{(11)}>0$ for $k\in \mathbb{T}$.

\emph{ii)}. In this case, $\widetilde{\overline{\mathbf{W}}}_k$ is of the following form
\begin{eqnarray}\label{W-0}
\widetilde{\overline{\mathbf{W}}}_k=\left(
\begin{array}{cc}
\overline{P}_{k+1}^{(11)}\mbox{Cov}(\Theta_k)&0\\[2mm]0&\overline{T}_{k+1}^{(22)}\mbox{Cov}(\Theta_k)
\end{array}
\right),~~~~k\in \mathbb{T}.
\end{eqnarray}
As $\mathbb{E}\Theta_k\in \mbox{Ran}\big{[}\mbox{Cov}(\Theta_k)\big{]}, k\in \mathbb{T}$, the condition (\ref{W-H-mv-2}) is satisfied. Hence, for any initial pair Problem (MV) admit open-loop self-coordination controls of precommitted self and of sophisticated selves, which coincide with the open-loop precommitted optimal control and open-loop time-consistent equilibrium control, respectively.

\emph{iii).} Let
$
\Xi_k^c=\Big{\{}\zeta \,\big{|}\, \mbox{Cov}(\Theta_k)\zeta =\mathbb{E}\Theta_k  \Big{\}}\neq \emptyset.
$
For given $\zeta_0\in \Xi^c_k$, there exists $d\geq 0$ such that
\begin{eqnarray}\label{L-0}
\Phi_k\zeta_0=d\mathbb{E}\Theta_k.
\end{eqnarray}
Note that (\ref{W-H-mv-2}) is equivalent to
\begin{eqnarray*}
\mbox{Ran}(\overline{\mathbf{H}}_{k})\subset \mbox{Ran}(\overline{\mathbf{W}}_{k}),~~~\overline{\mathbf{h}}_{k}\in \mbox{Ran}(\widetilde{\overline{\mathbf{W}}}_{k}),~~k\in \mathbb{T},
\end{eqnarray*}
and that
\begin{eqnarray*}
&&\mbox{Ran}(\overline{\mathbf{H}}_k)=\left\{
\left(
\begin{array}{c}
c \overline{P}_{k+1}^{(11)}\mathbb{E}\Theta_k\\[1.5mm]
0
\end{array}\right)\Bigg{|}\, c\in \mathbb{R}
\right\},\\[1mm]
&&\mbox{Ran}(\overline{\mathbf{W}}_k)=\left\{
\left(
\begin{array}{c}
L^{(1)}_ka+L_k^{(2)}b+\overline{P}_{k+1}^{(11)}\mathbb{E}\big{(}\Theta_k\Theta_k^T\big{)}a\\[1.5mm]
L_k^{(3)}a+L_k^{(4)}b+\mbox{Cov}(\Theta_k)\big{[}\overline{T}_{k+1}^{(21)}a+\overline{T}_{k+1}^{(22)}b \big{]}
\end{array}\right)\Bigg{|}\, a, b\in \mathbb{R}^m
\right\},\\[1mm]
&&\mbox{Ran}(\widetilde{\overline{\mathbf{W}}}_k)=\left\{
\left(
\begin{array}{c}
L^{(1)}_ka+L_k^{(2)}b+\overline{P}_{k+1}^{(11)}\mbox{Cov}(\Theta_k)a\\[1.5mm]
L_k^{(3)}a+L_k^{(4)}b+\mbox{Cov}(\Theta_k)\big{[}\overline{T}_{k+1}^{(21)}a+\overline{T}_{k+1}^{(22)}b \big{]}
\end{array}\right)\Bigg{|}\, a, b\in \mathbb{R}^m
\right\}.
\end{eqnarray*}
For $\zeta_0\in \Xi_k$, we have $\mbox{Cov}(\Theta_k)\zeta_0=\mathbb{E}\Theta_k$ and $L_k^{(i)}\zeta_0=d_i\mathbb{E}\Theta_k$ for some $d_i\in \mathbb{R}$. Letting $a=x_1\zeta_0, b=x_2\zeta_0$, the equation
\begin{eqnarray*}
\left(
\begin{array}{c}
L^{(1)}_ka+L_k^{(2)}b+\overline{P}_{k+1}^{(11)}\mathbb{E}\big{(}\Theta_k\Theta_k^T\big{)}a\\[1.5mm]
L_k^{(3)}a+L_k^{(4)}b+\mbox{Cov}(\Theta_k)\big{[}\overline{T}_{k+1}^{(21)}a+\overline{T}_{k+1}^{(22)}b \big{]}
\end{array}\right)=\left(
\begin{array}{c}
c \overline{P}_{k+1}^{(11)}\mathbb{E}\Theta_k\\[1.5mm]
0
\end{array}\right)
\end{eqnarray*}
becomes to finding $x_1,x_2$ such that
\begin{eqnarray}\label{linear-equation}
\left\{
\begin{array}{l}
d_1x_1+d_2x_2+\overline{P}_{k+1}^{(11)}(1+\mathbb{E}\Theta_k^T\zeta_0)x_1=c\overline{P}_{k+1}^{(11)},\\[1mm]
d_3x_1+d_4x_2+\overline{T}_{k+1}^{(21)}x_1+\overline{T}_{k+1}^{(22)}x_2=0
\end{array}
\right.
\end{eqnarray}
holds for given $c\in\mathbb{R}$. By some calculations, the determinant of coefficient matrix of (\ref{linear-equation}) is
\begin{eqnarray}\label{det-0}
D_{et}(d_1,d_2,d_3,d_4)=\Big{[}d_1+\overline{P}_{k+1}^{(11)}(1+\mathbb{E}\Theta_k^T\zeta_0)\Big{]}\Big{[}d_4 +\overline{T}_{k+1}^{(22)}\Big{]}- d_2\Big{[}d_3+\overline{T}_{k+1}^{(21)}\Big{]}.
\end{eqnarray}
If $c= 0$ or $\overline{P}_{k+1}^{(11)}=0$, $x_1$ and $x_2$ of (\ref{linear-equation}) can be both selected to be 0. For $c\neq 0$ and $\overline{P}_{k+1}^{(11)}\neq0$, we have following derivation. As $\overline{T}_{k+1}^{(22)}\neq 0$, let $d_1=-d_2=-d_3=d_4$ and (\ref{det-0}) becomes
\begin{eqnarray}\label{det}
D_{et}(d_1,-d_1,-d_1,d_1)=\Big{[}\overline{P}_{k+1}^{(11)}(1+\mathbb{E}\Theta_k^T\zeta_0)+\overline{T}_{k+1}^{(21)}+ \overline{T}_{k+1}^{(22)}\Big{]}d_1+\overline{P}_{k+1}^{(11)}(1+\mathbb{E}\Theta_k^T\zeta_0)\overline{T}_{k+1}^{(22)}.
\end{eqnarray}
As $1+\mathbb{E}\Theta_k^T\zeta_0=1+\zeta^T\mbox{Cov}(\Theta_k)\zeta_0\geq 1$, we must have $\overline{P}_{k+1}^{(11)}(1+\mathbb{E}\Theta_k^T\zeta_0)\overline{T}_{k+1}^{(22)}\neq 0$. Therefore, there exists some $d_1$ such that $D_{et}(d_1,-d_1,-d_1,d_1)\neq 0$ and (\ref{linear-equation}) is solvable.

Therefore, by selecting $\Upsilon_k, k\in \mathbb{T}_t$ with (\ref{L-0}) we can have $\mbox{Ran}(\overline{\mathbf{H}}_{k})\subset \mbox{Ran}(\overline{\mathbf{W}}_{k}), k\in \mathbb{T}$ and similarly $\overline{\mathbf{h}}_{k}\in \mbox{Ran}(\widetilde{\overline{\mathbf{W}}}_{k}), k\in \mathbb{T}$ can be proved.
This completes the proof. \hfill$\square$

\section{Conclusion}\label{section-conclusion}

For a time-inconsistent LQ optimal control, a Nash-type fictitious game framework is introduced, where the game is between the decision maker and an auxiliary control variable. The decision maker and auxiliary control variable are called real player and fictitious player, which look for time-consistent policy and precommitted optimal policy, respectively.
The equilibrium policy of real player is called an open-loop self-coordination control of the LQ problem.
As a generalization, a time-inconsistent nonzero-sum stochastic LQ dynamic game is studied, for which one player is to find the precommitted policy and the other player is to find the time-consistent policy. Necessary and sufficient conditions are derived to characterize the open-loop equilibrium of this nonzero-sum stochastic LQ dynamic game via Riccati-like equations, and as a byproduct, result to ensure the existence of open-loop self-coordination control of the original LQ optimal control is also obtained. To test the general theory, the mean-variance portfolio selection is investigated. For future research, the closed-loop self-coordination should be investigated.



\end{document}